\documentclass[10pt]{amsart}
\font\emailfont=cmtt9

\newcommand\Id{\mathrm{Id}}

\newcommand{\rk}{\mathrm{rk}}

\newtheorem{theorem}{Theorem}[section]
\newtheorem{prop}[theorem]{Proposition}
\newtheorem{cor}[theorem]{Corollary}

\newtheorem{lemma}[theorem]{Lemma}

\newtheorem{defn}[theorem]{Definition}

\def\endproof{\relax\ifmmode\expandafter\endproofmath\else
  \unskip\nobreak\hfil\penalty50\hskip.75em\hbox{}\nobreak\hfil\bull
  {\parfillskip=0pt \finalhyphendemerits=0 \bigbreak}\fi}
\def\endproofmath$${\eqno\bull$$\bigbreak}
\def\bull{\vbox{\hrule\hbox{\vrule\kern3pt\vbox{\kern6pt}\kern3pt\vrule}\hrule}}

\newcommand{\SpinC}{{\mathrm{Spin}}^c}
\newcommand{\Zmod}[1]{\Z/{#1}\Z}
\newcommand{\Q}{\mathbb{Q}}
\newcommand{\R}{\mathbb{R}}

\newcommand{\C}{\mathbb{C}}

\newcommand{\Z}{\mathbb{Z}}
\newcommand{\goesto}{\mapsto}
\newcommand{\cm}{\cdot}

\newcommand\spinc{\mathfrak s}

\newcommand\Field{\mathbb{F}}

\newcommand\fs{\mathfrak s}
\newcommand\spinct{\mathfrak t}

\newcommand\Gto{\check{G}}
\newcommand\Ffrom{\widehat F}
\newcommand\Fred{\overline F}
\newcommand\Gfrom{\widehat G}
\newcommand\fu{\mathfrak u}
\newcommand\pref{\sigma}

\newcommand\Fto{\check{F}}
\newcommand\Lto{\check L}

\newcommand\tHto{\Hto}

\newcommand\cycle{\eta}

\newcommand\mCP{{\overline{\mathbb{CP}}}^2}

\newcommand\PD{PD}

\providecommand{\dual}{\mathrm{D}}

\usepackage{amsmath,amsthm,amsfonts,amscd,flafter,graphicx}

\usepackage{rotating}   
\usepackage{pkmacros}

\title[{Monopoles and lens space surgeries}]
{Monopoles and lens space surgeries}

\author[Kronheimer]{Peter Kronheimer}
\address{Harvard University, Massachusetts
02138 \newline
\indent{\emailfont{kronheim@math.harvard.edu}}}
\thanks{PBK was partially supported by NSF grant number DMS-0100771.}

\author[Mrowka]{Tomasz Mrowka}
\address{Massachusetts Institute of Technology, Massachusetts 02139 \newline
\indent{Institute for Advanced Study, Princeton, New Jersey 08544} \newline
\indent{\emailfont{mrowka@math.mit.edu}}}
\thanks{TSM was partially supported by NSF grant numbers DMS-0206485,
DMS-0111298,
and FRG-0244663}

\author[Ozsv{\'a}th]{Peter Ozsv\'ath}
\address{Columbia University, New York 10027 \newline
\indent{Institute for Advanced Study, Princeton, New Jersey 08544}
\newline
\indent{\emailfont{petero@math.columbia.edu}}}
\thanks{PSO was partially supported by NSF grant numbers DMS-0234311,
DMS-0111298,
and FRG-0244663}

\author[Szab{\'o}]{Zolt{\'a}n Szab{\'o}} 
\address{Princeton University, New Jersey 08540 \newline
\indent{\emailfont{szabo@math.princeton.edu}}}
\thanks{ZSz was partially supported by NSF grant numbers DMS-0107792
and FRG-0244663,
and a Packard Fellowship.}

\begin{document}

\begin{abstract}
Monopole Floer homology is used to prove that real projective
three-space cannot be obtained from Dehn surgery on a non-trivial knot
in the three-sphere. To obtain this result, we use a surgery long
exact sequence for monopole Floer homology, together with a
non-vanishing theorem, which shows that monopole Floer homology
detects the unknot. In addition, we apply these techniques to give
information about knots which admit lens space surgeries, and to
exhibit families of three-manifolds which do not admit taut
foliations.
\end{abstract}

\maketitle
\section{Introduction}

Let $K$ be a knot in $S^3$. Given a rational number $r$, let
$S^3_r(K)$ denote the oriented three-manifold obtained from the knot
complement by Dehn filling with slope $r$.
The main purpose of this paper is to prove the following
conjecture of Gordon (see~\cite{GordonConjecture}, \cite{Gordon}):

\begin{theorem}
\label{thm:RP3}
Let $U$ denote the unknot in $S^3$, and let $K$ be any knot.
If there is an orientation-preserving diffeomorphism
$S^3_r(K)\cong S^3_r(U)$ for some rational number $r$,  then $K=U$.
\end{theorem}

To amplify the meaning of this result, we recall that $S^{3}_{r}(U)$
is the manifold $S^{1}\times S^{2}$ in the case $r=0$ and is a lens
space for all non-zero $r$. More specifically, with our conventions,
if $r=p/q$ in lowest
terms, with $p>0$, then $S^{3}_{r}(U) = L(p,q)$ as oriented manifolds.
The manifold $S^{3}_{p/q}(K)$ in general has first homology group
$\Z/p\Z$, independent of $K$. Because the lens space $L(2,q)$ is
$\mathbb{RP}^{3}$ for all odd $q$, the theorem implies (for example)
that $\mathbb{RP}^{3}$ cannot be obtained by Dehn
filling on a non-trivial knot.

Various cases of the Theorem~\ref{thm:RP3} were previously known.  The
case $r=0$ is the ``Property R'' conjecture, proved by
Gabai~\cite{GabaiKnots}, and the case where $r$ is
non-integral follows from the cyclic surgery theorem of Culler,
Gordon, Luecke, and Shalen~\cite{CGLS}.  The case where $r=\pm 1$ is a
theorem of Gordon and Luecke, see~\cite{GorLueckI}
and~\cite{GorLueckII}. Thus, the advance here is the case where $r$ is
an integer with $|r|>1$, though our techniques apply for any non-zero
rational $r$. In particular, we obtain an independent proof for the
case of the Gordon-Luecke theorem. (Gabai's result is an ingredient of
our argument.)

The proof of Theorem~\ref{thm:RP3} uses the Seiberg-Witten monopole
equations, and the monopole Floer homology package developed
in~\cite{KMBook}. Specifically, we use two properties of these
invariants. The first key property, which follows from the techniques
developed in \cite{KMcontact}, is a non-vanishing theorem for the
Floer groups of a three-manifold admitting a taut foliation. When
combined with the results of \cite{Gabai}, \cite{GabaiKnots}, this non-vanishing theorem
shows that Floer homology can be used to distinguish $S^{1}\times
S^{2}$ from $S^{3}_{0}(K)$ for non-trivial $K$.  The second property
that plays a central role in the proof is a surgery long exact
sequence, or exact triangle. Surgery long exact sequences of a related
type were introduced by Floer in the context of instanton Floer
homology, see~\cite{BraamDonaldson} and~\cite{FloerTriangles}.
The form of the surgery long exact sequence which is used in the
topological applications at hand is a natural analogue of a
corresponding result in the Heegaard Floer homology of~\cite{HolDisk}
and~\cite{HolDiskTwo}. In fact, the strategy of the proof presented
here follows closely the proof given in~\cite{BranchedDoubleCovers}.

Given these two key properties, the proof of Theorem~\ref{thm:RP3} has
the following outline. For integral $p$, we shall say that a knot $K$
is
\emph{$p$-standard} if $S^{3}_{p}(K)$ cannot be distinguished from
$S^{3}_{p}(U)$ by its Floer homology groups. (A more precise
definition is given in Section~\ref{sec:ProofOutline}, see
also Section~\ref{sec:NonVanishingProof}.) 
We can rephrase the non-vanishing theorem mentioned
above as the statement that, if $K$ is $0$-standard, then $K$
is unknotted. 
A surgery long exact sequence, involving the Floer homology groups of
$S^{3}_{p-1}(K)$, $S^{3}_{p}(K)$ and $S^{3}$, shows that if $K$ is
$p$-standard for $p>0$, then $K$ is also $(p-1)$ standard. By
induction, it follows that if $K$ is $p$-standard for some $p>0$, then
$K = U$. This gives the theorem for positive integers $p$. When $r>0$
is non-integral, we prove (again by using the surgery long exact sequence)
that if $S^3_r(K)$ is orientation-preservingly diffeomorphic to $S^3_r(U)$,
then $K$ is also $p$-standard, where $p$ is the smallest integer greater
than $r$. This proves Thoerem~\ref{thm:RP3} for all positive $r$.
The case of
negative $r$ can be deduced by changing orientations and replacing $K$
by its mirror-image.

As explained in Section~\ref{sec:FurtherApps}, the techniques
described here for establishing Theorem~\ref{thm:RP3} can be readily
adapted to other questions about knots admitting lens space
surgeries. For example, if $K$ denotes the $(2,5)$ torus knot, then it
is easy to see that $S^3_9(K)\cong L(9,7)$, and $S^3_{11}(K)\cong
L(11,4)$.  Indeed, a result described in Section~\ref{sec:FurtherApps}
shows that any lens space which is realized as integral surgery on a
knot in $S^3$ with Seifert genus two is diffeomorphic to one of these
two lens spaces.  Similar lists are given when $g=3$, $4$, and
$5$. Combining these methods with a result of Goda and Teragaito, we
show that the unknot and the trefoil are the only knots which admits a
lens space surgery with $p=5$.  In another direction, we give
obstructions to a knot admitting Seifert fibered surgeries, in terms
of its genus and the degree of its Alexander polynomial.

Finally, in Section~\ref{sec:Foliations}, we give some applications of
these methods to the study of taut foliations, giving several families
of three-manifolds which admit no taut foliation. One infinite family
of hyperbolic examples is provided by the $(-2,3,2n+1)$ pretzel knots
for $n\geq 3$: it is shown that all Dehn fillings with sufficiently
large surgery slope $r$ admit no taut foliation. The first examples of
hyperbolic three-manifolds with this property were constructed by
Roberts, Shareshian, and Stein in~\cite{RobertsShareshianStein}, see
also~\cite{CalegariDunfield}. In another direction, we show that if
$L$ is a non-split alternating link, then the double-cover of $S^3$
branched along $L$ admits no taut foliation. Additional examples
include certain plumbings of spheres and certain surgeries on the
Borromean rings, as described in this section.

\subsubsection*{Outline} The remaining sections of this paper are as
follows. In Section~\ref{sec:FloerOutline}, we give a summary of the
formal properties of the Floer homology groups developed
in~\cite{KMBook}. We do this in the simplest setting, where the
coefficients are $\Z/2$.  In this context we give precise statements
of the non-vanishing theorem and surgery exact sequence. With $\Z/2$
coefficients, the non-vanishing theorem is applicable only to knots
with Seifert genus $g>1$. In Section~\ref{sec:ProofOutline}, we use
the non-vanishing theorem and the surgery sequence to prove
Theorem~\ref{thm:RP3} for all integer $p$, under the additional
assumption that the genus is not $1$. (This is enough to cover all
cases of the theorem that do not follow from earlier known results,
because a result of Goda and Teragaito~\cite{GodaTeragaito} rules out
genus-$1$ counterexamples to the theorem.)

Section~\ref{sec:FloerDetails} describes some details of the
definition of the Floer groups, and the following two sections give
the proof of the surgery long-exact sequence
(Theorem~\ref{thm:TwoHandleSurgeries-A}) and the non-vanishing
theorem. In these three sections, we also introduce more general
(local) coefficients, allowing us to state the non-vanishing theorem
in a form applicable to the case of Seifert genus $1$. The surgery
sequence with local coefficients is stated as
Theorem~\ref{thm:TwoHandleSurgeries-Local}. In
Section~\ref{sec:NonVanishingProof}, we discuss
a refinement of the non-vanishing theorem using local coefficients.
At this stage we have the
machinery to prove Theorem~\ref{thm:RP3} for integral $r$ and any $K$,
without restriction on genus.
In Section~\ref{sec:Rational}, we explain how repeated applications of
the long exact sequence can be used to reduce the case of non-integral
surgery slopes to the case where the surgery slopes are integral, so
providing a proof of Theorem~\ref{thm:RP3} in the non-integral case
that is independent of the cyclic surgery theorem of \cite{CGLS}.

In Section~\ref{sec:FurtherApps}, we describe several further
applications of the same techniques to other questions involving
lens-space surgeries.  Finally, we give some applications of these
techniques to studying taut foliations on three-manifolds in
Section~\ref{sec:Foliations}.

\subsubsection*{Remark on orientations}
Our conventions about orientations and lens spaces have the following
consequences. If a $2$-handle is attached to the $4$-ball along an
attaching curve $K$ in $S^{3}$, and if the attaching map is chosen so
that the resulting $4$-manifold has intersection form $(p)$, then the
oriented boundary of the $4$-manifold is $S^{3}_{p}(K)$.  For positive
$p$, the lens space $L(p,1)$ coincides with $S^{3}_{p}(U)$ as an
oriented $3$-manifold. This is not consistent with the convention that
$L(p,1)$ is the quotient of $S^{3}$ (the oriented boundary of the unit
ball in $\C^{2}$) by the cyclic group of order $p$ lying in the center
of $U(2)$.

\subsubsection*{Acknowledgements} 
The authors wish to thank Cameron Gordon, John Morgan, and Jacob
Rasmussen for several very interesting discussions.  We are especially
indebted to Paul Seidel for sharing with us his expertise in
homological algebra. The formal aspects of the construction of the
monopole Floer homology groups described here have roots that can be
traced back to lectures given by Donaldson in Oxford in
1993. Moreover, we have made use of a Floer-theoretic construction of
Fr{\o}yshov, giving rise to a numerical invariant extending
the one which can be found in~\cite{Froyshov}. 

\section{Monopole Floer homology}
\label{sec:FloerOutline}


\subsection{The Floer homology functors}

We summarize the basic properties of the Floer groups constructed in
\cite{KMBook}.
In this section we will treat only monopole Floer homology with
coefficients in the field $\Field =\Z/2\Z$. Our three-manifolds will
always be smooth, oriented, compact, connected and without boundary
unless otherwise stated. To each such three-manifold $Y$, we associate
three vector spaces over $\Field$, 
\[
 \Hto_{\bullet}(Y), \quad \Hfrom_{\bullet}(Y), \quad \Hred_{\bullet}(Y).
\]
These are the monopole Floer homology groups, read ``HM-to'', 
``HM-from'', and ``HM-bar'' respectively.
They come equipped with linear maps $i_{*}$, $j_{*}$ and $p_{*}$ which
form a long exact sequence \begin{equation}\label{eq:ijp} \cdots
\stackrel{i_{*}}{\longrightarrow} \Hto_{\bullet}(Y)
\stackrel{j_{*}}{\longrightarrow} \Hfrom_{\bullet}(Y)
\stackrel{p_{*}}{\longrightarrow} \Hred_{\bullet}(Y)
\stackrel{i_{*}}{\longrightarrow} \Hto_{\bullet}(Y)
\stackrel{j_{*}}{\longrightarrow} \cdots .  \end{equation} A
\emph{cobordism} from $Y_{0}$ to $Y_{1}$ is an oriented, connected
$4$-manifold $W$ equipped with an orientation-preserving
diffeomorphism from $\partial W$ to the disjoint union of $-Y_{0}$ and
$Y_{1}$. We write $W: Y_{0} \to Y_{1}$. We can form a category, in
which the objects are three-manifolds, and the morphisms are
diffeomorphism classes of cobordisms. The three versions of
monopole Floer homology are
functors from this category to the category of vector spaces. That is,
to each $W:Y_{0}
\to Y_{1}$, there are associated maps
\[
\begin{aligned}
    \Hto(W) :\Hto_{\bullet}(Y_{0}) &\to
    \Hto_{\bullet}(Y_{1}) \\
    \Hfrom (W) :\Hfrom_{\bullet}(Y_{0}) &\to
    \Hfrom_{\bullet}(Y_{1}) \\
    \Hred (W) :\Hred_{\bullet}(Y_{0}) &\to
    \Hred_{\bullet}(Y_{1}) .
\end{aligned}
\]
The maps $i_{*}$, $j_{*}$ and $p_{*}$ provide natural transformations
of these functors. In addition to their vector space structure, the
Floer groups come equipped with a distinguished endomorphism, making
them modules over the polynomial ring $\Field[U]$. This module
structure is respected by the maps arising from cobordisms, as well as
by the three natural transformations.

These Floer homology groups are set up so as to be gauge-theory
cousins of the Heegaard homology groups $\mathit{HF}^{+}(Y)$,
$\mathit{HF}^{-}(Y)$ and $\mathit{HF}^{\infty}(Y)$ defined in
\cite{HolDisk}. Indeed, if $b_{1}(Y)=0$, then the monopole Floer
groups are conjecturally isomorphic to (certain completions of)
their Heegaard counterparts. 

\subsection{The non-vanishing theorem}

A \emph{taut foliation} $\Fol$ of an oriented $3$-manifold $Y$ is a
$C^{0}$ foliation of $Y$ with smooth, oriented $2$-dimensional leaves,
such that there exists a closed $2$-form $\omega$ on $Y$ whose
restriction to each leaf is everywhere positive. (Note that all
foliations which are taut in this sense are automatically
coorientable. 
There is a slightly weaker notion of tautness in the
literature which applies even in the non-coorientable case --
i.e. that there is a transverse curve which meets all the leaves. Of
course, when $H^1(Y;\Zmod{2})=0$, all foliations are coorientable, and
hence these two notions coincide.)  We write $e(\Fol)$ for the Euler
class of the $2$-plane field tangent to the leaves, an element of
$H^{2}(Y;\Z)$. The proof of the following theorem is based on the
techniques of
\cite{KMcontact} and makes use of the results of
\cite{EliashbergThurston}.  

\begin{theorem}
\label{thm:NonVanishing}
    Suppose $Y$ admits a taut foliation $\Fol$ and is not
    $S^{1}\times S^{2}$. If either (a) $b_{1}(Y)=0$, or (b) $b_{1}(Y)=1$
    and $e(\Fol)$ is non-torsion, then the image of
    $j_{*}:\Hto_{\bullet}(Y) \to \Hfrom_{\bullet}(Y)$ is non-zero.
\end{theorem}

The restriction to the two cases (a) and (b) in the statement of this
theorem arises from our use of Floer homology with coefficients
$\Field$. There is a quite general non-vanishing theorem for
$3$-manifolds satisfying the hypothesis in the first sentence; but for
this version (which is stated as Theorem~\ref{thm:NonVanishing-twist}
and proved in Section~\ref{sec:NonVanishingProof})
we need to use more general, local coefficients.

Note that $j_{*}$ for $S^2\times S^1$ is trivial in view of the following:

\begin{prop}
	\label{prop:PSC}
	If $Y$ is a three-manifold which admits a metric of positive
	scalar curvature, then the image of $j_{*}$ is zero.
\end{prop}

According to Gabai's theorem
\cite{GabaiKnots}, if $K$ is a non-trivial knot, then $S^{3}_{0}(K)$ admits
a taut foliation $\Fol$, and is not $S^{1}\times S^{2}$. If the Seifert genus
of $K$ is greater than $1$, then $e(\Fol)$ is non-torsion.
As a consequence, we have:

\begin{cor}
\label{cor:NonVanishing}
    The image of $j_{*} : \Hto_{\bullet}(S^{3}_{0}(K)) \to
    \Hfrom_{\bullet}(S^{3}_{0}(K))$ is non-zero if the Seifert genus
    of $K$ is $2$ or more, and is zero if $K$ is the unknot.
    \qed
\end{cor}

\subsection{The surgery exact sequence}
\label{subsec:Surgery-A}

Let $M$ be an oriented $3$-manifold with torus boundary. Let
$\gamma_{1}$, $\gamma_{2}$, $\gamma_{3}$ be three oriented simple
closed curves on $\partial
M$ with algebraic  intersection numbers
\[
(\gamma_{1}\cdot\gamma_{2}) = (\gamma_{2}\cdot \gamma_{3})  =
(\gamma_{3}\cdot \gamma_{1}) = -1.
\]
Define $\gamma_{n}$ for all $n$ so that $\gamma_{n} = \gamma_{n+3}$.
Let $Y_{n}$ be the closed $3$-manifold obtained by filling along
$\gamma_{n}$: that is, we attach $S^{1}\times D^{2}$ to $M$ so that
the curve $\{1\} \times \partial D^{2}$ is attached to $\gamma_{n}$.
There is a standard cobordism $W_{n}$ from $Y_{n}$ to $Y_{n+1}$. The
cobordism is obtained from $[0,1] \times Y_{n}$ by attaching a
$2$-handle to $\{1\}\times Y_{n}$, with framing $\gamma_{n+1}$.  Note
that these orientation conventions are set up so that
$W_{n+1}\cup_{Y_{n+1}} W_n$ always contains a sphere with
self-intersection number $-1$.

\begin{theorem}
\label{thm:TwoHandleSurgeries-A}
    There is an exact sequence
    \[
             \cdots \longrightarrow
             \Hto_{\bullet}(Y_{n-1})
             \stackrel{F_{n-1}}{\longrightarrow}
              \Hto_{\bullet}(Y_{n})
           \stackrel{F_{n}}{\longrightarrow}
               \Hto_{\bullet}(Y_{n+1}) \longrightarrow
             \cdots,
    \]
    in which the maps $F_{n}$ are given by the cobordisms $W_{n}$.
    The same holds for $\Hfrom_{\bullet}$ and $\Hred_{\bullet}$.
\end{theorem}

The proof of the theorem is given in
Section~\ref{sec:LongExactSequences}.

\subsection{Gradings and completions}
\label{subsec:Gradings}

The Floer groups are \emph{graded} vector spaces, but there are two
caveats: the grading is not by $\Z$, and a completion is involved. We
explain these two points.

Let $\Gr$ be a set with an action of $\Z$, not necessarily
transitive. We write $j \mapsto j+n$
for the action of $n\in \Z$ on $J$. A vector space $V$ is \emph{graded
by $\Gr$} if it is presented as a direct sum of subspaces $V_{j}$ indexed
by $J$. A homomorphism $h : V\to V'$ between vector spaces graded by
$J$ has degree $n$ if $h(V_{j})\subset V'_{j+n}$ for all $j$.

If $Y$ is an oriented $3$-manifold, we write $J(Y)$ for the
set of homotopy-classes of oriented $2$-plane fields (or equivalently
nowhere-zero vector fields) $\xi$ on $Y$. To define an
action of $\Z$, we specify that $[\xi] + n$ denotes
the homotopy class $[\tilde{\xi}]$ obtained from $[\xi]$ as
follows. Let $B^{3}\subset Y$ be a standard ball, and let
$\rho : (B^{3},\partial B^{3}) \to (SO(3), 1)$ be a map of
degree $-2n$, regarded as an automorphism of the trivialized
tangent bundle of the ball.
Outside the ball $B^{3}$, we take $\tilde{\xi} =
\xi$. Inside the ball, we define
\[
   \tilde{\xi}(y) = \rho(y)\xi(y).
\]
The structure of $J(Y)$ for a general three-manifold is as follows
(see \cite{KMcontact}, for example). A
$2$-plane field determines a $\SpinC$ structure on $Y$, so we can
first write
\[
          J(Y) = \bigcup_{\spinc\in \SpinC(Y)} J(Y,\spinc),
\]
where the sum is over all isomorphism classes of $\SpinC$ structures.
The action of $\Z$ on each $J(Y,\spinc)$ is transitive, and the stabilizer
is the subgroup of $2\Z$ given by the image of the map
\begin{equation}\label{eq:c1-map}
          x \mapsto \langle c_{1}(\spinc) , x\rangle
\end{equation}
from $H_{2}(Y;\Z)$ to $\Z$. In particular, if $c_{1}(\spinc)$ is torsion,
then $J(Y,\spinc)$ is an affine copy of $\Z$.

For each $j\in J(Y)$, there are subgroups
\[
\begin{aligned}
\Hto_{j}(Y) &\subset \Hto_{\bullet}(Y) \\ 
   \Hfrom_{j}(Y) &\subset \Hfrom_{\bullet}(Y) \\ 
   \Hred_{j}(Y) &\subset \Hred_{\bullet}(Y),
   \end{aligned}
\]
and there are internal direct sums which we denote by $\Hto_{*}$,
$\Hfrom_{*}$ and $\Hred_{*}$:
\[
\begin{aligned}
\Hto_{*}(Y) &= \bigoplus_{j} \Hto_{j}(Y) \subset
         \Hto_{\bullet}(Y) \\ 
         \Hfrom_{*}(Y) &= \bigoplus_{j} \Hfrom_{j}(Y) \subset
         \Hfrom_{\bullet}(Y) \\ 
         \Hred_{*}(Y) &= \bigoplus_{j} \Hred_{j}(Y) \subset
         \Hred_{\bullet}(Y) 
         \end{aligned}
\]
The $\bullet$ versions are obtained from the $*$ versions as follows.
For each $\spinc$ with $c_{1}(\spinc)$ torsion, pick an arbitrary
$j_{0}(\spinc)$ in $J(Y,\spinc)$. Define a decreasing sequence of subspaces
$\Hfrom[n] \subset \Hfrom_{*}(Y)$ by
\[
         \Hfrom[n] = \bigoplus_{\spinc}
         \bigoplus_{m\ge n} \Hfrom_{j_{0}(\spinc)-m}(Y),
\]
where the sum is over torsion $\SpinC$ structures. Make the same
definition for the other two variants. The groups $\Hto_{\bullet}(Y)$,
$\Hfrom_{\bullet}(Y)$ and $\Hred_{\bullet}(Y)$ are the completions of
the direct sums $\Hto_{*}(Y)$ etc.~with respect to these decreasing
filtrations. However, in the case of $\Hto$, the subspace $\Hto[n]$ is
eventually zero for large $n$, so the completion has no effect. 
From the decomposition of $J(Y)$ into orbits, we have direct sum
decompositions
\[
\begin{aligned}
    \Hto_{\bullet}(Y) = \bigoplus_{\spinc}
    \Hto_{\bullet}(Y,\spinc)\\
        \Hfrom_{\bullet}(Y) = \bigoplus_{\spinc}
    \Hfrom_{\bullet}(Y,\spinc)\\
    \Hred_{\bullet}(Y) = \bigoplus_{\spinc}
    \Hred_{\bullet}(Y,\spinc).
\end{aligned}
\]
Each of these decompositions has only finitely many non-zero terms.

The maps $i_{*}$, $j_{*}$ and $p_{*}$ are defined on the $*$ versions
and have degree $0$, $0$ and $-1$ respectively, while the endomorphism
$U$ has degree $-2$. The maps induced by cobordisms do not have a
degree and do not always preserve the $*$ subspace: they are continuous
homomorphisms between complete filtered vector spaces.

To amplify the last point above, consider a cobordism $W: Y_{0} \to
Y_{1}$. The homomorphisms $\Hto(W)$ etc.~can be written as sums
\[
            \Hto(W) = \sum_{\spinc} \Hto(W,\spinc),
\]
where the sum is over $\SpinC(W)$: for each $\spinc\in\SpinC(W)$, we have
\[
           \Hto(W,\spinc) : \Hto_{\bullet}(Y_{0},\spinc_{0}) \to
                            \Hto_{\bullet}(Y_{1},\spinc_{1}),
\]
where $\spinc_{0}$ and $\spinc_{1}$ are the resulting $\SpinC$
structures on the boundary components. 
The above sum is not 
necessarily finite,
but it is convergent.
The individual terms $\Hto(W,\spinc)$ have a well-defined degree, in
that for each $j_{0} \in J(Y_{0},\spinc_{0})$ there is a unique $j_{1}
\in J(Y_{1},\spinc_{1})$ such that
\[
       \Hto(W,\spinc)  : \Hto_{j_{0}}(Y_{0},\spinc_{0}) \to
                            \Hto_{j_{1}}(Y_{1},\spinc_{1}).
\]
The same remarks apply to $\Hfrom$ and $\Hred$.
The element $j_{1}$ can be characterized as follows. Let $\xi_{0}$ be
an oriented $2$-plane field in the class $j_{0}$, and let $I$ be an
almost complex structure on $W$ such that: (i) the planes $\xi_{0}$
are invariant under $I|_{Y_{0}}$ and have the complex orientation; and
(ii) the $\SpinC$ structure associated to $I$ is $\spinc$. Let $\xi_{1}$
be the unique oriented $2$-plane field on $Y_{1}$ that is invariant
under $I$. Then $j_{1} = [\xi_{1}]$. For future reference, we
introduce the notation
\[
          j_{0} \simspin{\spinc} j_{1}
\]
to denote the relation described by this construction.

\subsubsection{Remark}

Because of the completion involved in the definition of the Floer
groups, the $\Field[U]$-module structure of the groups
$\Hfrom_{*}(Y,\spinc)$ (and its companions) gives rise to an $\Field[[U]]$-module structure
on $\Hfrom_{\bullet}(Y,\spinc)$, whenever $c_{1}(\spinc)$ is torsion.
In the non-torsion case, the action of $U$ on $\Hfrom_{*}(Y,\spinc)$
is actually nilpotent, so again the action extends. In this way, each
of $\Hto_{\bullet}(Y)$, $\Hfrom_{\bullet}(Y)$ and $\Hred_{\bullet}(Y)$
become modules over $\Field[[U]]$, with continuous module
multiplication.

\subsection{Canonical mod 2 gradings}
\label{subsec:CanonModTwo}

The Floer groups have a canonical grading mod $2$. For a cobordism
$W:Y_{0} \to Y_{1}$, let us define
\[
       \iota(W) = \frac{1}{2}
       \bigl(\chi(W) + \sigma(W) - b_{1}(Y_{1}) +
       b_{1}(Y_{0})\bigr),
\]
where $\chi$ denotes the Euler number, $\sigma$ the signature, and
$b_1$ the first Betti number with real coefficients.  Then we have the
following proposition.

\begin{prop}
\label{prop:CanonicalMod2}
There is one and only one way to decompose the grading set $J(Y)$ for
all $Y$ into even and odd parts in such a way that the following two
conditions hold.
\begin{enumerate}
    \item The gradings $j\in J(S^{3})$ for which $\Hto_{j}(S^{3})$ is
    non-zero are even.
    \item If $W:Y_{0}\to Y_{1}$ is a cobordism and $j_{0}
    \simspin{\spinc} j_{1}$ for some $\SpinC$ structure $\spinc$ on $W$,
    then $j_{0}$ and $j_{1}$ have the same parity if and only if
    $\iota(W)$ is even.
\end{enumerate}
\end{prop}

This result gives provides  a canonical decomposition
\[
             \Hto_{\bullet}(Y) = \Hto_{\ev}(Y) \oplus \Hto_{\odd}(Y),
\]
with a similar decomposition for the other two flavors. With respect
to these mod $2$ gradings, the maps $i_{*}$ and $j_{*}$ in the long
exact sequence have even degree, while $p_{*}$ has odd degree. The
maps resulting from a cobordism $W$ have even degree if and only if
$\iota(W)$ is even.

\subsection{Computation from reducible solutions}
\label{subsec:Reducibles}

While the groups $\Hto_{\bullet}(Y)$ and $\Hfrom_{\bullet}(Y)$ are
subtle invariants of $Y$, the group $\Hred_{\bullet}(Y)$ by contrast
can be calculated knowing only the cohomology ring of $Y$. This is
because the definition of $\Hred_{\bullet}(Y)$ involves only the
\emph{reducible} solutions of the Seiberg-Witten monopole equations
(those where the spinor is zero). We discuss here the case that $Y$ is
a rational homology sphere.

When $b_{1}(Y)=0$, the number of different $\SpinC$ structures on $Y$ is
equal to the order of $H_{1}(Y;\Z)$, and $J(Y)$ is the union of the
same number of copies of $\Z$. The contribution to
$\Hred_{\bullet}(Y)$ from each $\SpinC$ structure is the same:

\begin{prop}
    \label{prop:Reducibles-Y}
Let $Y$ be a rational homology sphere and $\tspin$ a $\SpinC$
   structure on $Y$. Then
   \[
             \Hred_{\bullet}(Y,\tspin) \cong \Field[U^{-1},U]]
   \]
   as topological $\Field[[U]]$-modules, where the right-hand side
   denotes the ring of formal Laurent series in $U$ that are finite in the
   negative direction.
\end{prop}

The maps $\Hred_{\bullet}(W)$ arising from cobordisms between rational
homology spheres are also standard:

\begin{prop}
\label{prop:Reducibles}
 Suppose $W: Y_{0} \to Y_{1}$ is a cobordism between rational
    homology spheres, with $b_{1}(W)=0$, and suppose that the
    intersection form on $W$ is negative definite. Let $\spinc$ be a
    $\SpinC$ structure on $W$, and suppose $j_{0} \simspin{\spinc} j_{1}$.
    Then
    \[
                   \Hred(W,\spinc) : \Hred_{j_{0}}(Y_{0}) \to
                   \Hred_{j_{1}}(Y_{1})
    \]
    is an isomorphism. On the other hand, if the intersection form on
    $W$ is not negative definite, then $\Hred(W,\spinc)$ is zero, for
    all $\spinc$.
\end{prop}

The last part of the proposition above holds in a more general form.
Let $W$ be a cobordism between $3$-manifolds that are not necessarily
rational homology spheres, and let $b^{+}(W)$ denote the dimension of
a maximal positive-definite subspace for the quadratic form on the
image of $H^{2}(W,\partial W;\R)$ in $H^{2}(W;\R)$.

\begin{prop}
    \label{prop:Reducibles-II}
    If the cobordism $W : Y_{0} \to Y_{1}$ has $b^{+}(W) > 0$,
    then the map $\Hred(W)$ is zero.
\end{prop}

\subsection{Gradings and rational homology spheres}

We return to rational homology spheres, and cobordisms between them.
If $W$ is such a cobordism, then $H^{2}(W,\partial W;\Q)$ is
isomorphic to $H^{2}(W;\Q)$, and there is therefore a quadratic form
\[
        Q : H^{2}(W;\Q) \to \Q
\]
given by $Q(e) = (\bar{e}\smile\bar{e})[W,\partial W]$, where $\bar{e} \in
H^{2}(W,\partial W;\Q)$ is a class whose restriction to $W$ is $e$.
We will simply write $e^{2}$ for
$Q(e)$.

\begin{lemma}
    \label{lemma:Proto-Q-grading}
    Let $W, W': Y_{0} \to Y_{1}$ be two cobordisms between a pair of
    rational homology spheres $Y_{0}$ and
    $Y_{1}$. Let $j_{0}$ and $j_{1}$ be classes of oriented $2$-plane
    fields on the $3$-manifolds and suppose that
    \[
                     \begin{aligned}
                        j_{0}& \simspin{\spinc} j_{1} \\
                         j_{0}& \simspin{\spinc'} j_{1} \\
                     \end{aligned}
    \]
    for $\SpinC$ structure $\spinc$ and $\spinc'$ on the two
    cobordisms. Then
    \[
               c_{1}^{2}(\spinc) - 2\chi(W) - 3 \sigma(W) =
               c_{1}^{2}(\spinc') - 2\chi(W') - 3 \sigma(W'),
    \]
    where $\chi$ and $\sigma$ denote the Euler number and signature.  
\end{lemma}

\begin{proof}
    Every $3$-manifold equipped with a $2$-plane field $\xi$ is the
    boundary of some almost-complex manifold $(X,I)$ in such a way
    that $\xi$ is invariant under $I$; so bearing in mind the
    definition  of the relation $\simspin{\spinc}$, and using the
    additivity of all the terms involved, we can reduce the
    lemma to a statement about closed almost-complex manifolds. The
    result is thus a consequence of the fact that
    \[
          c_{1}^{2}(\spinc)[X] - 2\chi(X) - 3 \sigma(X) = 0
    \]
    for the canonical $\SpinC$ structure on a closed, almost-complex
    manifold $X$.
\end{proof}

Essentially the same point leads to the definition of the following $\Q$-valued
function on $J(Y)$, and the proof that it is well-defined:

\begin{defn}
    \label{def:Q-grading}
    For a three-manifold $Y$ with $b_{1}(Y) = 0$ and
    $j \in J(Y)$ represented by an oriented $2$-plane field $\xi$,
    we define $h(j) \in \Q$ by the formula
    \[
                     4 h(j) = c_{1}^{2}(X,I) - 2\chi(X) - 3\sigma(X) +
                     2,
    \]
    where $X$ is a manifold whose oriented
    boundary is $Y$, and $I$ is an almost-complex structure such that the $2$-plane field $\xi$ is
    $I$-invariant and has the complex orientation. The quantity
    $c_{1}^{2}(X,I)$ is to be interpreted again using the natural isomorphism
    $H^{2}(X,\partial X;\Q) \cong H^{2}(X;\Q)$. 
\end{defn}

The map $h: J(Y)\to \Q$ satisfies $h(j+1) = h(j)+1$.

Now let $\spinc$ be a $\SpinC$ structure on a rational homology sphere
$Y$, and consider the exact sequence
\begin{equation}\label{eq:ImISequence}
          0 \to \mathrm{im}(p_{*}) \hookrightarrow
          \Hred_{\bullet}(Y,\spinc) \stackrel{i_{*}}{\longrightarrow}
          \mathrm{im}(i_{*}) \to 0,
\end{equation}
where $p_{*} : \Hfrom_{\bullet}(Y,\spinc)  \to
\Hred_{\bullet}(Y,\spinc)$. The image of $p_{*}$ is a closed,
non-zero, proper $\Field[U^{-1}, U]]$-submodule of
$\Hred_{\bullet}(Y,\spinc)$; and the
latter is isomorphic to $\Field[U^{-1}, U]]$ by
Proposition~\ref{prop:Reducibles-Y}. The only such submodules of
$\Field[U^{-1}, U]]$ are the submodules $U^{r}\Field[[U]]$ for $r\in
\Z$. It follows that the short exact sequence above is isomorphic to
the short exact sequence
\[
       0 \to \Field[[U]] \to \Field[U^{-1}, U]]\to
       \Field[U^{-1}, U]]/\Field[[U]] \to 0.
\]
This observation leads to a
$\Q$-valued invariant of $\SpinC$ structures on rational homology
spheres, after Fr{\o}yshov \cite{Froyshov}:

\begin{defn}
    Let $Y$ be an oriented rational homology sphere and $\spinc$ a
    $\SpinC$ structure. We define (by either of two equivalent
    formulae)
    \[
    \begin{aligned}
    \Froy(Y,\spinc) &= \min\{ \,h(k) \mid 
                      \text{\normalfont $i_{*}: \Hred_{k}(Y,\spinc) \to
                                     \Hto_{k}(Y,\spinc)$ is
                                     non-zero}\,\},\\
                    &= \max\{ \,h(k) + 2 \mid 
                      \text{\normalfont
                      $p_{*}: \Hfrom_{k+1}(Y,\spinc) \to
                                     \Hred_{k}(Y,\spinc)$ is
                                     non-zero}\,\}.                  
                                     \end{aligned}
    \]
\end{defn}

When $j_{*}$ is zero, sequence \eqref{eq:ImISequence} determines
everything, and we have:

\begin{cor}\label{cor:Y-j-Froy}
Let $Y$ be a rational homology sphere for which
the map $j_{*}$ is zero. Then for each $\SpinC$ structure $\spinc$, 
the short exact sequence
\[
   0\to     \Hfrom_{\bullet}(Y,\spinc) \stackrel{p_{*}}{\longrightarrow}
        \Hred_{\bullet}(Y,\spinc)
      \stackrel{i_{*}}{\longrightarrow} \Hto_{\bullet}(Y,\spinc) \to 0
\]
is isomorphic as a sequence of topological $\Field[[U]]$-modules to the sequence
\[
        0\to   \Field[[U]] \to  \Field[U^{-1}, U]] \to \Field[U^{-1}, U]] /
        \Field[[U]]\to 0.
\]
Furthermore, if $j_{\min}$ denotes the lowest degree in which
$\Hto_{j_{\min}}(Y,\spinc)$ is non-zero, then $h(j_{\min})= \Froy(Y,\spinc)$.
\end{cor}

\subsection{The conjugation action}
\label{subsec:Conjugation}

Let $Y$ be a three-manifold, equipped with a spin bundle $W$. The
bundle ${\overline W}$ which is induced from $W$ with the conjugate
complex structure naturally inherits a Clifford action from the one on
$W$. This correspondence induces an involution on the set of $\SpinC$
structures on $Y$, denoted $\spinc\mapsto {\overline \spinc}$.

Indeed, this conjugation action descends to an action on the Floer homology 
groups:

\begin{prop}
Conjugation induces a well-defined involution on $\Hto_\bullet(Y)$,
sending $\Hto(Y,\spinc) \mapsto \Hto(Y,{\overline \spinc})$. Indeed,
conjugation induces involutions on the other two theories as well,
which are compatible with the maps $i_{*}$, $j_{*}$, and $p_{*}$.
\end{prop}

\newcommand\ft{\mathfrak t}
\section{Proof of Theorem~\ref{thm:RP3} in the simplest cases}
\label{sec:ProofOutline}

In this section, we prove Theorem~\ref{thm:RP3} for the case that the
surgery coefficient is an integer and the Seifert genus of 
$K$ is not $1$.

\subsection{The Floer groups of lens spaces}
\label{subsec:LensSpace-I}

We begin by describing the Floer groups of the
$3$-sphere. There is only one
$\SpinC$ structure on $S^{3}$, and $j_{*}$ is zero because there is a
metric of positive scalar curvature. Corollary~\ref{cor:Y-j-Froy} is
therefore applicable. It remains only to say what $j_{\min}$ is, or
equivalently what the Fr{\o}yshov invariant is.



Orient $S^{3}$ as the boundary of the unit ball in $\R^{4}$ and
    let $\mathrm{SU}(2)_{+}$ and $\mathrm{SU}(2)_{-}$ be the subgroups
    of $\mathrm{SO}(4)$ that act trivially on the 
    anti-self-dual and self-dual $2$-forms
    respectively.
    Let $\xi_{+}$
    and $\xi_{-}$ be $2$-plane fields invariant under
    $\mathrm{SU}(2)_{-}$ and $\mathrm{SU}(2)_{+}$ respectively. Our
    orientation conventions are set up so that $[\xi_{-}] =  [\xi_{+}]
    + 1$.

\begin{prop}\label{prop:FroyS3}
  The least $j\in J(S^{3})$ for which $\Hto_{j}(S^{3})$ is non-zero is
  $j = [\xi_{-}]$. The largest $j\in J(S^{3})$ for which
  $\Hfrom_{j}(S^{3})$ is non-zero is $[\xi_{+}] = [\xi_{-}] - 1$. The
  Fr{\o}yshov invariant of $S^{3}$ is therefore given by:
    \[
    \Froy(S^{3}) = h([\xi_{-}])
                      = 0.
    \]
\end{prop}

%

We next describe the Floer groups for the lens
space $L(p,1)$, realized as $S^{3}_{p}(U)$ for an integer $p>0$.  The short
description is provided by Corollary~\ref{cor:Y-j-Froy}, because
$j_{*}$ is zero. To give a longer answer, we must describe the
$2$-plane field in which the generator of $\Hto$ lies, for each
$\SpinC$ structure. Equivalently, we must give the Fr{\o}yshov
invariants.

We first pin down the grading set $J(Y)$ for $Y =
S^{3}_{p}(K)$ and $p>0$. For a general knot
$K$, we have a cobordism
\[
          W(p): S^{3}_{p}(K) \to S^{3},
\]
obtained by the addition of a single $2$-handle. The manifold $W(p)$ has
$H_{2}(W(p)) = \Z$, and a generator has self-intersection number $-p$. A
choice of orientation for a Seifert surface for $K$ picks out a
generator $h=h_{W(p)}$. 
For each integer $n$, there is a unique $\SpinC$ structure $\spinc_{n,p}$ on
$W(p)$ with
\begin{equation}\label{eq:spinc-n}
     \langle   c_{1}(\spinc_{n,p}) , h \rangle = 2n-p.
\end{equation}
We denote the $\SpinC$ structure on $S^{3}_{p}(K)$ which arises from
$\spinc_{n,p}$ by $\tspin_{n,p}$; it depends only on $n$ mod $p$. Define
$j_{n,p}$ to be the unique element of $J(S^{3}_{p}(K), \tspin_{n,p})$
satisfying
\[
        j_{n,p} \simspin{\spinc_{n,p}} [\xi_{+}],
\]
where $\xi_{+}$ is the $2$-plane field on $S^{3}$ described above.
Like $\tspin_{n,p}$, the class $j_{n,p}$ depends on our choice of
orientation for the Seifert surface.
Our convention implies that $j_{0,1} = [\xi_{+}]$ on $S^{3}_{1}(U)
= S^{3}$.  If
$n \equiv n'$ mod $p$, then $j_{n,p}$ and $j_{n',p}$ belong to the same
$\SpinC$ structure, so they differ by an element of $\Z$ acting on
$J(Y)$. The next lemma calculates that element of $\Z$.
\begin{lemma}\label{lemma:j-Difference}
We have
\begin{equation*}
        j_{n,p} - j_{n',p} = \frac{(2n-p)^{2} - (2n'-p)^{2}}{4p}.
\end{equation*}
\end{lemma}

\begin{proof}
    We can equivalently calculate $h(j_{n,p}) - h(j_{n',p})$. We can
    compare $h(j_{n,p})$ to $h([\xi_{+}])$ using the cobordism $W(p)$,
    which tells us
    \[
                     4 h(j_{n,p}) = 4 h([\xi_{+}])
                          - c_{1}^{2}(\spinc_{n,p}) + 2\chi(W(p)) +
                          3\sigma(W(p)),
    \]
    and hence
    \[
    \begin{aligned}
    4 h(j_{n,p})  &=  -4 + \frac{(2n-p)^{2}}{p} +2 - 3 \\
                  &= \frac{(2n-p)^{2}}{p} - 5.
    \end{aligned}
    \]
    The result follows.
\end{proof}

Now we can state the generalization of Proposition~\ref{prop:FroyS3}.

\begin{prop}\label{prop:FroyS3p}
  Let $n$ be in the range $0\le n \le p$.
  The least $j\in J(Y, \tspin_{n,p})$ for which
  $\Hto_{j}(S^{3}_{p}(U), \tspin_{n,p})$ is non-zero is
  $j_{n,p}+1$. The largest $j\in J(S^{3}_{p}(U),\tspin_{n})$ for which
  $\Hfrom_{j}(S^{3}_{p}(U),\tspin_{n,p})$ is non-zero is $j_{n,p}$.
  Equivalently, the
  Fr{\o}yshov invariant of $(S^{3}_{p}(U),\tspin_{n,p})$ is given by:
    \begin{equation}\label{eq:FroyLens}
    \begin{aligned}
    \Froy(S^{3}_{p}(U),\tspin_{n,p}) &= h(j_{n,p})+1 \\
                      &= \frac{(2n-p)^{2}}{4p} -\frac{1}{4}.
                      \end{aligned}
    \end{equation}
\end{prop}

The meaning of this last result may be clarified by the following
remarks. By Proposition~\ref{prop:Reducibles}, we have an isomorphism
\[
          \Hred(W(p),\spinc_{n,p}) : \Hred_{j_{n,p}}(S^{3}_{p}(U),
          \tspin_{n,p}) \to \Hred_{[\xi+]}(S^{3});
\]
and because $j_{*}$ is zero for lens spaces, the map \[p_{*} :
\Hfrom_{j_{n,p}+1}(S^{3}_{p}(U),
          \tspin_{n,p}) \to \Hred_{j_{n,p}}(S^{3}_{p}(U),
          \tspin_{n,p})\] is an isomorphism.
Proposition~\ref{prop:FroyS3p} is therefore equivalent to the
following corollary:

\begin{cor}
\label{cor:FroyS3p}
The map
    \[
            \Hfrom(W(p), \spinc_{n,p}) :
             \Hfrom_{\bullet}(S^{3}_{p}(U), \tspin_{n,p})\to
             \Hfrom_{\bullet}(S^{3})
\]
is an isomorphism, whenever $0\le n\le p$.
\end{cor}

A proof directly from the definitions is sketched in
Section~\ref{subsec:LensSpace-II}. See also
Proposition~\ref{prop:Inequality}, which yields a more general result
by a more formal argument.

We can now be precise about what it means for $S^{3}_{p}(K)$ to
resemble $S^{3}_{p}(U)$ in its Floer homology.

\begin{defn}
    For an integer $p>0$, we say that $K$ is $p$-standard if \begin{enumerate}
    \item the map $j_{*}:\Hto_{\bullet}(S^{3}_{p}(K)) \to
    \Hfrom_{\bullet}(S^{3}_{p}(K))$ is zero; and
    \item for $0\le n \le p$, the Fr{\o}yshov invariant of the $\SpinC$
    structure $\tspin_{n,p}$ on $S^{3}_{p}(K)$ is given by the same
    formula \eqref{eq:FroyLens} as in the case of the unknot.
    \end{enumerate} For $p=0$, for the sake of
    expediency, we say that $K$ is weakly $0$-standard if the map $j_{*}$ is
    zero for $S^{3}_{0}(K)$.
\end{defn}

Observe that $\ft_{n,p}$ depended on an orientation
Seifert surface for the knot $K$. Letting $\ft_{n,p}^+$ and $\ft_{n,p}^-$
be the two possible choices using the two orientations of the Seifert surface,
it is easy to see that  $\ft_{n,p}^+$ is the conjugate
of $\ft_{n,p}^-$. In fact, since the Fr{\o}yshov invariant is invariant
under conjugation, it follows that our notation of $p$-standard
is independent of the choice of orientation.

If $p>0$ and $j_{*}$ is zero, the second condition in the definition is
equivalent to the assertion that $\Hfrom(W(p),\spinc_{n,p})$ is an
isomorphism for $n$ in
the same range:

\begin{cor}
\label{cor:HfromOnto}
    If $K$ is $p$-standard and $p>0$, then
    \[
            \Hfrom(W(p), \spinc_{n,p}) :
             \Hfrom_{\bullet}(S^{3}_{p}(K), \tspin_{n,p})\to
             \Hfrom_{\bullet}(S^{3})
\]
is an isomorphism for $0\le n \le p$. Conversely, if $j_{*}$ is zero
for $S^{3}_{p}(K)$ and the above map is an isomorphism for $0\le n\le
p$, then $K$ is $p$-standard.
\end{cor}

%

The next lemma tells us that a counter-example to
Theorem~\ref{thm:RP3} would be a $p$-standard knot.

\begin{lemma}
\label{lemma:DiffeoImpliesStandard}
If $S^{3}_{p}(K)$ and $S^{3}_{p}(U)$ are orientation-preserving
diffeomorphic for some integer $p>0$, then $K$ is $p$-standard.
\end{lemma}

\begin{proof}
Fix an integer $n$, and let $\psi : S^{3}_{p}(K) \to S^{3}_{p}(U)$ be
a diffeomorphism. To avoid ambiguity, let us write $\tspin_{n,p}^{K}$ and
$\tspin_{n,p}^{U}$ for the $\SpinC$ structures on these two
$3$-manifolds, obtained as above. Because $j_{*}$ is zero for
$S^{3}_{p}(K)$ and $\Hred(W(p),\spinc_{n,p})$ is an isomorphism, the map
\[
          \Hfrom(W(p),\spinc_{n,p}) : \Hfrom_{\bullet}(S^{3}_{p}(K),
          \tspin_{n,p}^{K}) \to \Hfrom_{\bullet}(S^{3}) 
\]
is injective. Making a comparison with Corollary~\ref{cor:HfromOnto},
we see that
\[
             \Froy(S^{3}_{p}(K), \tspin_{n,p}^{K}) \le
             \Froy(S^{3}_{p}(U),
          \tspin_{n,p}^{U})
\]
for $0\le n\le p$. So
\[
            \sum_{n=0}^{p-1} \Froy(S^{3}_{p}(K), \tspin_{n,p}^{K}) \le
             \sum_{n=0}^{p-1} \Froy(S^{3}_{p}(U), \tspin_{n,p}^{U}) .
\]
On the other hand, the as $n$ runs from $0$ to $p-1$, we run through
all $\SpinC$ structures once each; and because the manifolds are
diffeomorphic, we must have equality of the sums. The Fr{\o}yshov
invariants must therefore agree term by term, and $K$ is therefore
$p$-standard.
\end{proof}

\subsection{Exploiting the surgery sequence}


When the surgery coefficient is an integer and the genus is not $1$, 
Theorem~\ref{thm:RP3} is now a
consequence of the following proposition and
Corollary~\ref{cor:NonVanishing}, whose statement we can rephrase as
saying that a weakly
$0$-standard knot has genus $1$ or is unknotted.

\begin{prop}
\label{prop:GoingDown}
	If $K$ is $p$-standard for some integer $p\ge 1$, then $K$
	is weakly $0$-standard. 
\end{prop}

\begin{proof}
Suppose that $K$ is $p$-standard, so that in particular, $j_{*}$ is
zero for $S^{3}_{p}(K)$. We apply
Theorem~\ref{thm:TwoHandleSurgeries-A} to the following sequence of
cobordisms
\[
        \cdots \stackrel{}{\longrightarrow} S^{3}_{p-1}(K)
              \stackrel{W_{0}}{\longrightarrow}  S^{3}_{p}(K)
              \stackrel{W_{1}}{\longrightarrow}  S^{3}
              \stackrel{W_{2}}{\longrightarrow} S^{3}_{p-1}(K)
              \stackrel{}{\longrightarrow} \cdots
\]
to obtain a commutative diagram with exact rows and columns,
\[
\small
   \begin{CD}
        @. @VV{j_{*}}V @VV{0}V @VV{0}V @. \\
    @>>> \Hfrom_{\bullet}(S^{3}_{p-1}(K)) @>\Hfrom(W_{0})>>
      \Hfrom_{\bullet}(S^{3}_{p}(K)) @>\Hfrom(W_{1})>> \Hfrom_{\bullet}(S^{3}) @ >>>
     \\
      @. @VV{p_{*}}V @VV{p_{*}}V @VV{p_{*}}V @.\\
      @>>> \Hred_{\bullet}(S^{3}_{p-1}(K)) @>\Hred(W_{0})>>
      \Hred_{\bullet}(S^{3}_{p}(K)) @>\Hred(W_{1})>> \Hred_{\bullet}(S^{3}) @ >>>
       \\
              @. @VV{i_{*}}V @VV{i_{*}}V @VV{i_{*}}V @. 
      \\
      @>>> \Hto_{\bullet}(S^{3}_{p-1}(K)) @>\Hto(W_{0})>>
      \Hto_{\bullet}(S^{3}_{p}(K)) @>\Hto(W_{1})>> \Hto_{\bullet}(S^{3}) @ >>>
       \\
        @. @VV{j_{*}}V @VV{0}V @VV{0}V @. 
   \end{CD}
\]
In the case $K=U$, the cobordism $W_{1}$ is diffeomorphic (preserving
orientation) to $N\setminus B^{4}$, where $N$ is a tubular
neighborhood of a $2$-sphere with self-intersection number $-p$; and
$W_{2}$ has a similar description, containing a sphere with
self-intersection $(p-1)$. In general, the cobordism $W_{1}$ is the
manifold we called $W(p)$ above.

%

\begin{lemma}
    The maps \[
    \begin{aligned}
    \Hred(W_{1}) :\Hred_{\bullet}(S^{3}_{p}(K)) &\to
    \Hred_{\bullet}(S^{3}) \\
     \Hfrom(W_{1}) :\Hfrom_{\bullet}(S^{3}_{p}(K)) &\to
    \Hfrom_{\bullet}(S^{3}) 
    \end{aligned}
\]
    are zero if $p=1$ and are surjective if
    $p\ge 2$.
\end{lemma}

\begin{proof}
We write
\[
          \Hfrom_{\bullet}(S^{3}_{p}(K)) = \bigoplus_{n=0}^{p-1}
          \Hfrom_{\bullet}(S^{3}_{p}(K),\tspin_{n,p}). 
\]
If $n$ in the range $0\le n \le p-1$,
the map $\Hfrom_{\bullet}(W_{1}, \spinc_{n,p})$ is an isomorphism by
Corollary~\ref{cor:HfromOnto}, which gives identifications
\[
\begin{CD}
    \Hfrom_{\bullet}(S^{3}_{p}(K),\tspin_{n,p})
    @>{\Hfrom(W_{1},\spinc_{n,p})}>> \Hfrom_{\bullet}(S^{3}) \\
    @VVV  @VVV \\
    \Field[[U]] @>1>> \Field[[U]].
\end{CD}
\]
For $n' \equiv n \bmod{p}$, under the same identifications,
$\Hfrom(W_{1},\spinc_{n',p})$ becomes multiplication by $U^{r}$, where
\[
    r = (j_{n',p} - j_{n,p})/2.
\]
This difference was calculated in Lemma~\ref{lemma:j-Difference}.
Taking the sum over all $\spinc_{n',p}$, we see that
\[
            \sum_{n'\equiv n \;(p)} \Hfrom(W_{1},\spinc_{n',p}) :
            \Hfrom_{\bullet}(S^{3}_{p}(K), \tspin_{n,p}) \to
            \Hfrom_{\bullet}(S^{3})
\]
is isomorphic (as a map of vector spaces) to the map 
    $\Field[[U]] \to \Field[[U]] 
$
given by multiplication by the series
\[
          \sum_{n'\equiv n \;(p)} U^{((2n'-p)^{2} - (2n-p)^{2})/8p}
          \in \Field[[U]] .
\]
When $n=  0$, this series is $0$ as the terms
cancel in pairs. For all other $n$ in the range $1\le n \le p-1$,
the series has leading coefficient $1$ (the contribution from
$n'=n$) and is therefore invertible. Taking the sum over all residue classes, we obtain the
result for $\Hfrom$. The case of $\Hred$ is similar, but does not
depend on Corollary~\ref{cor:HfromOnto}.
\end{proof}

We can now prove Proposition~\ref{prop:GoingDown} by induction on $p$.
Suppose first that $p\ge 2$ and let $K$ be $p$-standard.
The lemma above tells us that $\Hfrom(W_{1})$ is surjective, and from
the exactness of the rows it follows that $\Hfrom(W_{0})$ is
injective. Commutativity of the diagram shows that $\Hred(W_{0})\circ
p_{*}$ is injective, where $p_{*} : \Hfrom_{\bullet}(S^{3}_{p-1}(K)) \to
\Hred_{\bullet}(S^{3}_{p-1}(K))$. It follows that
\[
           j_{*} : \Hto_{\bullet}(S^{3}_{p-1}(K)) \to
           \Hfrom_{\bullet}(S^{3}_{p-1}(K))
\]
is zero, by exactness of the columns. To show that $K$ is
$(p-1)$-standard, we must examine its Fr{\o}yshov invariants.

Fix $n$ in the range $0\le n \le p-2$, and let
\[
               e \in \Hred_{j_{n,p-1}}(S^{3}_{p-1}(K),\tspin_{n,p-1})
\]
be the generator.
To show that the Fr{\o}yshov
invariants of $S^{3}_{p-1}(K)$ are standard is to show that
\[
            e \in
            \mathrm{image}\bigl(  p_{*} : 
            \Hfrom_{\bullet}(S^{3}_{p-1}(K)) \to
            \Hred_{\bullet}(S^{3}_{p-1}(K))  \bigr).
\]
From the diagram, this is equivalent to showing
\[
       \Hred(W_{0})(e) \in \mathrm{image}\bigl(  p_{*} : 
            \Hfrom_{\bullet}(S^{3}_{p}(K)) \to
            \Hred_{\bullet}(S^{3}_{p}(K))  \bigr).
\]
Suppose on the contrary that $\Hred(W_{0})(e)$ does not belong to the
image of $p_{*}$. This means that there is a $\SpinC$ structure
$\uspin$ on $W_{0}$ such that
\[
          j_{n,p-1} \simspin{\uspin} j_{m,p} + x
\]
for some integer $x>0$, and $m$ in the range $0\le m \le p-1$. There is
a
unique $\SpinC$ structure $\wspin$  on the composite cobordism
\[
    X =  W_{1}\circ W_{0}  : S^{3}_{p-1}(K) \to S^{3}
\]
whose restriction to $W_{0}$ is $\uspin$ and whose restriction to
$W_{1}$ is $\spinc_{m,p}$. We have
\[
           j_{n,p-1} \simspin{\wspin} [\xi_{+}] + x.
\]
On the other hand, the composite cobordism $X$ is diffeomorphic to the
cobordism $W(p-1) \# \mCP$ (a fact that we shall return to in
Section~\ref{sec:LongExactSequences}), and we can therefore write (in
a self-evident notation)
\[
        \wspin = \spinc_{n',p-1} \# \spinc
\]
for some $\SpinC$ structure $\spinc$ on $\mCP$, and some $n'$
equivalent to $n$ mod $p$. From
Lemma~\ref{lemma:Proto-Q-grading}, we have
\begin{multline*}
         c_{1}^{2}(\wspin)  - 2\chi(X) 
       - 3\sigma(X) 
         \\ = 4x + c_{1}^{2}(\spinc_{n,p-1})  - 2\chi(W(p-1))- 3\sigma(W(p-1))
\end{multline*}
or in other words
\[
     \frac{(2n-p+1)^{2} - (2n'-p+1)^{2}}{p-1} + c_{1}^{2}(\spinc) + 1 = 4x.
\]
But $n$ is in the range $0\le n \le p-1$ and $c_{1}^{2}(\spinc)$ has
the form $-(2k+1)^{2}$ for some integer $k$, so the left hand side is
not greater than $0$. This contradicts the assumption that $x$ is
positive, and completes the argument for the case $p\ge 2$.

In the case $p=1$, the maps $\Hred(W_{1})$, $\Hfrom(W_{1})$ and
$\Hto(W_{1})$ are all zero.  A diagram chase
again shows that $j_{*}$ is zero for $S^{3}_{0}(K)$, so $K$ is
weakly $0$-standard.
\end{proof}

\section{Construction of monopole Floer homology}
\label{sec:FloerDetails}

\subsection{The configuration space and its blow-up}

Let $Y$ be an oriented $3$-manifold, equipped with a Riemannian
metric. Let $\bonf(Y)$ denote the space of isomorphism classes of
triples $(\spinc, A, \Phi)$, where $\spinc$ is a $\SpinC$ structure,
$A$ is a $\SpinC$ connection of Sobolev class $L^{2}_{k-1/2}$ in the
associated spin bundle $S\to Y$, and $\Phi$ is an $L^{2}_{k-1/2}$ section
of $S$. Here $k-1/2$ is any suitably large Sobolev exponent, and we
choose a half-integer because there is a continuous restriction map $L^{2}_{k}(X)
\to L^{2}_{k-1/2}(Y)$ when $X$ has boundary $Y$.   The
space $\bonf(Y)$ has one component for each isomorphism class of
$\SpinC$ structure, so we can write
\[
            \bonf(Y) = \bigcup_{\spinc} \bonf(Y,\spinc).
\]
We call an element of $\bonf(Y)$ \emph{reducible} if $\Phi$ is zero
and irreducible otherwise. If we choose a particular $\SpinC$
structure from each isomorphism class, we can construct a space
\[
           \conf(Y) = \bigcup_{\spinc} \conf(Y,\spinc),
\]
where $\conf(Y,\spinc)$ is the space of all pairs $(A,\Phi)$, a
$\SpinC$ connection and section for the chosen $S$. Then we can regard
$\bonf(Y)$ as the quotient of $\conf(Y)$ by the gauge group $\G(Y)$ of
all maps $u : Y\to S^{1}$ of class $L^{2}_{k+1/2}$.

The space $\bonf(Y)$ is a Banach manifold except at
the locus of reducibles; the reducible locus $\bonf^{\red}$(Y) is itself a
Banach manifold, and the map
\[
\begin{gathered}
    \bonf(Y) \to \bonf^{\red}(Y)\\
    [\spinc, A,\Phi] \mapsto [\spinc,A,0]
\end{gathered}
\]
has fibers $L^{2}_{k-1/2}(S)/S^{1}$, which is a cone on a complex
projective space. We can resolve the singularity along the reducibles
by forming a real, oriented blow-up,
\[
            \pi : \bonf^{\sigma}(Y) \to \bonf(Y).
\]
We define $\bonf^{\sigma}(Y)$ to be the space of isomorphism classes
of quadruples $(\spinc, A, s, \phi)$, where $\phi$ is an element of
$L^{2}_{k-1/2}(S)$ with unit $L^{2}$ norm and $s\ge 0$. The map $\pi$ is
\[
      \pi : [\spinc, A,s,\phi] \mapsto [\spinc, A, s\phi].
\]
This blow-up is a Banach manifold with boundary: the boundary consists
of points with $s=0$ (we call these \emph{reducible}), and the
restriction of $\pi$ to the boundary is a map
\[
           \pi : \partial\bonf^{\sigma}(Y) \to \bonf^{\red}(Y)
\]
with fibers the projective spaces associated to the vector spaces
$L^{2}_{k-1/2}(S)$.

\subsection{The Chern-Simons-Dirac functional}

After choosing a preferred connection $A_{0}$ in a spinor bundle $S$ for
each isomorphism class of $\SpinC$ structure, we can define the
Chern-Simons-Dirac functional $\CSD$ on $\conf(Y)$ by
\[
               \CSD(A,\Phi)=
               -\frac{1}{8} \int_{Y} (
                       A^{\ttr} -  A^{\ttr}_{0}) \wedge 
	               (F_{ A^{\ttr}} + F_{ A^{\ttr}_{0}}) + \frac{1}{2}\int_{Y}\langle 
	               D_{A}\Phi,\Phi\rangle \,d\mathrm{vol}.
\]
Here $A^{\ttr}$ is the associated connection in the line bundle
$\Lambda^{2}S$.  The formal gradient of $\CSD$ with respect to the
$L^{2}$ metric $\| \Phi \|^{2} + \frac{1}{4} \| A^{\ttr} -
A_{0}^{\ttr}\|^{2}$ is a ``vector field'' $\tilde\V$ on $\conf(Y)$ that is
invariant under the gauge group and orthogonal to its orbits. We use
quotation marks, because $\tilde\V$ is a section of the $L^{2}_{k-3/2}$
completion of tangent bundle.  
Away
from the reducible locus, $\tilde\V$
descends to give a vector field (in the same sense) $\V$ on $\bonf(Y)$. Pulling back
by $\pi$, we obtain a vector field $\V^{\sigma}$ on the interior of
the manifold-with-boundary $\bonf^{\sigma}(Y)$. This vector field
extends smoothly to the boundary, to give a section
\[
          \V^{\sigma} : \bonf^{\sigma}(Y) \to \cT_{k-3/2}(Y),
\]
where $\cT_{k-3/2}(Y)$ is the $L^{2}_{k-3/2}$ completion of $T\bonf^{\sigma}(Y)$.
This vector field is tangent to the boundary at
$\partial\bonf^{\sigma}(Y)$. The Floer groups $\Hto(Y)$, $\Hfrom(Y)$
and $\Hred(Y)$ will be defined using the Morse theory of the vector
field $\V^{\sigma}$ on $\bonf^{\sigma}(Y)$.

\subsubsection{Example}
\label{ex:ReducibleCritPoints}
Suppose that $b_{1}(Y)$ is zero. For each
$\SpinC$ structure $\spinc$, there is (up to isomorphism) a unique
connection $A$ in the associated spin bundle with $F_{A^{\ttr}}=0$,
and there is a corresponding zero of the vector field $\V$ at the
point $\alpha = [\spinc, A_{0}, 0]$ in $\bonf^{\red}(Y)$.  The vector field
$\V^{\sigma}$ has a zero at the point
$[\spinc,A_{0},0,\phi]$ in $\partial\bonf^{\sigma}(Y)$ precisely when
$\phi$ is a unit eigenvector of the Dirac operator $D_{A}$.  If the
spectrum of $D_{A}$ is simple (i.e.~no repeated eigenvalues), then
the set of zeros of $\V^{\sigma}$ in the projective space
$\pi^{-1}(\alpha)$ is a discrete set, with one point for each
eigenvalue.

\subsection{Four-manifolds}

Let $X$ be a compact oriented Riemannian $4$-manifold (possibly with
boundary), and write
$\bonf(X)$ for the space of isomorphism classes of triples $(\spinc,
A,\Phi)$, where $\spinc$ is a $\SpinC$ structure, $A$ is a $\SpinC$
connection of class $L^{2}_{k}$ and $\Phi$ is an $L^{2}_{k}$ section
of the associated half-spin bundle $S^{+}$.  As in the $3$-dimensional
case, we can form a blow-up $\bonf^{\sigma}(X)$ as the space of
isomorphism classes of quadruples $(\spinc, A, s, \phi)$, where $s\ge
0$ and $\|\phi\|_{L^{2}(X)} =1$. If $Y$ is a boundary component of
$X$, then there is a partially-defined
restriction map
\[
          \begin{gathered}
            r: \bonf^{\sigma}(X) \dashrightarrow \bonf^{\sigma}(Y)\\
          \end{gathered}
\]
whose domain of definition is the set of configurations
$[\spinc,A,s,\phi]$ on $X$ with $\phi|_{Y}$ non-zero. The map $r$
is given by 
\[
                [\spinc, A, s, \phi] \mapsto [\spinc|_{Y}, A|_{Y}, s/c,
                c\phi|_{Y}],
\]
where $1/c$ is the $L^{2}$ norm of $\phi|_{Y}$. (We have identified
the spin bundle $S$ on $Y$ with the restriction of $S^{+}.)$
When $X$ is cylinder $I \times Y$, with $I$ a compact interval,
we have a similar restriction map
\[
           r_{t} : \bonf^{\sigma}(I\times Y) \dashrightarrow \bonf^{\sigma}(Y)
\]
for each $t\in I$.

If $X$ is non-compact, and in particular if $X = \R\times Y$, then our
definition of the blow-up needs to be modified, because 
the $L^{2}$ norm of $\phi$ need not be finite. Instead, we define 
$\bonf^{\sigma}_{\loc}(X)$ as the space of isomorphism classes of
quadruples $[\spinc, A, \psi, \R^{+}\phi]$, where $A$ is a $\SpinC$
connection of class $L^{2}_{l,\loc}$, the set $\R^{+}\phi$ is 
the closed ray generated by a non-zero spinor $\phi$ in
$L^{2}_{k,\loc}(X;S^{+})$, and $\psi$ belongs to the ray.
(We write $\R^{+}$ for the
\emph{non-negative} reals.)
This is the usual way to define the blow
up of a vector space at $0$, without the use of a norm. The
configuration is \emph{reducible} if $\psi$ is zero.


\subsection{The four-dimensional equations}
When $X$ is compact, the Seiberg-Witten
monopole equations for a configuration $\gamma = [\spinc,A,s,\phi]$
in $\bonf^{\sigma}(X)$
are the equations
\begin{equation}\label{eq:4dsw}
	\begin{aligned}
		\frac{1}{2}\rho(F^{+}_{ A^{\ttr}}) - s^{2}(\phi\phi^{*})_{0} &= 0 \\
		D^{+}_{A} \phi &= 0,
	\end{aligned}
\end{equation}	
where 
$\rho: \Lambda^{+}(X) \to i\mathfrak{su}(S^{+})$ is Clifford
multiplication and $(\phi\phi^{*})_{0}$ denotes the traceless part of
this hermitian endomorphism of $S^{+}$. When $X$ is non-compact, we
can write down essentially the same equations using the ``norm-free''
definition of the blow-up, $\bonf^{\sigma}_{\loc}(X)$.  In either
case, we write
these equations as
\[
             \SW(\gamma) = 0.
\]
In the compact case, we write
\[
          M(X) \subset \bonf^{\sigma}(X)
\]
for the set of solutions. We draw attention to the non-compact case by
writing $M_{\loc}(X) \subset \bonf^{\sigma}_{\loc}(X)$.

Take $X$ to be the cylinder $\R\times Y$, and suppose that $\gamma =
[\spinc, A, \psi,\R^{+}\phi]$ is an element of $M_{\loc}(\R\times Y)$. A
unique continuation result implies that the restriction of $\phi$ to
each slice $\{t\}\times Y$ is non-zero;
so there is a well-defined restriction
\[
           \check\gamma(t) = r_{t}(\gamma) \in \bonf^{\sigma}(Y)
\]
for all $t$. We have the following relation between the equations
$\SW(\gamma) = 0$ and the vector field $\V^{\sigma}$:

\begin{lemma}
\label{lem:SW-is-Flow}
    If $\gamma$ is in $M_{\loc}(\R\times Y)$, then the corresponding
    path $\check\gamma$ is a smooth path in the Banach
    manifold-with-boundary $\bonf^{\sigma}(Y)$ satisfying
    \[
            \frac{d}{dt} \check\gamma(t) = -\V^{\sigma}.
    \]
    Every smooth path $\check\gamma$ satisfying the above condition
    arises from some element of $M_{\loc}(\R\times Y)$ in this way.
\end{lemma}

We should note at this point that our sign convention is such that the
$4$-dimensional Dirac operator $D^{+}_{A}$ on the cylinder $\R\times
Y$, for a connection $A$ pulled back from $Y$, is equivalent to the
equation
\[
        \frac{d}{dt} \phi + D_{A}\phi = 0 
\]
for a time-dependent section of the spin bundle $S\to Y$.

Next we define the moduli spaces that we will use to construct the
Floer groups.

\begin{defn}
    Let $\fa$ and $\fb$ be two zeros of the vector field $\V^{\sigma}$
    in the blow-up $\bonf^{\sigma}(Y)$. We write $M(\fa, \fb)$ for the
    set of solutions $\gamma \in M_{\loc}(\R\times Y)$ such that the
    corresponding path $\check\gamma(t)$ is asymptotic to $\fa$ as
    $t\to-\infty$ and to $\fb$ as $t\to+\infty$.
\end{defn}

Let $W : Y_{0} \to Y_{1}$ be an oriented cobordism, and suppose the
metric on $W$ is cylindrical in collars of the two boundary
components. Let $W^{*}$ be the cylindrical-end manifold obtained by
attaching cylinders $\R^{-} \times Y_{0}$ and $\R^{+}\times Y_{1}$.
From a solution $\gamma$ in $M_{\loc}(W^{*})$, we obtain paths
$\check\gamma_{0}: \R^{-} \to \bonf^{\sigma}(Y_{0})$ and
$\check\gamma_{1}:\R^{+} \to \bonf^{\sigma}(Y_{1})$. The following
moduli spaces will be used to construct the maps on the Floer groups
arising from the cobordism $W$:

\begin{defn}
\label{def:W-moduli}
    Let $\fa$ and $\fb$ be zeros of the vector field $\V^{\sigma}$
    in $\bonf^{\sigma}(Y_{0})$ and $\bonf^{\sigma}(Y_{1})$
    respectively. We write $M(\fa,W^{*}, \fb)$ for the
    set of solutions $\gamma \in M_{\loc}(W^{*})$ such that the
    corresponding paths $\check\gamma_{0}(t)$ and
    $\check\gamma_{1}(t)$ are asymptotic to $\fa$ and $\fb$ as
    $t\to-\infty$ and  $t\to+\infty$ respectively.
\end{defn}

\subsubsection{Example}
\label{ex:ReducibleFlows}
In example \ref{ex:ReducibleCritPoints}, suppose the spectrum is
simple, let $\fa_{\lambda} \in \partial\bonf^{\sigma}(Y)$ be the critical point corresponding to
the eigenvalue $\lambda$, and  let $\phi_{\lambda}$  be a
corresponding eigenvector of $D_{A_{0}}$. Then the reducible locus
$M^{\red}(\fa_{\lambda},\fa_{\mu})$ in the moduli space
$M(\fa_{\lambda},\fa_{\mu})$ is the quotient by $\C^{*}$ of the set
of solutions $\phi$ to the Dirac
equation
\[
            \frac{d}{dt} \phi + D_{A_{0}}\phi = 0 
\]
on the cylinder, with asymptotics
\[
         \phi \sim
         \begin{cases}
         C_{0} e^{-\lambda t}\phi_{\lambda}, & \text{as
         $t\to-\infty$,} \\
         C_{1} e^{-\mu t}\phi_{\mu}, & \text{as
         $t\to+\infty$} ,
         \end{cases} 
\]
for some non-zero constants $C_{0}$, $C_{1}$.

\subsection{Transversality and perturbations}

Let $\fa\in\bonf^{\sigma}(Y)$ be a zero of $\V^{\sigma}$.  The
derivative of the vector field at this point is a Fredholm operator on
Sobolev completions of the tangent space,
\[
           \Deriv_{\fa}\V^{\sigma} : \cT_{k-1/2}(Y)_{\fa} \to
               \cT_{k-3/2}(Y)_{\fa} .
\]
Because of the blow-up, this operator is not symmetric (for any simple
choice of inner product on the tangent space); but its spectrum is
real and discrete. We say that $\fa$ is \emph{non-degenerate} as a zero
of $\V^{\sigma}$ if $0$ is not in the spectrum.
If $\fa$ is a non-degenerate zero, then it is isolated, and we can
decompose the tangent space as
\[
          \cT_{k-1/2} = \K^{+}_{\fa} \oplus \K^{-}_{\fa},
\]
where $\K^{\pm}_{\fa}$ and $\K^{-}_{\fa}$ is the closures of the sum of
the generalized eigenvectors belonging to positive (respectively,
negative) eigenvalues. The \emph{stable manifold} of $\fa$ is the set
\[
        \stable_{\fa} = \{\, r_{0}(\gamma) \mid \gamma\in M_{\loc}
        (\R^{-}\times Y),\ \lim_{t\to+\infty}r_{t}(\gamma) = \fa\,\}.
\]
The unstable manifold $\unstable_{\fa}$ is defined similarly. If $\fa$
is non-degenerate, these are locally closed Banach submanifolds of
$\bonf^{\sigma}(Y)$ (possibly with boundary), and their tangent spaces at $\fa$ are the spaces
$\K^{+}_{\fa}$ and $\K^{-}_{\fa}$ respectively. Via the map $\gamma\mapsto\gamma_{0}$, we can
identify $M(\fa,\fb)$ with the intersection
\[
             M(\fa,\fb) = \stable_{\fa} \cap \unstable_{\fa}.
\]

In general, there is no reason to expect that the zeros are all
non-degenerate. (In particular, if $b_{1}(Y)$ is non-zero then
the reducible critical points are never isolated.)
To achieve non-degeneracy we perturb the equations, replacing the
Chern-Simons-Dirac functional $\CSD$ by $\CSD + f$, where $f$ 
belongs to a suitable class $\Pert(Y)$ of gauge-invariant functions on
$\conf(Y)$. We write $\tilde\pertY$ for the gradient of $f$ on
$\conf(Y)$, and $\pertY^{\sigma}$ for the resulting vector field on
the blow-up. Instead of the flow equation of
Lemma~\ref{lem:SW-is-Flow}, we now look (formally) at the equation
\[
    \frac{d}{dt}\check \gamma(t) = - \V^{\sigma} -\pertY^{\sigma}.
\]
Solutions of this perturbed flow equation correspond to solutions
$\gamma\in \bonf^{\sigma}(\R\times Y)$ of an equation
$\SW_{\pertY}(\gamma) = 0$ on the $4$-dimensional cylinder.
We do not define the class of perturbations $\Pert(Y)$ here (see
\cite{KMBook}).

The first important fact is that we can choose a perturbation $f$ from
the class $\Pert(Y)$ so that all the zeros of $\V^{\sigma} +
\pertY^{\sigma}$ are non-degenerate. From this point on we suppose
that such a perturbation is chosen. We continue to write $M(\fa,\fb)$
for the moduli spaces, $\stable_{\fa}$ and $\unstable_{\fa}$ for the
stable and unstable manifolds, and so on, without mention of the
perturbation. The irreducible zeros will
be a finite set; but as in Example~\ref{ex:ReducibleCritPoints}, the
number of reducible critical points will be infinite.  In general,
there is one reducible critical point $\fa_{\lambda}$ in the blow-up
for each pair $(\alpha,\lambda)$, where $\alpha=[\spinc,A,0]$ is a zero
of the restriction of $\V + \pertY$ to $\bonf^{\red}(Y)$, and
$\lambda$ is an eigenvalue of a perturbed Dirac operator
$D_{A,\pertY}$. The point $\fa_{\lambda}$ is given by
$[\spinc,A,0,\phi_{\lambda}]$, where $\phi_{\lambda}$ is a
corresponding eigenvector, just as in the example.

\begin{defn}
    We say that a reducible critical point $\fa \in \partial\bonf^{\sigma}(Y)$
    is boundary-stable if the normal vector to the boundary at $\fa$ belongs to
    $\K^{+}_{\fa}$. We say $\fa$  is boundary-unstable if the normal
    vector belongs to
    $\K^{-}_{\fa}$.
\end{defn}

In our description above, the critical point $\fa_{\lambda}$ is
boundary-stable if $\lambda>0$ and boundary-unstable if $\lambda<0$.
If $\fa$ is boundary-stable, then $\stable_{\fa}$ is a
manifold-with-boundary, and $\partial\stable_{\fa}$ is the reducible
locus $\stable^{\red}_{\fa}$. The unstable manifold $\unstable_{\fa}$
is then contained in $\partial\bonf^{\sigma}(Y)$. If $\fa$ is
boundary-unstable, then $\unstable_{\fa}$ is a manifold-with-boundary,
while $\stable_{\fa}$ is contained in the $\partial\bonf^{\sigma}(Y)$.

The Morse-Smale condition for the flow of the vector field
$\V^{\sigma}+\pertY^{\sigma}$ would ask that the intersection
$\stable_{\fa}\cap\unstable_{\fa}$ is a transverse intersection of
Banach submanifolds in $\bonf^{\sigma}(Y)$, for every pair of critical
points. We cannot demand this condition, because if $\fa$ is
boundary-stable and $\fb$ is boundary-unstable, then $\unstable_{\fa}$
and $\stable_{\fb}$ are both contained in
$\partial\bonf^{\sigma}(Y)$. In this special case, the best we can ask
is that the intersection be transverse in the boundary.

\begin{defn}
    We say that the moduli space $M(\fa,\fb)$ is boundary-obstructed
    if $\fa$ and $\fb$ are both reducible, $\fa$ is boundary-stable
    and $\fb$ is boundary-unstable.
\end{defn}

\begin{defn}
   We say that a moduli space $M(\fa, \fb)$ is regular if the
   intersection $\stable_{\fa}\cap\unstable_{\fb}$ is transverse,
   either as an intersection in the Banach manifold-with-boundary
   $\bonf^{\sigma}(Y)$ or (in the boundary-obstructed case) as an
   intersection in $\partial\bonf^{\sigma}(Y)$.  We say the
   perturbation is regular if:
   \begin{enumerate}
   \item
   all the zeros of $\V^{\sigma}  +
   \pertY^{\sigma}$ are non-degenerate;
   \item
   all the moduli spaces are regular; and
   \item
   there are no reducible critical points in the components
   $\bonf^{\sigma}(Y,\spinc)$ belonging to $\SpinC$ structures
   $\spinc$ with $c_{1}(\spinc)$ non-torsion.
   \end{enumerate}
\end{defn}

The class $\Pert(Y)$ is large enough to contain regular perturbations,
and we suppose henceforth that we have chosen a perturbation of this
sort. The moduli spaces $M(\fa,\fb)$ will be either manifolds or
manifolds-with-boundary, and the latter occurs only if $\fa$ is
boundary-unstable and $\fb$ is boundary-stable.  We write
$M^{\red}(\fa,\fb)$ for the reducible configurations in the moduli
space

\subsubsection{Remark}
\label{rem:NoIrreducibles}
The moduli space $M(\fa,\fb)$ cannot contain
any irreducible elements if $\fa$ is boundary-stable or if $\fb$ is
boundary-unstable.

\medskip

We can decompose $M(\fa,\fb)$ according to the relative homotopy
classes of the paths $\check\gamma(t)$: we write
\[
         M(\fa,\fb) = \bigcup_{z} M_{z}(\fa,\fb),
\]
where the union is over all relative homotopy classes $z$ of paths
from $\fa$ to $\fb$.  For any points $\fa$ and $\fb$ and any relative
homotopy class $z$, we can define an integer $\gr_{z}(\fa,\fb)$ (as
the index of  a suitable Fredholm operator), so that
\[
       \dim M_{z}(\fa,\fb) =
       \begin{cases}
        \gr_{z}(\fa,\fb)+1, &\text{in the boundary-obstructed case,}\\
        \gr_{z}(\fa,\fb), &\text{otherwise,}
       \end{cases}
\]
whenever the moduli space is non-empty.
The quantity $\gr_{z}(\fa,\fb)$ is additive along composite paths. We
refer to $\gr_{z}(\fa,\fb)$ as the \emph{formal dimension} of the
moduli space $M_{z}(\fa,\fb)$.

Let $W:Y_{0}\to Y_{1}$ be a cobordism, and suppose $\pertY_{0}$ and
$\pertY_{1}$ are regular perturbations for the two $3$-manifolds. Form
the Riemannian manifold $W^{*}$ by attaching cylindrical ends as
before.  We perturb the equations $\SW(\gamma) = 0$ on the compact
manifold $W$ by a perturbation $\pertX$ that is supported in
cylindrical collar-neighborhoods of the boundary components. The term
perturbation $\pertX$ near the boundary component $Y_{i}$ is defined
by a $t$-dependent element of $\Pert(Y_{i})$, equal to $\pertY_{0}$ in
a smaller neighborhood of the boundary. We continue to denote the
solution set of the perturbed equations $\SW_{\pertX}(\gamma) = 0$ by
$M(W) \subset \bonf^{\sigma}(W)$. This is a Banach manifold with
boundary, and there is a restriction map
\[
         r_{0,1} : M(W) \to \bonf^{\sigma}(Y_{0}) \times
         \bonf^{\sigma}(Y_{1}).
\]
The cylindrical-end moduli space $M(\fa, W^{*},\fb)$ can be regarded
as the inverse image of $\unstable_{\fa} \times \stable_{\fb}$ under
$r_{0,1}$:
\[
\begin{CD}
    M(\fa, W^{*}, \fb) @>>>  \unstable_{\fa} \times \stable_{\fb} \\
    @VVV @VVV \\
     M(W) @>{r_{0,1}}>> \bonf^{\sigma}(Y_{0}) \times
         \bonf^{\sigma}(Y_{1}).
\end{CD}
\]
The moduli space $M(\fa, W^{*},\fb)$ is \emph{boundary-obstructed} if
$\fa$ is boundary-stable and $\fb$ is boundary-unstable.

\begin{defn}
    If $M(\fa, W^{*},\fb)$ is not boundary-obstructed, we say that
    the moduli space is regular if $r_{0,1}$ is transverse to
    $\unstable_{\fa} \times \stable_{\fb}$. In the boundary-obstructed
    case, $M(\fa, W^{*},\fb)$ consists entirely of reducibles, and we
    say that it is regular if the restriction
    \[
              r_{0,1}^{\red} : M^{\red}(W) \to
              \partial\bonf^{\sigma}(Y_{0}) \times
              \partial\bonf^{\sigma}(Y_{1})
    \]
    is transverse to $\unstable_{\fa} \times \stable_{\fb}$.
\end{defn}

One can always choose the perturbation $\pertX$ on $W$ so that the moduli
spaces $M(\fa, W^{*},\fb)$ are all regular. Each moduli space has a
decomposition
\[
             M(\fa, W^{*},\fb) = \bigcup_{z} M_{z}(\fa, W^{*},\fb)
\]
indexed by the connected components $z$ of the fiber
$r_{0,1}^{-1}(\fa,\fb)$ of the map
\[
     r_{0,1} : \bonf^{\sigma}(W) \dashrightarrow \bonf^{\sigma}(Y_{0}) \times
     \bonf^{\sigma}(Y_{1}).
\]
The set of these components is a principal homogeneous space for the
group $H^{2}(W, Y_{0}\cup Y_{1};\Z)$. We can define an integer
$\gr_{z}(\fa, W,\fb)$ which is additive for composite cobordisms, such
that the dimension of the non-empty moduli spaces is given by:
\[
       \dim M_{z}(\fa,W^{*},\fb) =
       \begin{cases}
        \gr_{z}(\fa,W,\fb)+1, &\text{in the boundary-obstructed case,}\\
        \gr_{z}(\fa,W,\fb), &\text{otherwise,}
       \end{cases}
\]

\subsection{Compactness}
We suppose now that a regular perturbation
$\pertY$ is fixed.  The moduli space $M(\fa,\fb)$ has an action of
$\R$, by translations of the cylinder $\R\times Y$.  We write
$\Mu(\fa,\fb)$ for the quotient by $\R$ of the non-constant
solutions:
\[
       \Mu(\fa,\fb) = \{ \, \gamma \in M(\fa,\fb) \mid
       \text{$r_{t}(\gamma)$ is non-constant}\,\} /\R.
\]
We refer to elements of $\Mu(\fa,\fb)$ as unparametrized trajectories.
The space $\Mu_{z}(\fa,\fb)$ has a compactification: the space of
broken (unparametrized) trajectories $\Mubk_{z}(\fa,\fb)$. This space
is the union of all products
\begin{equation}\label{eq:BrokenStratum}
            \Mu_{z_{1}}(\fa_{0},\fa_{1}) \times \dots \times
            \Mu_{z_{l}}(\fa_{l-1}, \fa_{l}),
\end{equation}
where $\fa_{0}=\fa$, $\fa_{l}=\fb$ and the composite of the paths
$z_{i}$ is $z$.

Because of the presence of boundary-obstructed trajectories, the
enumeration of the strata that contribute to the compactification is
more complicated than it would be for a Morse-Smale flow. For example:

\begin{lemma}
\label{lem:BoundaryStrata}
    If $\Mu_{z}(\fa,\fb)$ is zero-dimensional, then it is compact. If
    $\Mu_{z}(\fa,\fb)$ is one-dimensional and contains irreducible
    trajectories, then the non-empty products \eqref{eq:BrokenStratum}
    that contribute to the compactification $\Mubk_{z}(\fa,\fb)$ are
    of two types:
    \begin{enumerate}
        \item products $\Mu_{z_{1}}(\fa,\fa_{1})\times
        \Mu_{z_{2}}(\fa_{1},\fb)$ with two factors;

        \item products $\Mu_{z_{1}}(\fa,\fa_{1})\times
        \Mu_{z_{2}}(\fa_{1},\fa_{2}) \times\Mu_{z_{3}}(\fa_{2}, \fb)$
        with three factors, of which the middle one is
        boundary-obstructed.
    \end{enumerate}
\end{lemma}

The situation for the reducible parts of the moduli spaces is simpler.
If $\Mu^{\red}_{z}(\fa,\fb)$ has dimension one, then its compactification
involves only broken trajectories with two components,
\[
           \Mu^{\red}_{z_{1}}(\fa, \fa_{1}) \times
           \Mu^{\red}_{z_{2}}(\fa_{1}, \fb).
\]

In \cite{KMBook}, gluing theorems are proved that describe the
structure of the compactification $\Mubk_{z}(\fa,\fb)$ near a stratum
of the type \eqref{eq:BrokenStratum}. In the case of a $1$-dimensional
moduli space containing irreducibles (as in the lemma above),
the compactification is a $C^{0}$ manifold with boundary in a
neighborhood of the strata of the first type. At a point belonging to
a stratum of the second type (with three factors), the structure of
the compactification is more complicated: a neighborhood of such a
point can be embedded in the
positive quadrant $\R^{+}\times\R^{+}$ as the zero set of a continuous function
that is strictly positive on the positive $x$-axis, strictly
negative on the positive $y$-axis, and zero at the origin. We refer to
this structure (more general than a $1$-manifold with boundary) as a
codimension-1 $\delta$-structure. Despite the extra complication,
spaces with this structure share with compact $1$-manifolds the fact
that the number of boundary points is even:

\begin{lemma}
\label{lem:EvenEndpoints}
    Let $N = N_{1} \cup N_{0}$ be a compact space, containing an open
    subset $N_{1}$ that is a smooth $1$-manifold and a closed
    complement $N_{0}$ that is a finite set. Suppose $N$ has a
    codimension-1 $\delta$-structure in the neighborhood of each point
    of $N_{0}$. Then $|N_{0}|$ is even.
\end{lemma}

\subsubsection{Remark}
\label{rem:ExtraBoundary}
In the case that $\fa$ is boundary-unstable and $\fb$ is
boundary-stable, the space $\Mu_{z}(\fa,\fb)$ is already a
manifold-with-boundary before compactification: the boundary is
$\Mu^{\red}_{z}(\fa,\fb)$.

\medskip

The moduli spaces $M_{z}(\fa, W^{*},\fb)$ can be compactified in a
similar way. For example, if $M_{z}(\fa, W^{*},\fb)$ contains
irreducibles and is one-dimensional, then it has a compactification
obtained by adding strata that are products of either two or three
factors. Those involving two factors have one of the two possible
shapes
\begin{equation}
\label{eq:TwoFactorW}
\begin{gathered}
    \Mu_{z_{1}}(\fa, \fa_{1}) \times M_{z_{2}}(\fa_{1}, W^{*},\fb) \\
    M_{z_{1}}(\fa, W^{*}, \fb_{1}) \times \Mu_{z_{2}}(\fb_{1},\fb)
\end{gathered}
\end{equation}
where the $\fa$'s belong to $\bonf^{\sigma}(Y_{0})$ and the $\fb$'s
belong to $\bonf^{\sigma}(Y_{1})$.
Those involving three factors have one of the three possible shapes
\begin{equation}
\label{eq:ThreeFactorW}
 \begin{gathered}
    \Mu_{z_{1}}(\fa, \fa_{1}) \times \Mu_{z_{2}}(\fa_{1},\fa_{2})
    \times M_{z_{3}}(\fa_{2}, W^{*},\fb) \\
        \Mu_{z_{1}}(\fa, \fa_{1}) \times 
    \times M_{z_{2}}(\fa_{1}, W^{*},\fb_{1}) \times
     \Mu_{z_{3}}(\fb_{1},\fb) \\
    M_{z_{1}}(\fa, W^{*}, \fb_{1}) \times \Mu_{z_{2}}(\fb_{1},\fb_{2})
    \times \Mu_{z_{3}}(\fb_{2},\fb).
\end{gathered}
\end{equation}
In the case of three factors, the middle factor is
boundary-obstructed.
All these strata are finite sets, and the compactification $\Mbk(\fa,
W^{*}, \fb)$ has a codimension-1 $\delta$-structure at each point.

\subsection{Three Morse complexes}
\label{subsec:ThreeMorse}

Let $\Crit^{s}$, $\Crit^{u}$ and $\Crit^{o}$ denote the set of
critical points (zeros of $\V^{\sigma} + \pertY^{\sigma}$ in $\bonf^{\sigma}(Y)$) that
are boundary-stable, boundary-unstable, and irreducible respectively.
Let
\[
          C^{s}(Y),\quad C^{u}(Y), \quad C^{o}(Y)
\]
denote vector spaces over $\Field$, with bases $e_{\fa}$ indexed by
the elements $\fa$ of these three sets.  For every pair of critical
points $\fa$, $\fb$, we define
\[
             n_{z}(\fa,\fb) = 
             \begin{cases}
                | \Mu_{z}(\fa,\fb) | \bmod{2},
                &\text{if $\dim\Mu_{z}(\fa,\fb)=0$,}\\
                0, &\text{otherwise.}
             \end{cases}
\]
From these, we construct linear maps
\[
\begin{aligned}
    \doo : C^{o}(Y) &\to C^{o}(Y) \\
    \dos : C^{o}(Y) &\to C^{s}(Y) \\
    \duo : C^{u}(Y) &\to C^{o}(Y) \\
    \dIus : C^{u}(Y) &\to C^{s}(Y)
\end{aligned}
\]
by the formulae
\[
         \doo e_{\fa} = \sum_{\fb\in \Crit^{o}}\sum_{z}
      n_{z}(\fa,\fb) e_{\fb}, \quad (\fa\in \Crit^{o}),
\]
and so on.  The four maps correspond to the four cases in which a
space of trajectories can contain irreducibles: see
Remark~\ref{rem:NoIrreducibles} above. 

Along with the $n_{z}({\fa},{\fb})$, we define quantities $\bar
n_{z}({\fa}_{\fb})$ using the reducible parts of the moduli spaces:
\[
               \bar n_{z}(\fa,\fb) = 
             \begin{cases}
                | \Mu^{\red}_{z}(\fa,\fb) | \bmod{2},
                &\text{if $\dim\Mu^{\red}_{z}(\fa,\fb)=0$,}\\
                0, &\text{otherwise.}
             \end{cases}
\]
These are used similarly as the matrix entries of linear maps
\[
\begin{aligned}
    \dss : C^{s}(Y) &\to C^{s}(Y) \\
    \dsu : C^{s}(Y) &\to C^{u}(Y) \\
    \dus : C^{u}(Y) &\to C^{s}(Y) \\
    \duu : C^{u}(Y) &\to C^{u}(Y).
\end{aligned}
\]
Note that the maps $\dus$ and $\dIus$ are different. The
former counts reducible elements in zero-dimensional moduli spaces
$\Mu^{\red}_{z}(\fa, \fb)$, where the corresponding irreducible
moduli space $\Mu(\fa, \fb)$ will be $1$-dimensional.

\begin{lemma}\label{lem:FourIdentities}
    We have the following identities:
    \[
    \begin{gathered}
     \doo\doo + \duo\dsu\dos = 0 \\
    \dos\doo + \dss\dos + \dIus\dsu\dos = 0 \\
   \doo\duo + \duo\duu + \duo\dsu\dIus = 0 \\
    \dus + \dos\duo + \dss\dIus + \dIus\duu + \dIus\dsu\dIus =
       0. \\
     \end{gathered}
     \]
\end{lemma}

\begin{proof} All four identities are proved by enumerating the
end-points of all the
    $1$-dimensional moduli spaces $\Mubk_{z}(\fa,\fb)$ that contain
    irreducibles, using
    Lemma~\ref{lem:BoundaryStrata} and Lemma~\ref{lem:EvenEndpoints}.
    In the last identity of the four, the extra term $\dus$ is
    accounted for by  Remark~\ref{rem:ExtraBoundary}.
\end{proof}

Using the reducible parts of the moduli spaces, we obtain the simpler
result:

\begin{lemma}\label{lem:FourIdentities-red}
    We have the following identities:
   \[
   \begin{gathered}
    \dss\dss + \dus\dsu\ = 0\\
   \dss\dus + \dus\duu  = 0\\
   \duu\dsu + \dsu\dss  = 0\\
  \duu\duu + \dsu\dus  =
       0.
       \end{gathered}
       \]
    \qed
\end{lemma}

\begin{defn}
\label{def:Complexes}
    We construct three vector spaces with differentials,
    $(\Cto(Y), \dto)$, $(\Cfrom(Y), \dfrom)$ and
   $(\Cred(Y), \dred)$,
    by setting
    \[
\begin{aligned}
       \Cto(Y) &= C^{o}(Y) \oplus C^{s}(Y) \\
       \Cfrom(Y) &= C^{o}(Y) \oplus C^{u}(Y) \\
       \Cred(Y) &= C^{s}(Y) \oplus C^{u}(Y),
\end{aligned}
    \]
    and defining
    \[
\begin{aligned}
           \dto &= \begin{bmatrix}
                                 \doo &  \duo\dsu  \\
                                 \dos &  \dss + \dIus\dsu
                       \end{bmatrix}, &
          \dfrom &= \begin{bmatrix}
                                 \doo &  \duo \\
                                 \dsu \dos  &  \duu + \dsu\dIus
                       \end{bmatrix},   &                  
             \dred = \begin{bmatrix}
                                 \dss &  \dus \\
                                 \dsu   &  \duu
                       \end{bmatrix}.
\end{aligned}
    \]
\end{defn}

The proof that the differentials $\dto$, $\dfrom$ and $\dred$ each
have square zero follows by elementary manipulation of the identities
in the previous two lemmas. We define the Floer homology
groups
\[
        \Hto_{*}(Y), \quad \Hfrom_{*}(Y), \quad \Hred_{*}(Y)
\]
as the homology of the three complexes above.  Each of these is a sum
of subspaces contributed by the connected components
$\bonf^{\sigma}(Y,\spinc)$ of $\bonf^{\sigma}(Y)$, so that
\[
        \Hto_{*}(Y) = \bigoplus_{\spinc} \Hto_{*}(Y,\spinc)
\]
for example.  After choosing a base-point,  we can grade
the complex $\Cto(Y,\spinc)$ by $\Z/d\Z$, where $d\Z$ is the subgroup
of $\Z$ arising as the image of the map
\[
            z \mapsto \gr_{z}(\fa,\fa)
\]
from $\pi_{1}(\bonf^{\sigma}(Y,\spinc), \fa)$ to $\Z$. This image is
contained in $2\Z$ and coincides with the image of the map
\eqref{eq:c1-map}.  The $\bullet$ versions of the Floer groups are
obtained by completion, as explained in Section~\ref{subsec:Gradings}.

To motivate the formalism
a little, it may be helpful to say that the construction of these
complexes can also be carried out (with less technical difficulty) in
the case that we replace $\bonf^{\sigma}(Y)$ by a finite-dimensional
manifold with boundary, $(B,\partial B)$. In the finite-dimensional
case, the complexes compute respectively the ordinary homology groups,
\[
          H_{*}(B;\Field), \quad H_{*}(B,\partial B;\Field), \quad
          H_{*}(\partial B;\Field).
\]
The long exact sequence \eqref{eq:ijp} is analogous to the long exact
sequence of a pair $(B,\partial B)$. The maps $i_{*}$, $j_{*}$ and $p_{*}$ arise from
maps $i$, $j$, $p$ on the chain complexes of
Definition~\ref{def:Complexes}, given by the matrices
     \begin{align*}
         i &= \begin{bmatrix}
                0 &\duo \\
                1 & \dIus
              \end{bmatrix}, &
         j &= \begin{bmatrix}
                1 & 0 \\
                0 & \dsu
              \end{bmatrix}, &
         p &= \begin{bmatrix}
                \dos & \dIus \\
                0 & 1
              \end{bmatrix}  .
     \end{align*}
The exactness of the sequence is a formal consequence of the
identities.  

Up to canonical isomorphism, the Floer groups are
independent of the choice of metric and perturbation that are involved
in their construction. As in Floer's original argument
\cite{FloerOriginal}, this
independence follows from the more general construction of maps from
cobordisms, and their properties.

\subsection{Maps from cobordisms}
\label{subsec:MapsCobordisms}

Let $W:Y_{0}\to Y_{1}$ be a cobordism equipped with a Riemannian metric
and a regular perturbation $\pertX$ so that the moduli spaces
$M_{z}(\fa, W^{*},\fb)$ are regular. For each pair of critical points
$\fa$, $\fb$ belonging to $Y_{0}$ and $Y_{1}$ respectively, let
\[
            m_{z}(\fa,W,\fb) = 
             \begin{cases}
                | M_{z}(\fa,W^{*},\fb) | \bmod{2},
                &\text{if $\dim M_{z}(\fa,W^{*},\fb)=0$,}\\
                0, &\text{otherwise.}
             \end{cases}
\]
Define $\bar m_{z}(\fa,W,\fb) $ for reducible critical points similarly,
using $M^{\red}_{z}(\fa, W^{*},\fb)$. These provide the matrix entries
of eight linear maps
\[
\begin{aligned}
    \moo : C^{o}_{\bullet}(Y_{0}) &\to C^{o}_{\bullet}(Y_{1}) \\
    \mos : C^{o}_{\bullet}(Y_{0}) &\to C^{s}_{\bullet}(Y_{1}) \\
    \muo : C^{u}_{\bullet}(Y_{0}) &\to C^{o}_{\bullet}(Y_{1}) \\
     \mIus : C^{u}_{\bullet}(Y_{0}) &\to C^{s}_{\bullet}(Y_{1})
\end{aligned}
\]
and
\[
\begin{aligned}
    \mss : C^{s}_{\bullet}(Y_{0}) &\to C^{s}_{\bullet}(Y_{1}) \\
    \msu : C^{s}_{\bullet}(Y_{0}) &\to C^{u}_{\bullet}(Y_{1}) \\
    \mus : C^{u}_{\bullet}(Y_{0}) &\to C^{s}_{\bullet}(Y_{1}) \\
    \muu : C^{u}_{\bullet}(Y_{0}) &\to C^{u}_{\bullet}(Y_{1}),
\end{aligned}
\]
with definitions parallel to the those of the maps $\doo$ etc.~above.
For example,
\[
            \mss (e_{\fa}) = \sum_{\fb\in\Crit^{s}(Y_{1})} \sum_{z}
              \bar m_{z}(\fa,W,\fb)  e_{\fb}, \quad (\fa \in
              \Crit^{s}(Y_{0})).
\]
The bullets denote completion, which is necessary because, for a given
$\fa$,
there are infinitely many $\fb$
for which $\bar m^{\fa}_{\fb,z}$ may be non-zero for some $z$. (For a
given $\fa$ and $\fb$, only finitely many $z$ can contribute.)
Again, $\mIus$ and $\mus$ are different maps.
By enumerating the boundary points of $1$-dimensional moduli spaces
$\Mbk(\fa, W^{*},\fb)$ and appealing to Lemma~\ref{lem:EvenEndpoints},
we obtain identities involving these operators. For example, by
considering such moduli spaces for which the end-points $\fa$ and
$\fb$ are both irreducible, we obtain the identity
\begin{equation}\label{eq:One-m-Identity}
 \moo\doo  + \doo\moo
       +  \duo\dsu\mos+ \duo\msu\dos+\muo\dsu\dos =0.
\end{equation}
The five terms in this identity enumerate the boundary points of each
of the five types described in \eqref{eq:TwoFactorW} and \eqref{eq:ThreeFactorW}.

We combine these linear maps to define maps
\[
\begin{aligned}
\mto(W) : \Cto_{\bullet}(Y_{0}) \to \Cto_{\bullet}(Y_{1}),\\
     \mfrom(W) : \Cfrom_{\bullet}(Y_{0}) \to \Cfrom_{\bullet}(Y_{1}), \\
     \mred (W): \Cred_{\bullet}(Y_{0}) \to \Cred_{\bullet}(Y_{1}),
     \end{aligned}
\]
by the formulae
    \begin{gather*}
            \mto(W) = \begin{bmatrix}
                                 \moo &  \muo\dsu+ \duo\msu \\
                                 \mos&  \mss + \mIus\dsu +
                                 \dIus\msu
                       \end{bmatrix}, \\
             \mfrom(W) = \begin{bmatrix}
                                 \moo &  \muo \\
                                 \msu \dos+ \dsu\mos  &
                                 \muu + \msu\dIus+
                                 \dsu \mIus
                       \end{bmatrix},                       
     \end{gather*}
and
\[
  \mred(W) = \begin{bmatrix}
                                 \mss &  \mus \\
                                 \msu&  \muu
                       \end{bmatrix}.
\]
Identities such as \eqref{eq:One-m-Identity} supply the proof of:

\begin{prop}
    The maps $\mto(W)$, $\mfrom(W)$ and $\mred(W)$ are chain maps, and they
    commute with $i$, $j$ and $p$.
\end{prop}

We define $\Hto(W)$, $\Hfrom(W)$ and $\Hred(W)$  to be the maps on the
Floer homology groups arising from the chain maps $\mto(W)$, $\mfrom(W)$ and
$\mred(W)$.

\subsection{Families of metrics}

The chain maps $\mto(W)$ depend on a choice of Riemannian metric $g$
and perturbation $\pertX$ on $W$. Let $P$ be a smooth manifold,
perhaps with boundary, and let $g_{p}$ and $\pertX_{p}$ be a smooth
family of metrics and perturbations on $W$, for $p\in P$.  We suppose
that there are collar neighborhoods of the boundary components $Y_{0}$
and $Y_{1}$ on which all the $g_{p}$ are equal to the same fixed
cylindrical metrics and on which all the $\pertX_{p}$ agree with the
given regular perturbations $\pertY_{0}$ and $\pertY_{1}$. We can form
a \emph{parametrized} moduli space over $P$, as the union
\[
\begin{aligned}
M(\fa, W^{*},\fb)_{P} &= \bigcup_{p}\  \{ p \} \times M(\fa,
                  W^{*},\fb)_{p} \\ &\subset P \times
                  \bonf^{\sigma}_{\loc}(W^{*}).
                  \end{aligned}
\]
Regularity for such moduli spaces is defined as the transversality of
the map
\[
          r_{0,1} : M(W)_{P} \to \bonf^{\sigma}(Y_{0})\times
          \bonf^{\sigma}(Y_{1})
\]
to the submanifold $\unstable_{\fa}\times\stable_{\fb}$, with the
usual adaptation in the boundary-obstructed cases. If $P$ has boundary
$Q$, then we take regularity of $M(\fa, W^{*},\fb)_{P}$ to include
also the condition that $M(\fa, W^{*},\fb)_{Q}$ is regular.  Given any
family of metrics $g_{p}$ for $p\in P$, and any family of
perturbations $\pertX_{q}$ for $q\in Q$ such that the moduli
spaces $M(\fa, W^{*},\fb)_{Q}$ are regular, we can choose an extension of the family
$\pertX_{q}$ to all  of $P$ in such a way that all the moduli spaces
$M(\fa, W^{*},\fb)_{P}$  are regular also.

Now suppose that $P$ is compact, with boundary $Q$.  For each $\fa$
and $\fb$, define
\[
            m_{z}(\fa,W,\fb)_{P} = 
             \begin{cases}
                | M^{P}_{z}(\fa,W^{*},\fb) | \bmod{2},
                &\text{if $\dim M^{P}_{z}(\fa,W^{*},\fb)=0$,}\\
                0, &\text{otherwise.}
             \end{cases}
\]
Use the moduli spaces of reducible solutions to define
$\bar m_{z}(\fa,W,\fb)_{P}$ similarly. From these, we construct linear
maps
\[
       \mto(W)_{P} : \Cto_{\bullet}(Y_{0}) \to \Cto_{\bullet}(Y_{1})
\]
with companion maps $\mfrom(W)_{P}$ and $\mred(W)_{P}$.  If the
boundary $Q$ is empty, then these are chain maps, just as in the case
that $P$ is a point: the proof is by enumeration of the boundary
points in the compactifications of $1$-dimensional moduli spaces
$M_{z}(\fa, W^{*},\fb)_{P}$. (The compactification is constructed as the
parametrized union of the moduli spaces $M_{z}^{+}(\fa, W^{*},
\fb)_{p}$ over
$P$.)

If
$Q$ is non-empty, then the boundary of $M_{z}(\fa, W^{*}, \fb)_{P}$ has an additional
contribution, namely the moduli space $M_{z}(\fa, W^{*},\fb)_{Q}$.
Identities such as \eqref{eq:One-m-Identity} therefore have an
additional term: we have, for example (one of eight similar
identities),
\begin{equation}\label{eq:One-m-Identity-P}
 (\moo)^{P}\doo  + \doo(\moo)^{P}
       +  \duo\dsu(\mos)_{P}+ \duo(\msu)_{P}\dos+(\muo)_{P}\dsu\dos =
       (\moo)_{Q}.
\end{equation}
The maps $\mto(W)_{P}$ etc.~are no longer chain maps: instead, we have
\[
\begin{aligned}
        \dto \mto(W)_{P} + \mto(W)_{P} \dto &= \mto(W)_{Q} \\
\       \dfrom \mfrom(W)_{P} + \mfrom(W)_{P} \dfrom &= \mfrom(W)_{Q} \\
\       \dred \mred(W)_{P} + \mred(W)_{P} \dred &= \mred(W)_{Q} .
\end{aligned}
\]
Thus $\mto(W)_{Q}$ is chain-homotopic to zero, and $\mto(W)_{P}$
provides the chain-homotopy. If we take $P$ to be the interval $[0,1]$
and $Q$ to be the boundary $\{0,1\}$, we obtain:

\begin{cor}
    The chain maps $\mto(W)_{0}$ and $\mto(W)_{1}$ from
    $\Cto_{\bullet}(Y_{0})$ to $\Cto_{\bullet}(Y_{1})$, corresponding
    to two different choices of metric and regular perturbation on the
    interior of $W$, are chain homotopic, and therefore induce the
    same map on Floer homology. The same holds for the other two
    flavors.
\end{cor}

\subsection{Composing cobordisms}
\label{subsec:Composite-Cobordisms}

Let $W : Y_{0} \to Y_{2}$ be a composite cobordism
\[
             W: Y_{0}\stackrel{W_{0}}\longrightarrow Y_{1}
             \stackrel{W_{1}}\longrightarrow Y_{2}.
\]
Equip $W$ with a metric which is cylindrical near the two boundary
components as well as in a neighborhood of $Y_{1}\subset W$, and let
$\pertX$ be a perturbation that agrees with the regular perturbations
$\pertY_{i}$ near $Y_{i}$ for $i=0,1,2$.  For each $T\ge 0$, let
$W(T)\cong W$ be the Riemannian manifold obtained by cutting along
$Y_{1}$ and  inserting a cylinder
$[-T,T]\times Y_{1}$ with the product metric. We can form the
parametrized union
\begin{equation}\label{eq:Family-sans-infty}
               \bigcup_{T\ge 0} M_{z}( \fa, W(T)^{*}, \fb).
\end{equation}
Since the manifolds $W(T)$ are all copies of $W$ with varying metric,
this can be seen as an example of a parametrized moduli space
$M_{z}(\fa, W^{*}, \fb)_{P}$ of the sort we have been considering,
with $P\cong [0,\infty)$. We
add a fiber over $T=\infty$ by setting
\[
             W(\infty)^{*} = W_{0}^{*} \amalg W_{1}^{*},
\]
a disjoint union of the two cylindrical end manifolds. We define
\[
             M_{z}(\fa, W(\infty)^{*},\fb)
\]
to be the union of products
\begin{equation}\label{eq:T-infty-fiber}
         \bigcup_{\fc } \bigcup_{z_{0},z_{1}}
         M_{z_{0}}(\fa, W_{0}^{*}, \fc) \times M_{z_{1}}(\fc,
         W_{1}^{*}, \fb).
\end{equation}
The union is over all pairs of classes $z_{0}$, $z_{1}$ with composite
$z$ and all $\fc\in\Crit(Y_{1})$. Putting in this extra fiber, we
have a family
\[
         M_{z}(\fa, W^{*}, \fb)_{\bar{P}} = \bigcup_{T\in[0,\infty]} \ \{T\}
         \times M_{z}(\fa, W(T)^{*}, \fb),      
\]
parametrized by the space $\bar{P}\cong[0,\infty]$.
The moduli space just defined is a non-compact manifold with boundary. The
boundary consists of the union of the two fibers over $T=0$ and
$T=\infty$, together with the reducible locus $ M^{\red}_{z}(\fa,
W^{*}, \fb)_{\bar{P}}$ in the case that the moduli space contains both
reducibles and irreducibles. It is contained in a compact space
\[
         \Mbk_{z}(\fa, W^{*}, \fb)_{\bar{P}} = \bigcup_{T\in[0,\infty]} \ \{T\}
         \times \Mbk_{z}(\fa, W(T)^{*}, \fb),     
\]
where for $T=\infty$ a typical element of $\Mbk_{z}(\fa, W(T)^{*},
\fb)$ is a quintuple
\[
           (\gamma_{0}, \gamma_{01}, \gamma_{1},
           \gamma_{12},\gamma_{2}),
\]
where $\gamma_{i}$ is a broken trajectory for $Y_{i}$ (possibly with
zero components) and $\gamma_{01}$, $\gamma_{12}$ belong to the moduli
spaces of $W_{0}$ and $W_{1}$.

\begin{lemma}
\label{lem:BoundaryStrata-Composite}
    If $M_{z}(\fa,W^{*},\fb)_{P}$ is zero-dimensional, then it is compact. If
    $M_{z}(\fa,W^{*},\fb)_{P}$ is one-dimensional and contains irreducible
    trajectories, then the compactification
    $\Mbk_{z}(\fa,W^{*},\fb)_{\bar{P}}$ is a
    $1$-dimensional space with a codimension-1 $\delta$-structure at
    all boundary points. The
    boundary points are of the following types:
    \begin{enumerate}
        \item the fiber over $T=0$, namely the space $M_{z}(\fa,
        W^{*}, \fb)$ for $W(0)=W$;

        \item the fiber over $T=\infty$, namely the union of products
        \eqref{eq:T-infty-fiber};

        \item products of two factors, of one of the forms
 \begin{equation*}
\begin{gathered}
    \Mu_{z_{1}}(\fa, \fa_{1}) \times M_{z_{2}}(\fa_{1}, W^{*},\fb)_{P} \\
    M_{z_{1}}(\fa, W^{*}, \fb_{1})_{P} \times \Mu_{z_{2}}(\fb_{1},\fb)
\end{gathered}
\end{equation*}
       (cf.~\eqref{eq:TwoFactorW} above);

     \item products of three factors, of one of the forms
\begin{equation*}
 \begin{gathered}
    \Mu_{z_{1}}(\fa, \fa_{1}) \times \Mu_{z_{2}}(\fa_{1},\fa_{2})
    \times M_{z_{3}}(\fa_{2}, W^{*},\fb)_{P} \\
        \Mu_{z_{1}}(\fa, \fa_{1}) 
    \times M_{z_{2}}(\fa_{1}, W^{*},\fb_{1})_{P} \times
     \Mu_{z_{3}}(\fb_{1},\fb) \\
    M_{z_{1}}(\fa, W^{*}, \fb_{1})_{P} \times \Mu_{z_{2}}(\fb_{1},\fb_{2})
    \times \Mu_{z_{3}}(\fb_{2},\fb)
\end{gathered}
\end{equation*}
     (cf.~\eqref{eq:ThreeFactorW} above);

     \item 
     \label{item:ThreeFactorComposite}%
     parts of the fiber $M_{z}^{+}(\fa, W(\infty)^{*}, \fb)$
     over $T=\infty$ of one the forms
     \[
\begin{gathered}
            \Mu_{z_{1}}(\fa,\fa_{1})\times M_{z_{2}}(\fa_{1},
            W_{0}^{*}, \fc) \times M_{z_{2}}(\fc, W_{1}^{*},
            \fb) \\
     M_{z_{1}}(\fa,
            W_{0}^{*}, \fc_{1}) \times \Mu_{z_{2}}(\fc_{1}, 
            \fc) \times M_{z_{3}}(\fc, W_{1}^{*}, \fb) \\
             M_{z_{1}}(\fa,
            W_{0}^{*}, \fc) \times M_{z_{2}}(\fc, W_{1}^{*}, \fb_{1})
            \times \Mu_{z_{3}}(\fb_{1}, 
            \fb),
\end{gathered}
     \]
     where the middle factor is boundary-obstructed in each case;
     \item the reducible locus $M^{\red}_{z}(\fa, W^{*},\fb)_{P}$ in the
     case that the moduli space contains both irreducibles and
     reducibles (which requires $\fa$ to be boundary-unstable and
     $\fb$ to be boundary-stable).
    \end{enumerate}
\end{lemma}

Following a familiar pattern, we now count the elements in the
zero-dimensional moduli spaces, to obtain elements of $\Field$:
\[
            m_{z}(\fa,W,\fb)_{P} = 
             \begin{cases}
                | M_{z}(\fa,W^{*},\fb)_{P} | \bmod{2},
                &\text{if $\dim M_{z}(\fa,W^{*},\fb)_{P}=0$,}\\
                0, &\text{otherwise,}
             \end{cases}
\]
and
\[
            \bar m_{z}(\fa,W,\fb)_{P} = 
             \begin{cases}
                | M^{\red}_{z}(\fa,W^{*},\fb)_{P} | \bmod{2},
                &\text{if $\dim M^{\red}_{z}(\fa,W^{*},\fb)_{P}=0$,}\\
                0, &\text{otherwise.}
             \end{cases}
\]
These become the matrix entries of maps
\[
\begin{aligned}
\Koo : C^{o}_{\bullet}(Y_{0}) &\to C^{o}_{\bullet}(Y_{2}) \\
\Kos : C^{o}_{\bullet}(Y_{0}) &\to C^{s}_{\bullet}(Y_{2}) \\
\Kuo : C^{u}_{\bullet}(Y_{0}) &\to C^{o}_{\bullet}(Y_{2}) \\
\KIus : C^{u}_{\bullet}(Y_{0}) &\to C^{s}_{\bullet}(Y_{2}) \\
\end{aligned}
\]
and
\[
 \begin{aligned}
    \Kss : C^{s}_{\bullet}(Y_{0}) &\to C^{s}_{\bullet}(Y_{2}) \\
    \Ksu : C^{s}_{\bullet}(Y_{0}) &\to C^{u}_{\bullet}(Y_{2}) \\
    \Kus : C^{u}_{\bullet}(Y_{0}) &\to C^{s}_{\bullet}(Y_{2}) \\
    \Kuu : C^{u}_{\bullet}(Y_{0}) &\to C^{u}_{\bullet}(Y_{2}),
\end{aligned}
\]
just as we defined $\moo$ and its companions. From
Lemma~\ref{lem:BoundaryStrata-Composite} and
Lemma~\ref{lem:EvenEndpoints} we obtain identities involving these
operators, as usual. For example, as an operator $C^{o}(Y_{0}) \to
C^{o}(Y_{2})$, we have
\begin{multline*}
  \moo(W)  
      +  \moo(W_{1})\moo(W_{0})  + \Koo\doo  + \doo\Koo
      +\Kuo\dsu\dos  + \duo\Ksu\dos
      +  \duo\dsu\Kos
      \\+ \muo(W_{1})\msu(W_{0})\dos +
      \muo(W_{1})\dsu\mos(W_{0}) +
      \duo\msu(W_{1})\mos(W_{0}) = 0.
\end{multline*}
The ten terms in this identity correspond to the ten possibilities
listed in the first five cases of the lemma above. (The final
case of the lemma does not apply.)

We combine the pieces $\Koo$ etc.~to define a map
\[
         \Kto : \Cto_{\bullet}(Y_{0}) \to \Cto_{\bullet}(Y_{2})
\]
by the matrix
\[
\Kto = \begin{bmatrix}
                                 \Koo &  \Kuo\dsu
                                 +\muo(W_{1})\msu(W_{0})+\duo\Ksu  \\
                                 \Kos &  \Kss + \KIus\dsu+
                                 \mIus(W_{1})\msu(W_{0})+
                                 \dIus\Ksu
                       \end{bmatrix}.
\]

\begin{prop}
\label{prop:K-homotopy-Composite}
    We have the equality
    \[
             \dto\Kto + \Kto\dto = \mto(W_{1})\mto(W_{0}) + \mto(W)
    \]
    as maps $\Cto_{\bullet}(Y_{0}) \to \Cto_{\bullet}(Y_{2})$.  At the
    level of homology therefore, we have
    \[
                 \Hto(W_{1}) \circ \Hto(W_{0}) = \Hto(W).
    \]
    There are similar identities for the other two flavors of Floer
homology.
\end{prop}

\begin{proof}
    The chain-homotopy identity is a formal consequence of the
    ten-term identity above, together with its seven companions and
    the corresponding identities for the $\partial$ and $m$ operators.
\end{proof}

\subsection{The module structure}
\label{subsec:ModuleStructure}

We describe now a way to define the $\Field[U]$-modules structure on
Floer homology. A different and more general approach is taken in
\cite{KMBook}, but the result is the same, and the present version of
the definition (based on \cite{DonaldsonKronheimer}) is simpler to describe.
Let $W: Y_{0} \to Y_{1}$ be a cobordism, and let
$w_{1},\dots, w_{p} \in W$ be chosen points. Let $B_{1},\dots,B_{p}$
be standard ball neighborhoods of these points. The space
$\bonf^{\sigma}(B_{q})$  is a Hilbert manifold with boundary; and
because it arises as a free quotient by the gauge group $\G$ of
$L^{2}_{k+1}$ maps $u:
B_{q} \to S^{1}$, there is a natural line bundle $L_{q}$ on
$\bonf^{\sigma}(B_{q})$ associated to the homomorhpism $u \mapsto
u(w_{q})$ from $\G$ to $S^{1}$.

Because of unique continuation, there is a well-defined restriction
map
\[
        r_{q} : M(W) \to \bonf^{\sigma}(B_{q}),
\]
and hence also
\[
             r_{q} : M(\fa, W^{*},\fb) \to \bonf^{\sigma}(B_{q}),
\]
for all $\fa$ and $\fb$. Let $s_{q}$ be a smooth section of $L_{q}$,
and let $V_{q} \subset \bonf^{\sigma}(B_{q})$ be its zero set.
Omitting the restriction maps that are implied by our notation, we now
consider the moduli spaces
\[
          M_{z}(\fa, W^{*}, \{w_{1}, \dots , w_{p}\}, \fb)
          \subset M_{z}(\fa, W^{*}, \fb)
\]
defined as the intersection
\[
             M_{z}(\fa, W^{*}, \fb) \cap V_{1} \cap \cdots \cap V_{p}.
\]
We can choose the sections $s_{q}$ so that, for all $\fa$ and $\fb$,
their pull-backs of $s_{1},\dots,s_{q}$ to $M(\fa, W^{*},\fb)$ have
transverse zero sets. The above intersection is then a smooth
manifold.

We repeat verbatim the construction of the chain maps $\mto(W)$,
$\mfrom(W)$ and $\mred(W)$ from Section~\ref{subsec:MapsCobordisms},
but replacing $M_{z}(\fa, W^{*},\fb)$ by the lower-dimensional
moduli space
$M_{z}(\fa, W^{*}, \{w_{1}, \dots , w_{p}\}, \fb)$ throughout. In this
way, we construct maps that we temporarily denote by
\[
          \Hto(W, \{ w_{1}, \dots, w_{p}\})
          : \Hto_{\bullet}(Y_{0}) \to \Hto_{\bullet}(Y_{1}),
\]
with similar maps for the other two flavors. As a special case, we
define a map
\[
            U :  \Hto_{\bullet}(Y) \to \Hto_{\bullet}(Y)
\]
by taking $p=1$ and taking $W$ to be the cylinder $[0,1]\times Y$.
The proof of the composition law for composite cobordisms adapts to
prove that $U^{p}$ is equal to the map arising from the cylindrical
cobordism with $p$ base-points; and more generally,
\[
         \Hto(W, \{ w_{1}, \dots, w_{p}\}) = U^{p} \Hto(W),
\]
(a formula which then makes the notation $\Hto(W, \{ w_{1}, \dots, w_{p}\})$
obsolete).


\subsection{Local coefficients}
\label{subsec:TwistedSetup}

There is a variant of Floer homology, using local coefficients. We
continue to work over the field $\Field = \Z/2$, and we consider a
local system of $\Field$-vector spaces, $\Gamma$ on
$\bonf^{\sigma}(Y)$. This means that for each points $\fa$ in
$\bonf^{\sigma}(Y)$ we have a vector space $\Gamma_{\fa}$ over
$\Field$, and for each relative homotopy-class of paths $z$ from $\fa$
to $\fb$ we have an isomorphism
\[
             \Gamma(z) : \Gamma_{\fa} \to \Gamma_{\fb}.
\]
These should satisfy the composition law $\Gamma_{z} =
\Gamma_{z_{2}}\circ\Gamma_{z_{1}}$ for composite paths.
Given such a local system, and given as usual a Riemannian metric and
regular perturbation for $Y$, we introduce vector spaces
$C^{o}(Y;\Gamma)$, $C^{s}(Y;\Gamma)$ and $C^{u}(Y;\Gamma)$, defining
them as
\[
            \bigoplus_{\fa} \Gamma_{\fa},
\]
where the sum is over all critical points in $\bonf^{\sigma}(Y)$ that
are irreducible, boundary-stable or boundary-unstable respectively. We
define a map $\doo : C^{o}(Y;\Gamma) \to C^{o}(Y;\Gamma)$ 
by the formula
\[
       \doo(e) = \sum_{\fb\in \Crit^{o}} \sum_{z} n_{z}(\fa, \fb)
       \Gamma(z)(e),  \quad (e\in \Gamma_{\fa}),
\]
where $n_{z}(\fa, \fb)$ is defined as before. This map, along with its
companions $\dos$ etc., are then used to define the differential
\[
            \dto : \Cto(Y;\Gamma) \to \Cto(Y;\Gamma)
\]
for the complex $\Cto(Y;\Gamma) = C^{o}(Y;\Gamma) \oplus
C^{s}(Y;\Gamma)$.  Proceeding as before, we construct the Floer group
$\Hto_{\bullet}(Y;\Gamma)$, and also its companions
$\Hfrom_{\bullet}(Y;\Gamma)$ and $\Hred_{\bullet}(Y;\Gamma)$.

Let $W : Y_{0} \to Y_{1}$ now be a cobordism, and suppose local
systems $\Gamma_{i}$ are given on $Y_{i}$ for $i=0,1$. The
restriction maps
\[
        r_{i} : \bonf^{\sigma}(W) \dashrightarrow \bonf^{\sigma}(Y_{i}), \quad
        (i=1,2)
\]
are only partially defined, but the pull-backs $r_{i}^{*}(\Gamma_{i})$
provide well-defined local systems on $\bonf^{\sigma}(W)$. This is
because a local system $\Gamma$ on $\bonf^{\sigma}(Y)$ is, in a
canonical way, the pull-back of a local system on $\bonf(Y)$, and the
restriction maps to $\bonf(Y_{i})$ are everywhere defined.

\begin{defn}
    A $W$-morphism from the local system $\Gamma_{0}$ on
    $\bonf^{\sigma}(Y_{0})$ to the local system $\Gamma_{1}$ on
    $\bonf^{\sigma}(Y_{1})$ is an isomorphism of local systems,
    \[
                  \Gamma_{W} : r_{0}^{*}(\Gamma_{0}) \to
                  r_{1}^{*}(\Gamma_{1}).
    \]
\end{defn}

Given $\fa$ in $\bonf^{\sigma}(Y_{0})$ and $\fb$ in
$\bonf^{\sigma}(Y_{1})$, and given a choice $z$ of a connected component
in $r_{0,1}^{-1}(\fa,\fb)$, a $W$-morphism provides us with an
isomorphism
\[
     \Gamma_{W}(z) : \Gamma_{\fa} \to \Gamma_{\fb},
\]
which behaves as expected with respect to composition on either side
with paths in $\bonf^{\sigma}(Y_{i})$. We can use $\Gamma_{W}$ to
define maps
\[
       \moo : C^{o}(Y_{0};\Gamma_{0}) \to C^{o}(Y_{1};\Gamma_{1}),
\]
and so on, just as in Section~\ref{subsec:MapsCobordisms} above: for
example, we define
\[
            \moo (e) = \sum_{\fb\in\Crit^{o}(Y_{1})} \sum_{z}
              m_{z}(\fa,W,\fb)\,  \Gamma_{W}(z)(e), \quad (e \in
              \Gamma_{0,\fa}).
\]
The result is a map
\[
          \Hto(W;\Gamma_{W}) : \Hto_{\bullet}(Y_{0};\Gamma_{0}) \to
          \Hto_{\bullet}(Y_{1};\Gamma_{1}),
\]
with companion maps on $\Hfrom$ and $\Hred$.  The proof of
independence of the choice of metric and perturbation on $W$, and the
proof of the composition law (with the obvious notion of composition
of $W$-morphisms), carry over with straightforward modifications.

\subsubsection{Support of local systems}
\label{subsubsec:support}
Let $Y$ be a three-manifold, and fix an open subset $M\subset
Y$. There is a partially-defined restriction map $\rho_M \colon \bonf^\sigma(Y)
\dashrightarrow
\bonf^\sigma(M)$.  A local system $\Gamma$ over $\bonf^\sigma(M)$ induces a local
system $\rho_M^*(\Gamma)$ over $\bonf^\sigma(Y)$ by pull-back. (The
fact that the $\rho_{M}$ is only partially defined is again of no
conequence, as above.)
\begin{defn}
\label{defn:Support}
A local system over $\bonf^\sigma(Y)$ which is obtained as the pull-back of
one over $\bonf^\sigma(M)$ is said to be supported on $M$.
\end{defn}

Similarly, suppose we have a cobordism $W\colon Y_0\longrightarrow
Y_1$, equipped with an open set $B\subset W$, and let $M_0=B\cap Y_0$
and $M_1=B\cap Y_1$. Let $\Gamma_0$ and $\Gamma_1$ be local systems on
the $M_i$. Again, we have partially-defined restriction maps
\[
\rho_0\colon \bonf^\sigma(B) \dashrightarrow \bonf^\sigma(M_{i}),
\quad (i=1,2)
\]
which induce well-defined local systems $\rho_i^*(\Gamma_i)$.  over
$\bonf^\sigma(B)$. A $B$-morphism of local systems is an isomorphism
$$\Gamma_B\colon \rho_0^*(\Gamma_0) \longrightarrow \rho_1^{*}(\Gamma_1)$$
of local systems over $\bonf^\sigma(B)$. Using the 
restriction map
$$\rho\colon \bonf^\sigma(W) \dashrightarrow \bonf^\sigma(B),$$
we can pull back a $B$-morphism $\Gamma_B$ of local systems to 
obtain a $W$-morphism
\[
\rho_B^*(\Gamma_B) \colon r_0^*(\rho_0^*(\Gamma_0)) \longrightarrow
r_1^*(\rho_1^*(\Gamma_1)).
\]
Such a $W$-morphism is said to be {\em supported on $B$}.

\subsubsection{Example}
\label{ex:LocalSystem}
Let $\eta$ be a $C^{\infty}$ singular
$1$-cycle in $Y$ with real coefficients. Given a relative homotopy
class of paths $z$ from $\fa$ to $\fb$ in $\bonf^{\sigma}(Y)$, let us
choose a representative path $\tilde z$, and let $[A_{\tilde z}, s,
\phi_{\tilde z}]$ be the corresponding element of $\bonf^{\sigma}([0,1]\times
Y)$. Define 
\[
           f_{\eta}(z) = (i/2\pi) \int_{[0,1]\times \eta}
           F_{A^{\ttr}_{\tilde z}} .
\]
This depends only on $\eta$ and $z$.  

Let $\Kfield$ be an integral domain of characteristic
$2$, and let
\[
        \mu : \R \to \Kfield^{\times}
\]
be a homomorphism from the additive group $\R$ to the multiplicative
group of units in $\Kfield$. 
We can construct a local
system $\Gamma_{\eta}$ on $\bonf^{\sigma}(Y)$ by declaring that
$\Gamma_{\eta,\fa}$ is $\Kfield$ for all $\fa$, and that
\[
       \Gamma_{\eta}(z) : \Gamma_{\eta,\fa} \to \Gamma_{\eta,\fb}
\]
is multiplication by the unit $\mu (f_{\eta}(z))$ in
$\Kfield^{\times}$. For definiteness, we henceforth take $\Kfield$ to
be the field of fractions of the group ring $\Field[\R]$, and $\mu$ to
be the natural inclusion
\[
          \mu : \R \to \Field[\R]^{\times} \subset \Kfield^{\times}.
\]

Now let $W: Y_{0}\to Y_{1}$ be a cobordism, let $\eta_{0}$, $\eta_{1}$
be $1$-cycles as above, and let $\Gamma_{\eta_{i}}$ be the
corresponding local systems.  Suppose we are given a $C^{\infty}$
singular $2$-chain $\nu$ in $W$, with
\[
           \partial \nu = \eta_{1} - \eta_{0}.
\]
Given a component $z$ in $r^{-1}_{0,1}(\fa, \fb)$, we choose a
representative $[A_{z}, s, \phi_{z}]$ in $\bonf^{\sigma}(W)$, and we
extend our notation above by setting
\[
            f_{\nu}(z) =(i/2\pi) \int_{\nu} F_{A_{z}}.
\]
We can define a $W$-morphism $\Gamma_{W,\nu} : \Gamma_{\eta_{0}} \to
\Gamma_{\eta_{1}}$ by specifying that the isomorphism
\[
             \Gamma_{W,\nu}(z) : \Gamma_{\eta_{0},\fa} \to
             \Gamma_{\eta_{1},\fb}
\]
is given by multiplication by $\mu(f_{\nu}(z))$. There are
corresponding maps
\[
\begin{aligned}
    \Hto(W;\Gamma_{W,\nu}) : \Hto_{\bullet}(Y_{0};\Gamma_{\eta_{0}}) &\to
          \Hto_{\bullet}(Y_{1};\Gamma_{\eta_{1}}) \\
    \Hfrom(W;\Gamma_{W,\nu}) : \Hfrom_{\bullet}(Y_{0};\Gamma_{\eta_{0}}) &\to
          \Hfrom_{\bullet}(Y_{1};\Gamma_{\eta_{1}}) \\
    \Hred(W;\Gamma_{W,\nu}) : \Hred_{\bullet}(Y_{0};\Gamma_{\eta_{0}}) &\to
          \Hred_{\bullet}(Y_{1};\Gamma_{\eta_{1}}) \\
\end{aligned}
\]

We can consider these constructions as defining functors on an
extension of our cobordism category. We have a category whose objects
are pairs $(Y,\eta)$, where $Y$ is a $3$-manifold (compact, connected
and oriented as usual) and $\eta$ is a $C^{\infty}$ singular $1$-cycle
with real coefficients. The morphisms are diffeomorphism classes of
pairs $(W,\nu)$, where $\nu$ is a $2$-cycle and $\partial \nu =
\eta_{1} - \eta_{0}$.

From the definitions, it follows that if $\tilde{\nu} = \nu +
\partial\theta$ for some $C^{\infty}$ $3$-chain $\theta$, then
$\Gamma_{W,\nu}$ and $\Gamma_{W,\tilde\nu}$ are equal. As a
consequence, there are isomorphisms (for example)
\[
        \Hto_{\bullet}(Y;\Gamma_{\eta}) \cong
        \Hto_{\bullet}(Y;\Gamma_{\eta'})
\]
whenever $[\eta] = [\eta']$ in $H_{1}(Y;\R)$.  However, to specify a
particular isomorphism, one must express $\eta-\eta'$ as a boundary.
Indeed, suppose that $\partial \nu_1 = \eta-\eta'=\partial \nu_2$,
then the two isomorphisms differ  by the 
automorphism which on $\Hto_{\bullet}(Y,\spinct;\Gamma_{\eta})$ 
is given by multiplication by $\mu({\langle
c_1(\spinct),[\nu-\nu']\rangle})$. In particular, when $Y$ is a rational
homology three-sphere and $\eta$ is a cycle, then there is a canonical identification
$$\Hto(Y;\Gamma_\cycle)\cong \Hto(Y)\otimes\Kfield.$$

If the $1$-cycle $\eta$ is contained in $M\subset
Y$, then the local coefficient system
$\Gamma_{\eta}$ is supported on $M$, in the sense of the definition
above. Moreover, suppose $W\colon
Y_{0}\to Y_{1}$ is a cobordism, $B\subset W$ is an open set with
$M_{i}=B\cap Y_{i}$, and $\eta_{i}$ are $1$-cycles in $M_{i}$. Let
$\nu$ be a $2$-cycle with $\partial\nu = \eta_{1}-\eta_{0}$. Then
the $W$-morphism $\Gamma_{W,\nu} :\Gamma_{\eta_{0}} \to
\Gamma_{\eta_{1}}$  is supported on
$B$ if $\nu$ is contained in $B$.

We conclude this section by noting the following result. (The only
particular property of our choice of $\Kfield$ and $\mu$ that is used
here is the fact that $1-\mu(t)$ is a unit, for all non-zero $t$.)

\begin{lemma}
	\label{lemma:VanishTwist}
    If $[\eta]$ is non-zero in $H_{1}(Y;\R)$, then
    $\Hred_{\bullet}(Y;\Gamma_{\eta})=0$ and hence
    $j_{*}:\Hto_{\bullet}(Y;\Gamma_{\eta}) \to
    \Hfrom_{\bullet}(Y;\Gamma_{\eta})$ is an isomorphism.
\end{lemma}

\begin{proof}
    If $c_{1}(\spinc)$ is non-torsion, there are no reducible critical
    points, so $\Hred_{\bullet}(Y;\Gamma_{\eta})$ has contributions
    only from those $\spinc$ with torsion first Chern class.  For each
    such $\spinc$, the space of reducibles in $\bonf(Y,\spinc)$ has
    the homotopy type of the torus $H^{1}(Y;\R)/H^{1}(Y;\Z)$, for it
    deformation-retracts onto the torus $T$ of $\SpinC$ connections
    $A$ in $S$ with $F_{A^{\ttr}}=0$.  The lemma is a consequence of
    the fact that the ordinary homology group $H_{*}(T;
    \Gamma_{\eta})$ with local coefficients is zero. Details are given
    in \cite{KMBook}.
\end{proof}

\subsection{Duality and pairings}
\label{subsec:Duality}

Along with the chain complex $(\Cto_{*}(Y), \dto)$ and its companions,
we have the corresponding cochain complexes,
\[
    (\Cto^{*}(Y), \dto^{*}), \quad  (\Cfrom^{*}(Y), \dfrom^{*}), \quad
     (\Cred^{*}(Y), \dred^{*}),
\]
and the monopole Floer cohomology groups $\Hto^{*}(Y)$,
$\Hfrom^{*}(Y)$, $\Hred^{*}(Y)$. These are modules over $\Field[U]$,
with $U$ now acting with degree $2$. To form the $\bullet$ versions
$\Hto^{\bullet}(Y)$, we should now complete in the direction of
increasing degree, so that they again become modules over
$\Field[[U]]$. We can do the same with local coefficients, and we have
non-degenerate pairings of $\Kfield$-vector spaces
\begin{equation}\label{eq:HoCo-pairing}
       \Hto^{j}(Y;\Gamma_{-\eta}) \times \Hto_{j}(Y;\Gamma_{\eta}) \to
       \Kfield,
\end{equation}
for any $C^{\infty}$ real $1$-cycle $\eta$.

Let $-Y$ denote the oriented manifold obtained from $Y$ by reversing
the orientation. The spaces $\bonf^{\sigma}(Y)$ and
$\bonf^{\sigma}(-Y)$ can be canonically identified, though the change
of orientation changes the sign of the functional $\CSD$, and so
changes the vector field $\V^{\sigma}$ to $-\V^{\sigma}$. If $\pertY$
is a regular perturbation for $Y$, we can select $-\pertY$ as regular
perturbation for $-Y$.  The notion of boundary-stable and
boundary-unstable are interchanged when the vector field changes sign,
so we have identifications
\[
     \begin{aligned}
        \Crit^{o}(Y) &= \Crit^{o}(-Y) \\
        \Crit^{s}(Y) &= \Crit^{u}(-Y) \\
        \Crit^{u}(Y) &= \Crit^{s}(-Y).
     \end{aligned}
\]
The moduli space $M(\fa,\fb)$ for the cylinder $\R\times Y$ is the
same as the moduli space $M(\fb,\fa)$ for the cylinder $\R\times
(-Y)$. So the operator $\dos$ on $-Y$, for example, becomes
$(\duo)^{*}$ for $Y$.  In this way, the boundary map $\dfrom$ for $-Y$
becomes the operator $\dto^{*}$ for $Y$, and so on. Thus we obtain the
following proposition, which we state also for local coefficients.

\begin{prop}
    There are canonical isomorphisms
    \[
    \begin{aligned}
    \dual: \Hto_{j}(-Y) &\to \Hfrom^{j}(Y) \\
              \dual: \Hfrom_{j}(-Y) &\to \Hto^{j}(Y) \\
              \dual: \Hred_{j}(-Y) &\to \Hred^{j}(Y).
              \end{aligned}
    \]
    If $\eta$ is a real $1$-cycle in $Y$, and $\Gamma_{\eta}$ the
    corresponding local system with fiber $\Kfield$, then there are
    also isomorphisms
    \[
    \begin{aligned}
    \dual: \Hto_{j}(-Y;\Gamma_{\eta}) &\to \Hfrom^{j}(Y;\Gamma_{\eta}) \\
              \dual: \Hfrom_{j}(-Y;\Gamma_{\eta}) &\to \Hto^{j}(Y;\Gamma_{\eta}) \\
              \dual: \Hred_{j}(-Y;\Gamma_{\eta}) &\to \Hred^{j}(Y;\Gamma_{\eta}).
              \end{aligned}
              \]
    Here $j$ belongs to the grading set $J(Y)$, which we can identify
    canonically with $J(-Y)$, because the notion of an oriented
    $2$-plane field on $Y$ makes no reference to the orientation of
    the manifold.
\end{prop}

Note that the canonical identification $J(Y) = J(-Y)$ does not respect
the action of $\Z$: there is a sign change, so $j+n$ becomes $j-n$ if
the orientation of $Y$ is reversed.

If we combine the isomorphisms $\dual$ with the pairing
\eqref{eq:HoCo-pairing}, we obtain a non-degenerate pairing of vector
spaces over $\Field$,
\[
\langle\;\mathord{-}\;, \;\mathord{-}\;\rangle_{\dual}
        : \Hfrom_{*}(-Y) \times \Hto_{*}(Y) \to \Field.
\]
With local coefficients, there is a non-degenerate pairing of vector
spaces over $\Kfield$:
\[
\langle\;\mathord{-}\;, \;\mathord{-}\;\rangle_{\dual}
        : \Hfrom_{*}(-Y;\Gamma_{-\eta}) \times
        \Hto_{*}(Y;\Gamma_{\eta}) \to \Kfield.
\]

\subsection{Calculations for lens spaces}
\label{subsec:LensSpace-II}

Let $S^{3}_{p}(U)$ be the lens space obtained by integer surgery on
the unknot, for some $p>0$, and let $W(p) : S^{3}_{p}(U) \to S^{3}$ be
the surgery cobordism, as described in
Section~\ref{subsec:LensSpace-I}.  Corollary~\ref{cor:FroyS3p} states
that the map
\begin{equation}\label{eq:Wp-onto-again}
            \Hfrom(W(p), \spinc_{n,p}) :
             \Hfrom_{\bullet}(S^{3}_{p}(U), \tspin_{n,p})\to
             \Hfrom_{\bullet}(S^{3})
\end{equation}
is an isomorphism if $0\le n\le p$, an assertion which is equivalent
to the calculation of the Fr{\o}yshov
invariant of $(S^{3}_{p}(U), \tspin_{n,p})$ as given in
Proposition~\ref{prop:FroyS3}.   We shall now provide a proof of this
result.

Recall that we have injective maps
\[
\begin{gathered}
    p_{*} : \Hfrom_{\bullet}(S^{3}_{p}(U),
    \tspin_{n,p})\hookrightarrow
    \Hred_{\bullet}(S^{3}_{p}(U), \tspin_{n,p})\\
     p_{*} : \Hfrom_{\bullet}(S^{3})\hookrightarrow
    \Hred_{\bullet}(S^{3})
\end{gathered}
\]
because these $3$-manifolds have positive scalar curvature
for their standard metrics, c.f. Proposition~\ref{prop:PSC}. The
complex $\Cred_{*}(S^{3}_{p}(U))$ can be described using the material
of Example~\ref{ex:ReducibleCritPoints}. In $\bonf(S^{3}_{p}(U),
\tspin_{n,p})$, we have a unique critical point $\alpha =
[0,A_{n},0]$, with $A_{n}^{\ttr}$ flat. After a choice of small
perturbation, we can assume that the perturbation $D_{A_{n},\pertY}$
of the Dirac operator $D_{A_{n}}$ has simple spectrum; and there is
then one non-degenerate critical point $\fa_{\lambda}$ in the blow-up
$\bonf^{\sigma}(S^{3}_{p}(U), \tspin_{n,p})$ for each eigenvalue
$\lambda$ of $D_{A_{n},\pertY}$.  The differentials in the Floer
complexes are all zero, and we have
\[
\begin{aligned}
\Hred_{*}(S^{3}_{p}(U),\tspin_{n,p}) &=
       \Cred_{*}(S^{3}_{p}(U),\tspin_{n,p}) \\
       &=\bigoplus_{\lambda} \Field \, e_{\fa_{\lambda}}.
       \end{aligned}
\]
The image $p_{*} \Hfrom_{*}(S^{3}_{p}(U),\tspin_{n,p})$ is the
subspace
\[
\Cfrom_{*}(S^{3}_{p}(U),\tspin_{n,p}) 
       =\bigoplus_{\lambda < 0} \Field \, e_{\fa_{\lambda}},
\]
which is generated as an $\Field[U]$-module by the element
$e_{\fa_{\lambda_{-1}}}$, where $\lambda_{-1}$ is the first negative
eigenvalue. Let $\fb_{\mu_{-1}}$ similarly denote the critical point
in $\bonf^{\sigma}(S^{3})$ corresponding to the generator of
$p_{*}\Hfrom_{*}(S^{3})$; so $\mu_{-1}$ is the first negative
eigenvalue of the perturbed Dirac operator on $S^{3}$.

The assertion that the map \eqref{eq:Wp-onto-again} is surjective is
equivalent to the assertion that
\[
     \Hfrom(W(p), \spinc_{n,p})(e_{\fa_{\lambda_{-1}}}) =
      e_{\fb_{\mu_{-1}}}.
\]
The moduli space $M^{\red}_{\spinc_{n,p}}(\fa_{\lambda_{-1}}, W(p),
\fb_{\mu_{-1}})$ can be identified with the space of equivalence
classes of pairs $[A, \phi]$, where
\begin{itemize}
\item
$A$ is a $\SpinC$ connection for $\spinc_{n,p}$ on
$W(p)^{*}$, satisfying $F^{+}_{A^{\ttr}}=0$ and such that $A^{\ttr}$ is
asymptotic to a flat connection on both ends; and
\item
$\phi$ is a section of $S^{+}$ on $W(p)^{*}$ satisfying a
small-perturbation of the Dirac equation $D_{A,\pertX}^{+}\phi=0$ and
having asymptotics
\[
\begin{aligned}
        | \phi | &=  O( e^{-\lambda_{-1}t} ), & t&\to-\infty,\\
        | \phi | &=  O( e^{-\mu_{-1}t} ), & t&\to+\infty,
\end{aligned}
\]
on the two ends of $W(p)^{*}$.
\end{itemize}
(The equivalence relation is generated by the action of gauge
treansformations and scaling $\phi$ by non-zero complex scalars.)

The
connection $A$ satisfying the first condition is unique up to gauge
transformation. Given $A$, the spinors $\Phi$ satisfying the second
condition form an open subset of a projective space of complex
dimension $2d$, where $d$ is the $L^{2}$ index of the Dirac operator
$D^{+}_{A}$ on $W(p)^{*}$. The moduli space is a point if $d=0$. So
the surjectivity of \eqref{eq:Wp-onto-again} is eventually equivalent
to the next lemma.

\begin{lemma}
    Let $A$ be a $\SpinC$ connection for a $\SpinC$ structure
    $\spinc$ on $W(p)^{*}$, and suppose that
    that $F_{A^{\ttr}}$ is anti-self-dual and is
    asymptotically zero on both ends. Then  the index of the Fredholm
    operator \[ D^{+}_{A} : L^{2}_{1,A}(W(p)^{*},S^{+})) \to
  L^{2}(W(p)^{*}, S^{-})\] is zero if the pairing of $c_{1}(\spinc)$
  with a generator $h$ of $H_{2}(W(p);\Z)$ satisfies
  \[
               -p \le \langle c_{1}(\spinc),h\rangle \le p.
  \]
\end{lemma}

\begin{proof}
    If we change the orientation of $W(p)$, make a conformal change, and add two points, we
    obtain a K\"ahler orbifold $\bar{W}$, the weighted projective
    space obtained as the quotient of $\C^{3}\setminus \{0\}$ by the
    action of $\C^{*}$ with weights $(1,1,p)$.  Because of its
    conformal invariance, we can equivalently
    study the index of the Dirac operator on $\bar{W}$. Let the
    K\"ahler form $\omega$ be normalized so as to have integral $1$ on
    $h$, and let $L_{k}$ be the orbifold line bundle on $\bar{W}$ with
    curvature $c_{1}(L_{k}) = k[\omega]$. Let $\spinc_{0}$ be the canonical
    $\SpinC$ structure (which has $c_{1}(\spinc_{0})=-(p+2)[\omega]$)
    and let $\spinc_{k}$ be the $\SpinC$ structure on
    $\bar{W}$ obtained by tensoring $\spinc_{0}$ with $L_{k}$. Let
    $A_{k}$ be a $\SpinC$ connection for $\spinc_{k}$ compatible with
    the holomorphic structure.

    The kernel of
    $D^{+}_{A_{k}}$ is isomorphic to the orbifold Dolbeault cohomology
    group
    \[
             H^{0,0}(\bar{W}; L_{k}) \oplus  H^{0,2}(\bar{W}; L_{k}).
    \]
    If $k<0$, then $H^{0,0}$ vanishes because $L_{k}$ then has
    negative degree. If $k > -p-2$, then $H^{0,2}$ vanishes, by Serre
    duality. So the kernel of the Dirac operator is zero for $-p-2 < k
    < 0$; or equivalently $-p - 2 < c_{1}(\spinc_{k}) < p+2$. The
    cokernel is $H^{0,1}(\bar{W};L_{k})$, which vanishes for all $k$.
\end{proof}
\newcommand\ab{\mathrm{ab}}

\section{Proof of the surgery long exact sequence}
\label{sec:LongExactSequences}

\begin{figure}
\mbox{\vbox{\includegraphics{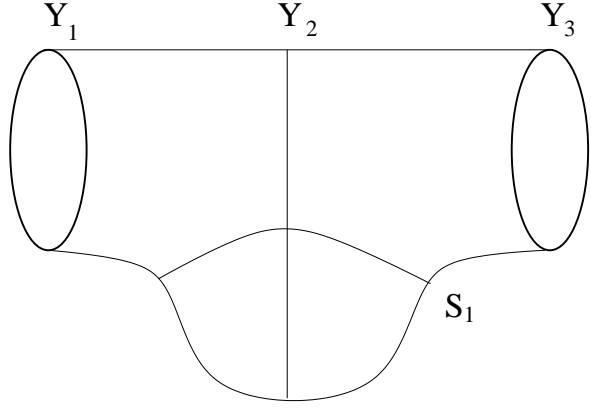}}}
\caption{\label{fig:Square}
{\bf{Breaking up the composite.}}
This indicates two ways of breaking up the composite cobordism $X_1$.
$Y_2$ separates $X_1$ into $W_1$ and $W_2$, while $S_1$
separates it into $B_1$ and $Z_1$.}
\end{figure}

We return to the notation of Theorem~\ref{thm:TwoHandleSurgeries-A}.
Before starting the proof in earnest, we explain the simple argument
which shows that
the composites 
$\Hto_{\bullet}(W_{n+1})\circ \Hto_{\bullet}(W_{n})$ are zero. We use
the composition law to equate this map to $\Hto(X_{n})$, where $X_{n}$
is the composite cobordism,
\[
         X_{n} = W_{n} \cup_{Y_{n+1}} W_{n+1} 
\]
from $Y_{n}$ to $Y_{n+2}$.  
Recall that each cobordism $W_{n}$ arises from the addition of a
single $2$-handle. The core of the $2$-handle in $W_{n+1}$ attaches to
the cocore of the $2$-handle in $W_{n}$ to form an embedded
$2$-sphere,
\[
        E_{n} \subset X_{n}.
\]
This sphere has self-intersection number $-1$, and the boundary
$S_{n}$ of a tubular neighborhood of $E_{n}$ is an
embedded $3$-sphere, giving an alternative decomposition
\[
   X_{n}   = B_n\#_{S_n} Z_n,
\]
where here
$B_n$ is another cobordism from $Y_n$ to $Y_{n+2}$ (obtained by a two-handle
addition), punctured at point, and $Z_n$ is the tubular neighborhood.  (See
Figure~\ref{fig:Square} for a schematic sketch of $X_{1}$.)

There is a diffeomorphism $\tau : X_{n} \to X_{n}$ which is the
identity on $B_{n}$ with
\[
           \tau_{*} [E_{n}] = -[E_{n}] 
\]
in $H_{2}(X_{n};\Z)$. If $\spinc$ is a $\SpinC$ structure on $X_{n}$, then
$\langle c_{1}(\spinc),[E_{n}]\rangle$ is odd; so $\tau$ acts on
$\SpinC(X_{n})$ without fixed points. We can write
\[
              \Hto(X_{n}) = \sum_{\spinc} \Hto(X_{n},\spinc),
\]
and by diffeomorphism-invariance, the contributions from $\spinc$ and
$\tau^{*}(\spinc)$ are equal. Since $\Field$ has characteristic $2$,
the sum vanishes.

The proof that the sequence is exact is considerably harder than the
proof that the composites are zero.  
We use the following straightforward result from homological algebra.

\begin{lemma}
\label{lemma:MCone}
Let $\{C_n\}_{n\in\Zmod{3}}$ be a collection of chain complexes 
over $\Zmod{2}$ and let
$$\{f_n\colon C_n\longrightarrow C_{n+1}\}_{n\in\Zmod{3}}$$
be a collection of chain maps  with the following 
two properties:
\begin{enumerate}
\item
\label{item:NearlyChain}
the composite
$f_{n+1}\circ f_n:C_n\longrightarrow C_{n+2}$ is chain-homotopic to zero, by a 
chain homotopy $H_{n}$:
   \[
          \partial \circ H_{n} + H_{n} \circ \partial = f_{n+1} \circ
          f_{n};
   \]
\item
\label{item:Composite}
the sum 
$$\psi_n=f_{n+2}\circ H_n+ H_{n+1}\circ f_n\colon C_n \longrightarrow
C_{n}$$ 
(which is a chain map) induces isomorphisms on homology,
$(\psi_{n})_{*} : H_{*}(C_{n}) \to
           H_{*}(C_{n})$.
\end{enumerate}
Then the sequence
    \[
             \cdots \longrightarrow
             H_{*}(Y_{n-1})
             \stackrel{(f_{n-1})_{*}}{\longrightarrow}
             H_{*}(Y_{n})
             \stackrel{(f_{n})_{*}}{\longrightarrow}
              H_{*}(Y_{n+1})
             \longrightarrow
             \cdots
    \]
is exact.    
\end{lemma}

We wish to apply the lemma to the chain maps $\mto(W_{n})$; and while
we know that the composites $\mto(W_{n+1})\mto(W_{n})$ induce the zero
map on Floer homology, we need an explicit chain-homotopy in order to
apply the lemma. That is our goal in the next subsection.

\subsection{The first chain homotopy}
\label{subsec:FirstHomotopy}

We now construct the required null-homotopy of
$\mto(W_{n+1})\mto(W_{n})$.  Take $n=1$, and equip $X_{1}$ with a
metric $g$ which is product-like near both the separating
hypersurfaces $Y_{2}$ and $S_{1}$. From this, we construct a family of
metrics $Q(S_{1}, Y_{2})$ parametrized by $T\in\R$ as follows. When
the parameter $T$ for the family is negative, we insert a cylinder
$[T,-T]\times S_1$ normal to $S_1$, and when it is positive, we insert
a cylinder $[-T,T]\times Y_2$ normal to $Y_2$.  We can arrange that
the metric on $S_{1}$ has positive scalar curvature and is close to
the round metric.

There is a corresponding parametrized moduli space
\[
        M_{z}(\fa, X_{1}^{*}, \fb)_{Q}.
\]
As in Section~\ref{subsec:Composite-Cobordisms}, we can complete the
family of Riemannian manifolds: at $T=-\infty$ we obtain the disjoint
union
\[
     X_{1}(-\infty)^{*} = B_{1}^{*} \amalg Z_{1}^{*},     
\]
and at $T=+\infty$ we obtain
\[
     X_{1}(+\infty)^{*} = W_{1}^{*} \amalg W_{2}^{*}.  
\]
The manifold $B_{1}^{*}$ has three cylindrical ends. There is now a
moduli space
\[
           M_{z}(\fa, X_{1}^{*}, \fb)_{\bar Q} \to \bar{Q}(S_{1},
           Y_{2}) 
\]
over $\bar Q(S_{1}, Y_{2}) = [-\infty, \infty]$ and its
compactification $\Mbk_{z}(\fa, X_{1}^{*},\fb)_{\bar Q}$, involving
broken trajectories.

Define quantities $m_{z}(\fa, X_{1},\fb)_{Q} \in \Field$ by counting elements in
zero-dimensional moduli spaces $M_{z}(\fa, X_{1}^{*},\fb)_{Q}$ in the now
familiar way, and define $\bar m_{z}(\fa, X_{1},\fb)_{Q}$ similarly,
using $M^{\red}_{z}(\fa, X_{1}^{*},\fb)_{Q}$.  We use these as the matrix
entries of the linear map
\begin{equation}\label{eq:Hoo-def}
               \Hoo : C^{o}(Y_{1}) \to C^{o}(Y_{3})
\end{equation}
and its seven companions $\Hos$, $\Huo$, $\HIus$, $\Hss$, $\Hsu$,
$\Hus$ and $\Huu$;
and from these we construct a map $\check{H}_{1}$ by the same formula that defined
the chain-homotopy $\Kto$ in
Section~\ref{subsec:Composite-Cobordisms}:
\[
          \check{H}_{1} =  \begin{bmatrix}
                               \Hoo &  \Huo\dsu
                                 +\muo(W_{2})\msu(W_{1})+\duo\Hsu \\
                                 \Hos&  \Hss + \HIus\dsu+
                                 \mIus(W_{2})\msu(W_{1})+
                                 \dIus\Hsu
                       \end{bmatrix}.
\]

\begin{prop}
\label{prop:Homotopies}
    If the chosen perturbation on $S_{1}\cong S^{3}$ is sufficiently
    small, then we have
    \[
                 \dto \circ \check{H}_{1} + \check{H}_{1}\circ \dto
                        = \mto(W_{2})\circ \mto(W_{1})
    \]
    as chain maps from $\Cto_{\bullet}(Y_{1})$ to
    $\Cto_{\bullet}(Y_{3})$.
\end{prop}

\begin{proof}
    The formula closely resembles the formula involving $\Kto$ from
    Proposition~\ref{prop:K-homotopy-Composite}. The chain homotopy $\Kto$ was
    defined using the family of metrics parametrized by the positive
    half, $[0,\infty]$, of the family $\bar Q$. The fiber over $T=0$
    contributed the extra term $\mto(W)$ in the previous formula.

    To prove the present proposition, we proceed as before,
    obtaining identities involving $\Hoo$ and its companions by
    examining $1$-dimensional moduli spaces $\Mbk_{z}(\fa,
    X_{1}^{*},\fb)_{\bar Q}$ and counting their boundary points. The
    new phenomena occur in examining the fiber of $T=-\infty$.

    A typical element of $\Mbk_{z}(\fa,X_{1}^{*},\fb)_{\bar Q}$ in the
    fiber over $T=-\infty$ is a quintuple
    \[
              (\breve{\gamma}_{Y_{1}}, \breve{\gamma}_{S_{1}},
             \breve{\gamma}_{Y_{3}}, \gamma_{B_{1}}, \gamma_{Z_{1}}),
     \]
     where the first three are broken trajectories on the
     corresponding cylinders (each of these may be empty) and
     $\gamma_{B_{1}}$ and $\gamma_{Z_{1}}$ are solutions on the
     corresponding cylindrical-end manifolds.

     To understand which of these decompositions occur, we must
     understand the Floer complex for the three-sphere $S_1$. Since
     the $S_1$ has positive scalar curvature and is simply-connected,
     there is a unique (reducible) critical point in
     $\bonf(S_1)$. After a small perturbation, the critical points in
     $\bonf^\sigma(S_1)$ still lie over a single reducible
     configuration. 
     We label these critical points in
     $\bonf^{\sigma}(S_{1})$ as $\fa_{\lambda_{i}}$, where
     $\lambda_{i}$ are the eigenvalues of a self-adjoint Fredholm
     operator obtained as a small perturbation of the Dirac operator
     on $S^{3}$. (See Example~\ref{ex:ReducibleCritPoints}.)  We
     assume that the $\lambda_{i}$ are strictly increasing as $i$ runs
     through $\Z$ and that $\lambda_{0}$ is the first positive
     eigenvalue:
     \begin{equation}\label{eq:lambda-Order}
                 \cdots \lambda_{-2} < \lambda_{-1} < 0 < \lambda_{0} <
                 \lambda_{1} < \cdots.
     \end{equation}
     It is a consequence of this description that
     $\breve{\gamma}_{S_{1}}$ {\em a priori} 
     live in even-dimensional moduli spaces.
     In fact, by counting dimensions, we see that the trajectories
     $\breve{\gamma}_{S_{1}}$ in this fiber are empty.

     We claim that the possibilities for $\gamma_{Z_{1}}$ come in
     pairs. Specifically, regard $Z_{1}$ as a manifold with boundary
     $S_{1}$ (so that $-S_{1}$ is a boundary component of $B_{1}$).
     For each critical point, we have moduli spaces $M_{z}(Z_{1},
     \fa_{\lambda_{i}})$. The choice of $z$ is equivalent in this
     instance to a choice of $\SpinC$ structure $\spinc$ on $Z_{1}$,
     which in turn is determined by its first Chern class. We write
     $z_{k}$ for the component corresponding to the $\SpinC$ structure
     $\spinc$ with $\langle c_{1}(\spinc),[E_1]\rangle = 2k-1$.

     \begin{lemma}
         \label{lem:CP2-even}
    The following hold for a sufficiently small perturbation on
    $S_{1}$.
     \begin{enumerate}
     \item
        The moduli spaces $M_{z}(Z_{1}, \fa_{\lambda_{i}})$ contain no
        irreducibles. They are empty for $i\ge 0$.
      \item
        For $i < 0$, the
        moduli space $M_{z_{k}}(Z_{1}, \fa_{\lambda_{i}})$ consists of
        a single point when it has formal dimension equal to zero.
      \item
        The formal
        dimensions of $M_{z_{k}}(Z_{1}, \fa_{\lambda_{i}})$ and
        $M_{z_{1-k}}(Z_{1}, \fa_{\lambda_{i}})$ are the same.
        \end{enumerate}
     \end{lemma}

     \begin{proof} The formal dimensions of all the moduli spaces
     $M_{z_k}(Z_{1},\fa_{\lambda_{0}})$ are the same as the moduli
     spaces for the corresponding $\SpinC$ structures over $\mCP$. It
     follows that the formal dimension of
     $M_{z_k}(Z_{1},\fa_{\lambda_{i}})$ is $-k(k-1)-2i-1$ if $i\geq 0$,
     and $-k(k-1)-2i-2$ if $i<0$,. Thus, the
     formal dimensions of all the moduli spaces
     $M_{z}(Z_{1},\fa_{\lambda_{i}})$ with $i\ge 0$ are negative, and
     these moduli spaces are therefore empty. If $i < 0$, then
     $\fa_{\lambda_{i}}$ is boundary-unstable, so the corresponding
     moduli space contains no irreducibles. That the zero-dimensional
     moduli spaces are points is the same phenomenon that underlies
     Proposition~\ref{prop:Reducibles}. Specifically, The last
     statement is a consequence of the diffeomorphism $\tau: Z_{1} \to
     Z_{1}$ which is the identity on the cylindrical end and sends
     $[E_{1}]$ to $-[E_{1}]$.  \end{proof}

     From the lemma, it follows that the number of end-points of a
     $1$-dimensional moduli space $\Mbk_{z}(\fa, X_{1}^{*}, \fb)_{\bar
     Q}$ which lie over $T=-\infty$ is even. The identities which we
     obtain from these moduli spaces are therefore the same as the
     identities for $\Koo$ etc.~in
     Section~\ref{subsec:Composite-Cobordisms}, but without the term
     from $T=0$. For example, we have
     \begin{multline*}
     \moo(W_{1})\moo(W_{0})  + \Hoo\doo  + \doo \Hoo
      +\Huo\dsu\dos  + \duo\Hsu\dos
      +  \duo\dsu \Hos
      \\+ \muo(W_{1})\msu(W_{0})\dos +
      \muo(W_{1})\dsu\mos(W_{0}) +
      \duo\msu(W_{1})\mos(W_{0}) = 0.          
     \end{multline*}
      The chain identity in the proposition follows from this identity
      and its companions.
\end{proof}

\subsection{The second chain homotopy}

Proposition~\ref{prop:Homotopies} gives us the first 
chain homotopy required by Lemma~\ref{lemma:MCone}. Our next goal is
to to construct the second homotopy required by the lemma. This will
be constructed by counting points in moduli spaces
associated to a two-parameter family of metrics on a four-manifold.

%

Specifically, consider the four-manifold $V_1$ obtained as $$V_1=W_1\cup_{Y_2}
W_2 \cup_{Y_3} W_3.$$ This four-manifold contains the $2$-spheres
$E_{1}$ and $E_{2}$, and the $3$-spheres $S_{1}$ and $S_{2}$ which
bound their tubular neighborhoods.  The spheres $E_{1}$ and $E_{2}$
intersect transversely in a single point, with intersection number
$1$. The $3$-spheres $S_{1}$ and $S_{2}$ intersect transversely in a
$2$-torus.
Let $N_{1}$ be a regular neighborhood of $E_{1}\cup E_{2}$ containing
the $3$-spheres $S_{1}$ and $S_{2}$. The boundary of $N_{1}$ is
a separating hypersurface $R_{1}$ in $V_{1}$, diffeomorphic to
$S^{1}\times S^{2}$. The manifold $N_{1}$ is diffeomorphic to the
complement of the neighborhood of a standard circle in $\mCP$, and
gives a decomposition
\[
           V_{1} = U_{1} \cup_{R_{1}} N_{1}.
\]
The manifold $U_{1}$ is obtained topologically by removing a
neighborhood of the curve $\{0\}\times K$ from the cylindrical
cobordism $[-1,1]\times Y_{1}$, where $K$ is the core of the solid
torus $S^{1}\times D^{2}$ that was used in the Dehn filling to create
$Y_{1}$.

\begin{figure}
\mbox{\vbox{\includegraphics{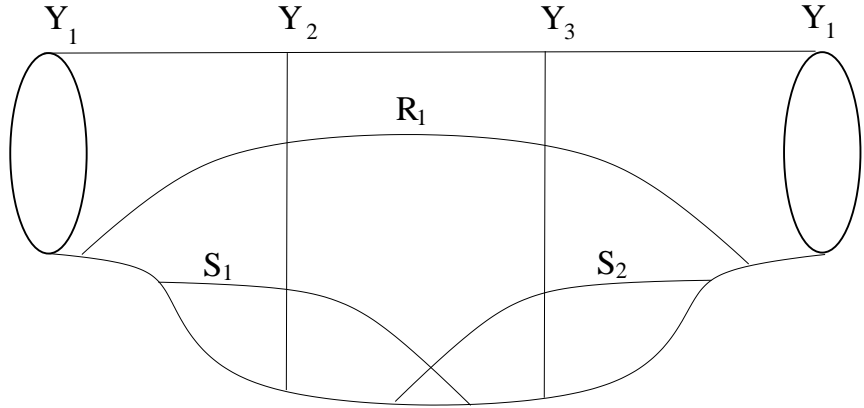}}}
\caption{\label{fig:Pentagon}
{\bf{Breaking up the triple-composite.}}
This indicates the five hypersurfaces which separate $V_1$.}
\end{figure}

\begin{figure}
\mbox{\vbox{\includegraphics{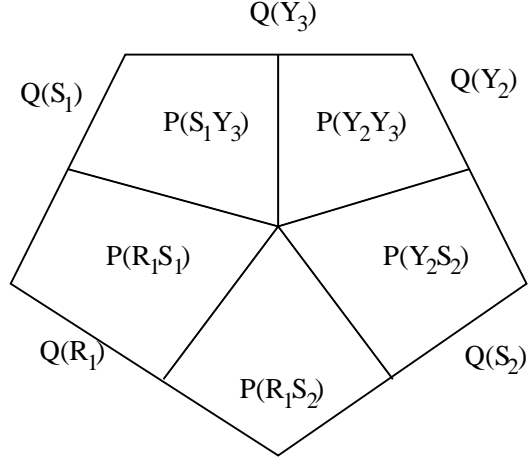}}}
\caption{\label{fig:Metrics}
{\bf{The two-parameter family of metrics.}}  This is a schematic
illustration of the two-parameter family of metrics $\bar P$,
parameterized by a pentagon. The five regions parameterize the five
two-parameter families of metrics where the metrics are varied normal
to two of the five three-manifolds.  Any two two-parameter families
meet along an edge which parameterizes metrics where only one of the
five three-manifolds is pulled out. The five edges on the boundary
parameterize metrics where one of the five three-manifolds is
stretched completely out.}
\end{figure}

In all,
we have five separating hypersurfaces
$Y_2$, $R_1$, $Y_3$, $S_2$, and $S_1$,
as pictured in Figure~\ref{fig:Pentagon}. 
These are arranged cyclically so that 
any one intersects only its two neighbors.  For any two of these
surfaces, say $S$ and $S'$, which do \emph{not} intersect, we can
construct a $2$-parameter family of metrics $P(S,S')$ parametrized by
$\R^{+}\times \R^{+}$, by inserting cylinders $[-T_{S}, T_{S}]\times
S$ and $[-T_{S'},T_{S'}]\times S'$. In the usual way, we complete this
family to obtain a family of Riemannian manifolds over the ``square''
\[
        \bar{P}(S, S') \cong [0,\infty]\times [0,\infty].
\]
There are five such families of metrics, corresponding to the five pairs
of disjoint separating surfaces. The squares fit together along five
common edges, corresponding to families of metrics where just one
of the lengths $T_{S}$ is non-zero.
In this way we set up a  two-parameter family  of metrics $\bar P= \bar
P(R_1,Y_2,Y_3,S_1,S_2)$, as the union of five squares $\bar
P(S, S')$, as shown in Figure~\ref{fig:Metrics}. 
For each of the five hypersurfaces $S$, there are two
edges of $\bar P$ where $T_{S} = \infty$. 
We denote the union of these
two edges by $\bar Q_{S}$. Thus,
\[
       \partial\bar P = \bar{Q}_{S_{2}}\cup \bar{Q}_{Y_{2}} \cup
       \bar{Q}_{Y_{3}} \cup \bar{Q}_{S_{2}} \cup \bar{Q}_{R_{1}}.
\]
By a small adjustment, we can arrange throughout the family that the metrics on $R_1$,
$S_1$, and $S_2$ are standard.

For each pair of critical points $\fa$, $\fb$ in $\Crit(Y_{1})$, we
now have a (parametrized) moduli space $M_{z}(\fa, V_{1}^{*}, \fb)_{P}$ and its
compactification $\Mbk_{z}(\fa, V_{1}^{*},\fb)_{\bar{P}}$. If this moduli
space is zero-dimensional, then it is compact. As usual, we define
quantities
\[
           m_{z}(\fa, V_{1}, \fb)_{\bar{P}},\quad \bar{m}_{z}(\fa, V_{1}, \fb)_{\bar{P}}
\]
by counting points (mod $2$) in zero-dimensional moduli spaces
$M_{z}(\fa, V_{1}^{*}, \fb)_{P}$ and $M^{\red}_{z}(\fa, V_{1}^{*}, \fb)_{P}$
respectively. These are the matrix entries of maps such as
\[
              \Goo : C^{o}_{\bullet}(Y_{1}) \to
              C^{o}_{\bullet}(Y_{1}) 
\]
and its seven companions $\Gos$, $\Guo$, $\GIus$, $\Gss$, $\Gsu$, $\Gus$
and $\Guu$. 

Now suppose that $M_{z}(\fa, V_{1}^{*}, \fb)_{P}$ has dimension $1$. In
the first instance, let us suppose that $\fa$ and $\fb$ are in
$\Crit^{o}(Y_{1})$. As in the earlier settings, we will obtain an
identity
\[
             A^{o}_{o} = 0
\]
for an operator $A^{o}_{o} : C^{o}_{\bullet}(Y_{1}) \to
C^{o}_{\bullet}(Y_{1})$ by enumerating mod $2$ the endpoints of the compactification
$\Mbk_{z}(\fa, V_{1}^{*}, \fb)_{\bar{P}}$, and summing over all $\fa$,
$\fb$ and $z$.  First, there are the endpoints which lie over the
interior of $P \subset \bar{P}$. These arise from strata with either
two factors,
 \begin{equation*}
\begin{gathered}
    \Mu_{z_{1}}(\fa, \fa_{1}) \times M_{z_{2}}(\fa_{1},
    V_{1}^{*},\fb)_{P} \\
    M_{z_{1}}(\fa, V_{1}^{*}, \fb_{1})_{P} \times
    \Mu_{z_{2}}(\fb_{1},\fb),
\end{gathered}
\end{equation*}
or three:
\begin{equation*}
 \begin{gathered}
    \Mu_{z_{1}}(\fa, \fa_{1}) \times \Mu_{z_{2}}(\fa_{1},\fa_{2})
    \times M_{z_{3}}(\fa_{2}, V_{1}^{*},\fb)_{P} \\
        \Mu_{z_{1}}(\fa, \fa_{1}) 
    \times M_{z_{2}}(\fa_{1}, V_{1}^{*},\fb_{1})_{P} \times
     \Mu_{z_{3}}(\fb_{1},\fb) \\
    M_{z_{1}}(\fa, V_{1}^{*}, \fb_{1})_{P} \times \Mu_{z_{2}}(\fb_{1},\fb_{2})
    \times \Mu_{z_{3}}(\fb_{2},\fb),
\end{gathered}
\end{equation*}
just as in Lemma~\ref{lem:BoundaryStrata-Composite}. In the case of
three factors, the middle one is boundary-obstructed.  Together, these
terms contribute
\begin{subequations}
\label{subeqns:Aoo}
\begin{equation}\label{eq:Aoo1}
            \Goo \doo + \doo \Goo +  \duo \dsu \Gos + \duo \Gsu \dos +
            \Guo \dsu \dos
\end{equation}
to the operator $A^{o}_{o}$.  The remaining terms of $A^{o}_{o}$ come
from boundary points in the moduli space that lie over one of the five
parts $\bar{Q}_{S}$ of the boundary $\partial\bar{P}$.  In the case
that $S=S_{1}$ or $S_{2}$, the contribution from $\bar{Q}_{S}$ is
zero. This is because when $T_{S_{1}}$ or $T_{S_{2}}$ becomes
infinite, the manifold $V_{1}$ splits off either $Z_{1}^{*}$ or
$Z_{2}^{*}$, so we can apply Lemma~\ref{lem:CP2-even} to see that the
total number of endpoints over $\bar{Q}_{S_{1}}$ and $\bar{Q}_{S_{2}}$
is even. (This is the same mechanism involved in
the proof of Proposition~\ref{prop:Homotopies}.)  

Next we analyze the endpoints which lie over $\bar{Q}_{Y_{3}}$. When
$T_{Y_{3}}=\infty$, the manifold $V_{1}^{*}$ decomposes as a disjoint
union $X_{1}^{*}
\cup W_{3}^{*}$, where $X_{1}^{*}$ is the composite cobordism
above. In the family parametrized by $\bar{Q}_{Y_{3}}$, the metric on
$W_{3}^{*}$ is constant, while the family of metrics on the component
$X_{1}^{*}$ is the same family $\bar{Q}$ that appeared as
$\bar{Q}(S_{1}, Y_{2})$ in Section~\ref{subsec:FirstHomotopy}.
Endpoints lying over the interior part of the edge, $Q_{Y_{3}}\subset
\bar{Q}_{Y_{3}}$ may belong to strata with two
factors, which have the form
\[
           M_{z_{1}} (\fa, X_{1}^{*},\fa_{1})_{\bar{Q}} \times
               M_{z_{2}}(\fa_{1}, W_{3}^{*}, \fb);
\]
or they may belong to strata with three factors, one of which is
boundary obstructed, as in case~\ref{item:ThreeFactorComposite} of
Lemma~\ref{lem:BoundaryStrata-Composite}.  Altogether, these terms
contribute four terms to $A^{o}_{o}$,
\begin{multline}
        \label{eq:Aoo2}
       \moo(W_{3})\Hoo(X_{1}) + 
       \muo(W_{3})\Hsu(X_{1})\dos \\+
      \muo(W_{3})\dsu\Hos(X_{1}) +
      \duo\msu(W_{3})\Hos(X_{1})
\end{multline}
where the operators $H^{*}_{*} = H^{*}_{*}(X_{1})$ are those defined
at \eqref{eq:Hoo-def}.  The contributions from endpoints lying over
${Q}_{Y_{2}}$ are similar: we obtain
\begin{multline}
        \label{eq:Aoo3}
       \Hoo(X_{2})\moo(W_{1}) +
       \Huo(X_{2})\msu(W_{1})\dos \\+
      \Huo(X_{2})\dsu\mos(W_{1}) +
      \duo\Hsu(X_{2})\mos(W_{1}).
\end{multline}
There is also one possible type of endpoint that occurs at the vertex
of $\bar{P}$ where $\bar{Q}_{Y_{3}}$ and $\bar{Q}_{Y_{2}}$ meet: these
lie in a moduli space
\[
        M_{z_{1}}(\fa, W_{1}^{*}, \fa_{1}) \times
         M_{z_{2}}(\fa_{1}, W_{2}^{*}, \fa_{2})
         \times
          M_{z_{3}}(\fa_{2}, W_{3}^{*}, \fb),
\]
where the middle factor is boundary-obstructed. These contribute a
term
\begin{equation}\label{eq:Aoo3a}
          \muo(W_{3}) \msu(W_{2}) \mos(W_{1})
\end{equation}
to $A^{o}_{o}$.

When $T_{R_{1}}=\infty$, we have a decomposition of $V_{1}^{*}$ into
two pieces
\[
          N_{1}^{*}\cup   U_{1}^{*},
\]
where $U_{1}$ and $N_{1}$ are as above. The manifold $U_{1}^{*}$ has
three ends. We regard $U_{1}$ as a cobordism from $R_{1}\amalg Y_{1}$
to $Y_{1}$, and $N_{1}$ as a manifold with oriented boundary $R_{1}$.
The $1$-parameter family of metrics $\bar{Q}_{R_{1}}$ is constant on $U_{1}$, and we
have moduli spaces
\[
         M_{z}(N_{1}, \fa')_{\bar Q}    \text{\quad
             and\quad}  M_{z}(\fa', \fa, U_{1}^{*}, \fb)
\]
or $\fa'\in \Crit(R_{1})$ and $\fa,\fb\in\Crit(Y_{1})$. Here
$\bar{Q}\cong[-\infty, \infty]$
is the family of metrics $\bar{Q}(S_{1}, S_{2})$ on $U_{1}$ which
stretches along $S_{1}$ when $T$ is negative and $S_{2}$ when $T$ is
positive.  On $N_{1}$ we can count points in zero-dimensional moduli
spaces $M_{z}(N_{1}, \fa')_{\bar Q}$, and so define elements
\[
\begin{aligned}
n_{s} &\in C^{s}_{\bullet}(R_{1}), \\
          n_{o} &\in C^{o}_{\bullet}(R_{1})\\
        \bar{n}_{s} &\in C^{s}_{\bullet}(R_{1}) \\
        \bar{n}_{u} &\in C^{u}_{\bullet}(R_{1}) .
        \end{aligned}
\]
(The last two count points in zero-dimensional moduli spaces
$M^{\red}_{z}(N_{1}, \fa')$.) The situation simplifies slightly, on
account of the following lemma.

\begin{lemma}
        \label{lem:R1-lemma1}
  If the perturbation on $R_{1}$ is sufficiently small, then there are
  no irreducible critical points (so $n_{o}$ is zero), and
  no irreducible trajectories on $\R\times R_{1}$. The perturbation
  can be chosen so that the one-dimensional reducible
  trajectories come in pairs, so $\dss$, $\dsu$, $\dus$ and $\duu$ are
  all zero. The invariant $\bar{n}_{s}$ is zero also.
\end{lemma}

\begin{proof}
        We postpone the proof to Section~\ref{subsec:Calculation}
        below, where we also calculate $\bar{n}_{u}$. 
\end{proof}

The zero-dimensional moduli spaces  $M_{z}(\fa',\fa,U_{1}^{*},\fb)$
provide the matrix entries of maps
\[
            \muoo : C^{u}_{\bullet}(R_{1}) \otimes
            C^{o}_{\bullet}(Y_{1}) \to C^{o}_{\bullet}(Y_{1})
\]
as well as $\muuo$, $\muos$ and $\mIuus$, while the reducible parts of
these moduli spaces define $\msss$, $\mssu$, $\msus$, $\msuu$,
$\muss$, $\musu$, $\muus$ and $\muuu$. Of these eight maps defined by
zero-dimensional moduli spaces of reducible solutions, the maps
$\msss$, $\msus$ and $\musu$ arise from boundary-obstructed moduli
spaces. The moduli spaces $M_{z}(\fa', \fa,U_{1}^{*},\fb)$
contributing to $\mssu$ are \emph{doubly} boundary obstructed (or
boundary obstructed with corank 2, in the notation of \cite{KMBook}):
these zero-dimensional moduli spaces have formal dimension
$\gr_{z}(\fa',\fa, U_{1}, \fb) = -2$.

We can now enumerate the end-points belonging to $\bar{Q}_{R_{1}}$ in
the $1$-dimensional moduli spaces $\Mbk_{z}(\fa, V_{1}^{*},
\fb)_{\bar{P}}$ that contribute to $A^{o}_{o}$.
First there are points belonging to strata with two
factors, of the form
\[
           M_{z_{1}}(N_{1}^{*}, \fa')_{\bar{Q}} \times M_{z_{2}}(\fa',
           \fa, U_{1}^{*}, \fb),
\]
where $\fa'$ is necessarily boundary-unstable (so the solution on
$N_{1}^{*}$ is reducible). Next we should look for points belonging to strata
with three factors, one of which is boundary-obstructed; but when
$\fa$ and $\fb$ are irreducible, there are no such contributions.
Finally, there are points belonging to strata with four factors, one of
which is doubly boundary-obstructed. These have the form
\[
          M_{z_{1}}(N_{1}^{*},\fa')_{\bar Q} \times
          M_{z_{2}}(\fa,\fa_{1}) \times
          M_{z_{3}}(\fa',\fa_{1},U_{1}^{*},\fb_{1})\times
          M_{z_{4}}(\fb_{1}, \fb), 
\]
where $\fa'$ is boundary-stable, $\fa_{1}$ is boundary-stable, and
$\fb_{1}$ is boundary-unstable. From these we obtain the final two
terms in $A^{o}_{o}$:
\begin{equation}\label{eq:Aoo4}
          \muoo(\bar{n}_{u}\otimes\mathord{\cdot}) +
         \duo \mssu(n_{s}\otimes \dos(\mathord{\cdot})).
\end{equation}
The identity $A^{o}_{o}=0$ has sixteen terms, from
\eqref{eq:Aoo1}--\eqref{eq:Aoo4}:
\end{subequations}
\begin{multline*}
                \Goo \doo + \doo \Goo +  \duo \dsu \Gos + \duo \Gsu \dos +
            \Guo \dsu \dos \\
            +    
       \moo(W_{3})\Hoo(X_{1}) + 
       \muo(W_{3})\Hsu(X_{1})\dos \\+
      \muo(W_{3})\dsu\Hos(X_{1}) +
      \duo\msu(W_{3})\Hos(X_{1}) \\
      +  \Hoo(X_{2})\moo(W_{1}) +
       \Huo(X_{2})\msu(W_{1})\dos \\+
      \Huo(X_{2})\dsu\mos(W_{1}) +
      \duo\Hsu(X_{2})\mos(W_{1})\\ +  \muo(W_{3}) \msu(W_{2}) \mos(W_{1})
      +
          \muoo(\bar{n}_{u}\otimes\mathord{\cdot}) +
         \duo \mssu(n_{s}\otimes \dos(\mathord{\cdot})) = 0.
\end{multline*}
There are three similar identities, $A^{o}_{s}=0$, $A^{u}_{o}=0$ and
$A^{u}_{s}=0$, coming from the three other types of $1$-dimensional
moduli spaces that contain irreducibles. 
In full, these are
\begin{equation}
	\begin{aligned}{}
	A^{o}_{s} &=
		\dss \Gos + \Gss \dos 
		+ \dos \Goo + \Gos \doo 
		+ \dIus \dsu\Gos + \dIus \Gsu \dos 
		+ \GIus \dsu \dos \\ &
		+ \msss(n_{s}\otimes \dos\mathord{\cdot})
		+ \dIus \mssu(n_{s}\otimes\dos \mathord{\cdot}) \\&
		+ \muos({\overline n}_{u}\otimes \mathord{\cdot}) \\&
		+ \Hss(X_{2}) \mos(W_{1})
		+ \Hos(X_{2}) \moo(W_{1})\\&
		+ \dIus \Hsu(X_{2}) \mos(W_{1})
		+ \HIus(X_{2})\dsu\mos(W_{1})
		+ \HIus(X_{2})\msu(W_{1})\dos \\&
		+ \mss(W_{3})\Hos(X_{1})
		+ \mos(W_{3})\Hoo(X_{1}) \\&
		+ \dIus \msu(W_{3}) \Hos(X_1)
		+ \mIus(W_{3})\dsu \Hos(X_{1})
		+ \mIus(W_{3})\Hsu(X_{1})\dos \\ &
		+\mIus(W_{3})\msu(W_{2})\mos(W_{1}) \\
	\end{aligned}
	\nonumber
\end{equation}
\begin{equation}
     \begin{aligned}[]
	A^{u}_{o} &=
		\doo\Guo 
		+ \Goo\duo 
		+ \duo\Guu 
		+ \Guo\duu
		+ \duo\dsu\GIus 
		+ \duo\Gsu\dIus 
		+ \Guo\dsu\dIus \\ &
		+ \duo\mssu(n_{s}\otimes \dIus\mathord{\cdot})
		+ \duo\msuu(n_{s}\otimes \mathord{\cdot}) \\ &
		+ \muuo({\overline n}_{u}\otimes \mathord{\cdot}) \\ &
		+ \Hoo(X_{2})\muo(W_{1})
		+ \Huo(X_{2})\muu(W_{1}) \\ &
		+ \duo\Hsu(X_{2})\mIus(W_{1}) 
		+ \Huo(X_{2})\dsu\mIus(W_{1})
		+ \Huo(X_{2})\msu(W_{1})\dIus \\ &
		+ \moo(W_{3})\Huo(X_{1})
		+ \muo(W_{3})\Huu(X_{1}) \\ &
		+ \duo\msu(W_{3})\HIus(X_{1})
		+ \muo(W_{3})\dsu\HIus(X_{1})
		+ \muo(W_{3})\msu(X_{1})\dIus \\ &
		+ \mIus(W_{3})\msu(W_{2})\mos(W_{1}) \\
	\end{aligned}
	\nonumber
\end{equation}
\begin{equation}
	\begin{aligned}[]
	A^{u}_{s} &=
		\Gus
		+ \dss \GIus
		+ \Gss\dIus
		+ \dIus \Guu
		+ \GIus \duu 
		+ \dIus\dsu\GIus
		+ \dIus\Gsu\dIus
		+ \GIus\dsu\dIus  \\&
		+ \dos \Guo + \Gos \duo\\&
		+ \dIus\mssu(n_{s}\otimes \dIus\mathord{\cdot})
		+ \dIus\msuu(n_{s}\otimes\mathord{\cdot})
		+ \msus(n_{s}\otimes \cdot) +
			\msss(n_{s}\otimes \dIus\mathord{\cdot}) \\ &
		+ \mIuus({\overline n}_{u}\otimes \mathord{\cdot}) \\ &
		+ \Hss(X_{2})\mIus(W_{1})
		+ \HIus(X_{2})\muu(W_{1}) \\ &
		+ \dIus \Hsu(X_{2})\mIus(W_{1})
		+ \HIus(X_{2})\dsu\mIus(W_{1}) \\&
		+ \HIus(X_{2})\msu(W_{1})\dIus
 		+ \Hos(X_{2})\muo(W_1) \\ &
		+ \mss(W_{3})\HIus(X_{1})
		+ \mIus(W_{3})\Huu(X_{1})\\ &
		+ \dIus\msu(W_{3})\HIus(X_{1})
		+ \mIus(W_{3})\dsu\HIus(X_{1}) \\&
		+ \mIus(W_{3})\Hsu(X_{1})\dIus 
		+ \mos(W_{3})\Huo(X_{1})\\ &
		+ \mIus(W_{3})\msu(W_{2})\mIus(W_{1})
     \end{aligned}
	\nonumber
\end{equation}
 
There are four simpler identities involving only the reducible moduli
spaces: these are the vanishing of expressions
$\bar{A}^{s}_{s}$, $\bar{A}^{s}_{u}$,
$\bar{A}^{u}_{s}$ and $\bar{A}^{u}_{u}$, where for example
\begin{multline*}
     \bar{A}^{s}_{s} =  \Gss \dss + \Gus\dsu + \dss \Gss + \dus\Gsu   \\
            +    
       \mss(W_{3})\Hss(X_{1}) +  \mus(W_{3})\Hsu(X_{1}) +
        \Hss(X_{2})\mss(W_{1}) + \Hus(X_{2})\msu(W_{1}) 
       \\
      +
          \muss(\bar{n}_{u}\otimes\mathord{\cdot}) .
\end{multline*}

We define an operator
\[
             \check{L} : \Cto_{\bullet}(Y_{1}) \to
             \Cto_{\bullet}(Y_{1})
\]
by combining some of the contributions from $\bar{Q}_{R_{1}}$:
we write
\begin{equation}
    \check{L} = 
    \begin{bmatrix}
          \loo &  \luo\dsu+ \duo\lsu \\
                                 \los&  \lss + \lIus\dsu +
                                 \dIus\lsu
    \end{bmatrix},
\end{equation}
where
\begin{equation*}
    \loo = \muoo(\bar{n}_{u} \otimes \mathord{\cdot} ),
\end{equation*}
and so on. (The term $\bar{n}_{u}$ appears in the definition of all of
these.)
In words, when $\gr_{z}(\fa,V_1,\fb)=-1$ 
the $e_{\fb}$ component of $\check{L}(e_{\fa})$ counts points in the
zero-dimensional strata of $\Mubk_{z}(\fa,\fb)$ which are broken along
a critical point in $\Crit^{u}(R_1)$.

We define $\check{G} : \Cto_{\bullet}(Y_{1}) \to
\Cto_{\bullet}(Y_{1})$ by the formula
\[
      \check{G} =     \begin{bmatrix}
          a&  b\\
         c&  d
    \end{bmatrix},
\]
where 
\begin{equation*}
    \begin{aligned}
        a &= \Goo \\
        b &= \duo\Gsu + \Guo\dsu + \muo\Hsu + \Huo\msu +
        \duo(\mssu(n_{s}\otimes \mathord{\cdot})) \\
        c &= \Gos \\
        d &= \Gss + \dIus\Gsu + \GIus\dsu + \mIus\Hsu + \HIus\msu
            + \dIus\mssu(n_{s}\otimes\mathord{\cdot}).
    \end{aligned}
\end{equation*}
Here, we have written $\muo\Hsu$ for example as an abbreviation for
$\muo(W_{3})\Hsu(X_{1})$, because no ambiguities arise in the formulae.
In words, if $\gr(\fa,V_1^*,\fb)=-2$, the $e_{\fb}$ component
of $\check{G}(e_{\fa})$ counts points in the zero-dimensional
strata of
$\Mubk_{z}(\fa,\fb)$.

\begin{prop}
\label{prop:Homotopy2}
    We have the identity
    \[
               \dto\circ\check{G} + \check{G}\circ\dto = \mto_{3}\circ\check{H}_{1} +
               \check{H}_{2}\circ\mto_{1} + \check{L},
    \]
    where $\mto_{3} = \mto(W_{3})$ and $\check{H}_{1}$,
    $\check{H}_{2}$ are the operators from
    Proposition~\ref{prop:Homotopies}, using
    $X_{1}$ and $X_{2}$ respectively.
\end{prop}

\begin{proof}
	In addition to the identities $A^*_*=0$ and ${\overline
	A}^*_*=0$, there are the identities arising from pieces of the
	cobordism. For example, we can consider the three-ended
	manifold $U_1^*$. On this, once again, we can enumerate ends
	of the one-dimensional moduli spaces
	$M_{z}(\fa',\fa,U_1^*,\fb)$.  These give relations which are
	formally similar to the relations coming from a two-ended
	cobordism from $Y_1$ to $Y_2$ (since differentials for the
	Floer homology of $R_1$ are trivial,
	c.f. Lemma~\ref{lem:R1-lemma1}). Looking at ends of
	irreducible moduli spaces, we get four relations of the type
	$B^{uo}_{o}$, $B^{uu}_{s}$, $B^{uu}_{o}$ and $B^{uu}_{s}$ For
	example, the relation of the form $B^{uu}_{s}=0$, can be
	written out as:
\begin{equation}
	\begin{aligned} B^{uu}_{s} &=
	\muus(\mathord{\cdot}\otimes\mathord{\cdot}) +
	\muss(\mathord{\cdot}\otimes \dIus(\mathord\cdot)) +
	\mIuus(\mathord{\cdot}\otimes \duu(\mathord\cdot)) \\ & +
	\dss\mIuus({\mathord\cdot}\otimes{\mathord\cdot})
	+\dIus\muuu({\mathord\cdot}\otimes{\mathord\cdot}) 
	\end{aligned}
\end{equation}
	In addition to these, we have eight relations
	coming from looking at ends of reducible moduli spaces
	${\overline B}^{**}_{*}=0$ (where here each $*$ can be either
	symbol $u$ or $s$).
	There are the relations coming from ends of moduli
	spaces for the $X_1$ and $X_2$ 
	(with its two-parameter families of metrics), the ends of the
	moduli spaces for $W_i$, and finally, the ends of moduli spaces
	for $\R\times Y_i$.
	Putting all these together, we get the proposition.
\end{proof}

We can put the above outlined proof of
Proposition~\ref{prop:Homotopy2} into a more conceptual framework.
Throughout the following discussion, we fix two critical points $\fa$
and $\fb$ with $\gr(\fa,V_1^*,\fb)=-1$ (where both $\fa$ and $\fb$ are
in $\Crit^{s}\cup \Crit^{o}$, if we are considering the case
of $\Hto$, for example). 
We count the ends of those one-dimensional strata in
$\Mbk_{z}(\fa,V_1^*,\fb)_{\overline P}$.	
Clearly, these ends count points in the zero-dimensional strata in 
$\Mbk_{z}(\fa,V_1^*,\fb)_{\overline P}$.
We claim that the total sum of these zero-dimensional strata
counts the $e_{\fb}$ component
of the image of $e_{\fa}$ under the map
\begin{equation}
\label{eq:DesiredEquation}
\dto\circ\check{G} + \check{G}\circ\dto + \mto_{3}\circ \check{H}_1 + 
               \check{H}_{2}\circ\mto_{1} + \check{L}, 
\end{equation}
which must therefore be zero. 

The verification can be broken into the following three steps.  Recall
that the strata of $\Mbk_{z}(\fa,V_1^*,\fb)_{\overline P}$ consist of
fibered products over various critical points of moduli spaces. We
call these critical points with multiplicity (if the same critical point
appears more than once) {\em break points} for the stratum.  We say
that a stratum in $\Mbk_{z}(\fa,V_1^*,\fb)_{\overline P}$ has a {\em
good break} if at least one
of its break points  lies in
$(\Crit^{s}\cup
\Crit^{o})(Y_i)$ (with $i\in\{1,2,3\}$) or in $\Crit^{u}(R_1)$. 
One must verify first that the non-empty, zero-dimensional strata in
$\Mbk_{z}(\fa,V_1^*,\fb)_{\overline P}$ which have a good break have,
in fact, a unique (i.e. with multiplicity one) good break.  This
follows from a straightforward dimension count, after listing all 
possible good breaks. Indeed, in
view of the definitions of the maps $\Gto$, $\check{H}$, $\mto$, and
$\Lto$, the above dimension counts show that  
the strata for which the good break occurs along $Y_1$ are
counted in $\dto\circ\Gto+\Gto\circ
\dto$, those where it occurs along $Y_3$ are counted in $\mto_{3}\circ
\check{H}_{1}$, those where it occurs along $Y_2$ are counted in
$\check{H}_{2}\circ\mto_{1}$, and those where it occurs along $R_1$ are
counted in $\Lto$. Second, one verifies that any of the
zero-dimensional strata with one good break appear uniquely as
boundaries of one-dimensional moduli spaces in
$\Mbk_{z}(\fa,V_1^*,\fb)_{\overline P}$.  Finally, one verifies that
any of the zero-dimensional strata in
$\Mbk_{z}(\fa,V_1^*,\fb)_{\overline P}$ which have no good break, and
hence are not accounted for in Equation~\eqref{eq:DesiredEquation},
are counted exactly twice: they appear in the boundaries of two
distinct one-dimensional strata in $\Mbk_{z}(\fa,V_1^*,\fb)_{\overline
P}$.

\subsection{Calculation}
\label{subsec:Calculation}

The plan of the proof is to deduce
Theorem~\ref{thm:TwoHandleSurgeries-A} from Lemma~\ref{lemma:MCone}.
We have already constructed the chain homotopies referred to in the
first part of the lemma: these are the chain homotopies $\check{H}_{n}$.
To verify the hypothesis in the second part of the lemma, it is enough
to verify that $\check{L}$ induces isomorphisms in homology, because
of Proposition~\ref{prop:Homotopy2}.  That is the content of the next
proposition.

\begin{prop}
 \label{prop:L-Calc}
    The map $\check{L} : \Cto_{\bullet}(Y_{1}) \to
    \Cto_{\bullet}(Y_{1})$ induces isomorphisms in homology. Indeed,
    the resulting map on $\Hto_{\bullet}(Y_{1})$ is multiplication by
    the power series
    \[
             \sum_{k\ge 0} U^{k(k+1)/2},
    \]
    which has leading coefficient $1$.
\end{prop}

We begin by examining the Floer complexes of the manifold $R_{1} =
S^{1}\times S^{2}$, equipped with a standard metric and small regular
perturbation $\pertY$ from the class $\Pert(R_{1})$. With no
perturbation, the critical points $(A,\Phi)$ of the Chern-Simons-Dirac
functional all belong to the $\SpinC$ structure $\spinc_{0}$ with
$c_{1}(\spinc_{0}) = 0$. They
are simply the reducible solutions $(A,0)$ with $A^{\ttr}$ flat. In
$\bonf(R_{1})$, they form a circle. For any flat connection $A$, the
corresponding Dirac operator $D_{A}$ on $R_{1}$ has no kernel, so
there is no spectral flow between any two points in the circle.  We can choose
our small perturbation to restrict to a standard Morse function on this
circle, with one maximum $\alpha^{1}$ and one minimum $\alpha^{0}$. We
also need to arrange that the corresponding perturbed Dirac operators
a these two points have simple eigenvalues, and we choose the
perturbation small enough so as not to introduce any spectral flow on
the paths joining $\alpha^{1}$ to $\alpha^{0}$.  In the blow-up
$\bonf^{\sigma}(R_{1})$, each of these critical points gives rise to a
collection of critical points $\fa^{0}_{i}$ and $\fa^{1}_{i}$,
corresponding to the eigenvalues $\lambda_{i}$ of the perturbed Dirac
operator. We again assume these eigenvalues are labeled in increasing
order, with $\lambda_{0}$ the first positive eigenvalue, as in
\eqref{eq:lambda-Order}.  The points $\fa^{0}_{i}$ and $\fa^{1}_{i}$
are boundary-stable when $i\ge 0$ and boundary-unstable when $i<0$.

The trajectories on $\R\times R_{1}$ are all reducible, because
$R_{1}$ has positive scalar curvature. Their images in $\bonf(R_{1})$
are therefore either constant paths at $\alpha^{0}$ or $\alpha^{1}$,
or one of the two trajectories $\gamma$, $\gamma'$ from $\alpha^{1}$
to $\alpha^{0}$ on the circle. For each $i$, there are two
trajectories $\gamma_{i}$ and $\gamma'_{i}$ lying over $\gamma$ and
$\gamma'$ respectively in a moduli space $M_{z}(\fa^{1}_{i},
\fa^{0}_{i})$. These are the only trajectories belonging to
$1$-dimensional moduli spaces, so all boundary maps are zero in the
Floer complexes.   Thus $\Hto_{\bullet}(R_{1}) \cong
\Cto_{\bullet}(R_{1})$, which has generators
\[
            e^{0}_{i}, e^{1}_{i}, \qquad (i\ge 0)
\]
corresponding to the critical points $\fa^{0}_{i}$ and $\fa^{1}_{i}$,
while $\Cfrom_{\bullet}(R_{1})$ has generators
\[
            e^{0}_{i}, e^{1}_{i}, \qquad (i< 0).
\]
We identify $\Gr(Y,\spinc_{0})$ with $\Z$ in such a way that
$e^{0}_{0}$ belongs to $\Cto_{0}(R_{1})$. Then
$e^{\mu}_{i}$ is in grading $\mu + 2i$ for $i$ positive and in grading
$\mu + 2i + 1$ for $i$ negative.

The manifold $N_{1}$ has boundary $R_{1}$, and its homology is
generated by the classes $[E_{1}]$ and $[E_{2}]$ of the two spheres. A
$\SpinC$ structure $\tspin$ on $N_{1}$ whose restriction to $R_{1}$ is
$\spinc_{0}$ is uniquely determined by the evaluation of
$c_{1}(\tspin)$ on $[E_{1}]$. For each $k\in \Z$, we write
$\tspin_{k}$ for the $\SpinC$ structure whose first Chern class
evaluates to $2k+1$ on $[E_{1}]$.   The $\SpinC$ structures
$\tspin_{k}$ and $\tspin_{-1-k}$ are complex conjugates. We write
\[
    M_{k}(N_{1}^{*}, \fa')_{\bar{Q}}
\]
for the union of the moduli spaces belonging to components $z$ which
give rise the $\SpinC$ structure $\tspin_{k}$.  The family $\bar{Q}$
is the same $1$-parameter family $\bar{Q}(S_{1}, S_{2})$ that appeared above.
The following lemma and its two corollaries are straightforward.

\begin{lemma}
    The dimension of the moduli space $M_{k}(N_{1}^{*},
    \fa^{\mu}_{i})_{\bar{Q}}$ is given by
    \[
                \dim M_{k}(N_{1}^{*},
    \fa^{\mu}_{i})_{\bar{Q}} =
    \begin{cases}
         -\mu - k(k+1) - 2i, &i\ge0,\\
         -\mu - k(k+1) - 2i - 1, &i< 0.\\
    \end{cases}
    \]
    \qed
\end{lemma}

\begin{cor}
    The only non-empty moduli spaces $M_{k}(N_{1}^{*},\fa')_{\overline
    Q}$ with $\fa'$ boundary-stable occur when $k=0$ or $-1$ and $\fa'
    = \fa^{0}_{0}$, in which case the moduli space is
    zero-dimensional.  The moduli spaces $M^{\red}_{k}(N_{1}^{*},
    \fa')_{\bar{Q}}$ are empty for all boundary-stable $\fa'$.
    \qed
\end{cor}

\begin{cor}
    The zero-dimensional moduli spaces $M_{k}(N_{1}^{*}, \fa')_{\bar{Q}}$, with
    $\fa'$ boundary-unstable, are the moduli spaces
    \begin{equation}\label{eq:ContributorsToNu}
       M_{k}(N_{1}^{*}, \fa^{1}_{i_{k}})_{\bar{Q}}, \qquad i_{k} = -1 -
       k(k+1)/2.
    \end{equation}  \qed
\end{cor}

At this point, we have verified all parts of
Lemma~\ref{lem:R1-lemma1}, and in addition we can now express
$\bar{n}_{u} \in C^{u}_{\bullet}(R_{1})$ as
\[
     \bar{n}_{u} = \sum_{k\in \Z} a_{k} \,e^{1}_{-1-k(k+1)/2} 
\]
where $a_{k}$ counts points in the moduli space
\eqref{eq:ContributorsToNu} (which consists entirely or reducibles,
because $\fa^{1}_{i_{k}}$ is boundary-unstable).

\begin{lemma}
  \label{lem:KeyOdd}
    For all $k\in \Z$, the sum $a_{k} + a_{-1-k}$ is $1$ mod $2$.
\end{lemma}

\begin{proof}
  Let $g$ be any Riemannian metric on $N_{1}^{*}$ that is standard on
  the end. Fix a $\SpinC$ structure $\tspin_{k}$ on $N_{1}$. Because
  $N_{1}^{*}$ has no first homology and no self-dual, square-integrable harmonic $2$-forms,
  there is a unique $\SpinC$ connection
    \[
             A = A(k,g)
    \]
  in the associated spin bundle
  $S^{+} \to N_{1}$ with $L^{2}$ curvature satisfying the abelian
  anti-self-duality equation $F^{+}_{A^{\ttr}}=0$.  On the cylindrical end,
  $A^{\ttr}$ is asymptotically flat, so $A$ defines a point
  \[
                \theta_{k}(g) \in \mathcal{S},
  \]
  where $\mathcal{S}\subset \bonf(R_{1})$ is the circle of flat $\SpinC$
  connections.  This depends only on $k$ and $g$.

  Fix a Spin structure on $R_{1}$ whose associated $\SpinC$ structure
  is $\spinc_{0}$. This fixes an isomorphism between $\spinc_{0}$ and
  its complex conjugate. Complex conjugation now gives an involution
  on the circle, $\sigma : \mathcal{S} \to \mathcal{S}$, with two
  fixed points $s_{+}$ and $s_{-}$.  The isomorphism between
  $\spinc_{0}$ and its conjugate extends to an isomorphism
  $\bar{\tspin}_{k}\to \tspin_{-1-k}$, and we therefore have
  \begin{equation}\label{eq:Conjugation-Circle} \theta_{-1-k}(g) =
  \sigma \theta_{k}(g).  \end{equation}

  Consider now the family of metrics $\bar{Q} = \bar{Q}(S_{1}, S_{2})$
  on $N_{1}$. As $T$ goes to $-\infty$, the manifold $N_{1}^{*}$
  decomposes into two pieces, one of which has cylindrical ends
  $\R^{-}\times S_{1}$ and $\R^{+}\times R_{1}$. This piece, call it
  $T_{1}^{*}$, carries no $L^{2}$ harmonic $2$-forms (it is a
  punctured $S^{2}\times D^{2}$ with cylindrical ends), so the map
  $\theta_{k}$ extends to continuously to $T=-\infty$ and
  $\theta_{k}(-\infty)$ is one of $s_{+}$ or $s_{-}$. The same applies
  to $T=+\infty$, so we obtain a map \[ \theta_{k} :\bar{Q} \to
  \mathcal{S} \] with \[ \theta_{k}(\pm\infty) \in \{ s_{+}, s_{-} \}.
  \] We also have $\theta_{k}(+\infty) = \theta_{-1-k}(+\infty)$, and
  the same with $T=-\infty$. Thus the two maps \[ \theta_{k},
  \theta_{-1-k} : [-\infty, \infty] \to \mathcal{S} \] together define
  a mod 2 $1$-cycle in $\mathcal{S}$.  The statement of the lemma
  follows from the assertion that this $1$-cycle has non-zero degree
  mod $2$.

  Specifically, let $\Theta^{k}\colon S^1\longrightarrow \mathcal{S}$
  be the cycle obtained by joining $\theta_k$ and $\theta_{-1-k}$. We
  show that for generic $x\in\mathcal{S}$, the space
  $\Theta_k^{-1}(x)$ is cobordant to \begin{equation} \label{eq:Ends}
  \left(M_k^{\mathrm{ab}}(N_1^*,\fa^1)_{\bar{Q}}\bigcup
  M_{-1-k}^{\mathrm{ab}}(N_1^*,\fa^1)_{\bar{Q}}\right)\times
  M^{\ab}(\fa^1,([0,1]\times R_1)^*,{\mathcal{S}})
  \times_{\mathcal{S}}{\{x\}\times{S}} \end{equation} Here,
  $M_k^{\mathrm{ab}}(N_1^*,\fa^1)_{\bar{Q}}$ denotes the moduli space
  of solutions to the perturbed abelian anti-self-duality equiations.
  Similarly,
  $M^{\mathrm{ab}}(\fa^1,([0,1]\times R_1)^*,{\mathcal{S}})$ denotes
  the moduli space of solutions to the perturbed abelian
  anti-self-duality equations, where we use the perturbation at the
  $t=\infty$ and no perturbation and the $t=-\infty$ end. In
  particular, this moduli space admits a map by taking the limit as
  $t\goesto \infty$ to ${\mathcal{S}}$.  
  The claimed compact cobordism is
  induced by taking a one-parameter family of perturbations indexed by
  $T\in \R^+$ on $N_1^*$ which are supported on ever-longer pieces of
  the attached cylinder. The fiber over $T=0$ of this cobordism is
  $\{x\}\times_{\mathcal S} \left(M^{\mathrm{ab}}_{k}(N_1^*,\mathcal{S})
  \cup M^{\mathrm{ab}}_{-1-k}(N_1^*,\mathcal{S})\right)$, whose
  number of points coincides with the stated degree, while the fiber
  over $T=\infty$ is the space described in~\eqref{eq:Ends}.  In
  particular, if the stated degree is odd, then so is the number of
  points in $M_k^{\mathrm{ab}}(N_1^*,\fa^1) \cup
  M_{-1-k}(N_1^*,\fa^1)$. But $M_k^{\mathrm{ab}}(N_1^*,\fa^1)$ is
  identified with $M_k(N_1^*,\fa^1_{i_k})$.

  It remains now to show that the degree of $\Theta^k$ is non-zero modulo two.
  This in turn is equivalent to saying that
  \[
              \theta_{k}(\infty) \ne \theta_{k}(-\infty),
  \]
  because of the relationship \eqref{eq:Conjugation-Circle}.  So we
  must prove that $\theta_{k}: \bar{Q} \to \mathcal{S}$ is a path
  joining $s_{-}$ to $s_{+}$.

  To get a concrete model for $\theta_{k}$, choose a standard closed
  curve $\delta$ representing the generator of $H_{1}(R_{1})$, and let
  $\Sigma\subset N_{1}^{*}$ be a topological open disk with a
  cylindrical end $\R^{+}\times\delta$.  To pin it down, we make
  $\Sigma$ disjoint from $Z_{1}\subset N_{1}$ and have geometric
  intersection $1$ with $E_{2} \subset Z_{2}$. If we write $A =
  A(k,g)$ again for the anti-self-dual connection, then
  \[
       \theta_{k}(g) = \exp \frac{1}{2}\int_{\Sigma} F_{A^{\ttr}}
  \]
  is a model for the map $\theta_{k}$ as a map from the circle. (The
  factor of $1/2$ is there because of the relationship between $A$ and
  $A^{\ttr}$.)  When $T=-\infty$, the surface $\Sigma$ is contained in
  the piece $T_{1}^{*}$ on which the connection $A^{\ttr}$ has become
  flat. So with this model, $\theta_{k}(-\infty) = 1$. When
  $T=+\infty$, the surface $\Sigma$ decomposes into two pieces: one is
  a cylinder contained in $T_{2}^{*}$, which contributes nothing to
  the integral; and the other is a disk $\Delta$ with cylindrical end contained
  in $Z_{2}^{*}\cong \mCP \setminus B^{4}$.  The curvature
  $F_{A^{\ttr}}$ has exponential decay on the end of $Z_{2}^{*}$ and
  its integral on $\Delta$ is equal to its integral on any compact
  surface $\Delta'\subset Z_{2}^{*}$ having the same intersection with
  $E_{2}$. Since $\Delta$ has intersection $1$ with $E_{2}$ and
  $c_{1}(\tspin_{k})$ evaluates to $-(2k+1)$ on $E_{2}$, we have
  \[
                    \int_{\Delta} F_{A^{\ttr}} = (2\pi/i) (2k+1).
  \]
   So $\theta_{k}(\infty) = -1$, and we have the result.
\end{proof}

For each integer $i<0$, we can define a map
\[
          \check{L}[i] : \Cto_{\bullet}(Y_{1}) \to \Cto_{\bullet}(Y_{1}) 
\]
by repeating the definition of $\check{L}$ above, but replacing
$\bar{n}_{u}$ in the formulae by the basis vector $e^{1}_{i}$. Because
there are no differentials on $R_{1}$, the map $\check{L}[i]$ is a
chain map, and from the formula for $\bar{n}_{u}$, we have
\[
\begin{aligned}
    \check{L} &= \sum_{k\in\Z} a_{k} \, \check{L}[-1-k(k+1)/2] \\
    &= \sum_{k\ge 0} \check{L}[-1-k(k+1)/2],
\end{aligned}
\]
where in the second line we have used Lemma~\ref{lem:KeyOdd}.
Proposition~\ref{prop:L-Calc} now follows from:

\begin{lemma}
\label{lem:Key-U}
    The map $\check{L}[i] : \Cto_{\bullet}(Y_{1}) \to
    \Cto_{\bullet}(Y_{1})$ for $i<0$ gives rise to the map
    $\Hto_{\bullet}(Y_{1}) \to \Hto_{\bullet}(Y_{1})$ given by
    multiplication by $U^{-i-1}$.
\end{lemma}

\begin{proof}
    We use the fact that the manifold $U_{1}$ (whose moduli spaces
    define $\check{L}$) can be realised as the complement in
    $[-1,1]\times Y_{1}$ of the tubular neighborhood a curve $K$: thus
    \[
             [-1, 1] \times Y_{1} = U_{1} \cup_{R_{1}} N(K),
    \]
    where $N(K) \cong S^{1}\times B^{3}$ and $R_{1}$ is the oriented
    boundary of $N(K)$. Referring to the definition of the action of
    $U^{p}$ from Section~\ref{subsec:ModuleStructure}, we choose $p$
    basepoints $w_{1}, \dots w_{p}$ in the interior of $N(K)$, and use
    these together with the cylindrical cobordism $[-1, 1]\times
    Y_{1}$ to define a chain map
    \begin{equation}\label{eq:mto-U}
       \mto([-1,1]\times Y_{1}, \{ w_{1}, \dots, w_{p}\})
        : \Cto_{\bullet}(Y_{1}) \to \Cto_{\bullet}(Y_{1}) 
    \end{equation}
    which induces the map $U^{p}$ on $\Hto_{\bullet}(Y_{1})$.

    We can use any Riemannian metric on the cylinder in the construction of this chain
    map. We choose a metric in which $N(K)$ has positive scalar
    curvature, the metric on $R_{1}$ is standard, and $R_{1}$ has a
    product neighborhood. We then consider the family of metrics
    parametrized by
    $Q = [0,\infty)$ obtained by inserting a cylinder $[-T,T]\times
    R_{1}$.  By an argument similar to our previous analysis, we
    obtain in this way a chain-homotopy between the chain map
    \eqref{eq:mto-U} and the chain map
    \[
              \sum_{j<0} b_{j} \, \check{L}[j], 
    \]
    where
    \[
                 \sum_{j<0} b_{j} e^{1}_{j} = \bar{n}_{u}(N(K))
    \]
    is the element of $\Cfrom_{\bullet}(R_{1}) \cong
    \Hfrom_{\bullet}(R_{1})$ obtained by counting points in moduli
    spaces on $N(K)^{*}$ with $p$ base-points. That is,
    \[
             b_{j} = | M( N(K), \fa^{1}_{j}) \cap V_{1} \cap \cdots
             \cap V_{p} | \bmod{2},
    \]
    or zero if the intersection is not zero-dimensional. An
    examination of dimensions shows that the only contributions occur
    when $j = -1 -p$. The moduli consists of reducibles, so the
    calculation of $b_{j}$ is straightforward: we have $b_{j}=1$ when
    $j=-1-p$.  Thus the sum above has just one term, and the map
    $U^{p}$ is equal to the map arising from the chain map
    $\check{L}[-1-p]$.
\end{proof}

With the verification of Proposition~\ref{prop:L-Calc}, the proof of
Theorem~\ref{thm:TwoHandleSurgeries-A} is complete for the case of
$\Hto_{\bullet}$. The other two case have similar proofs. In the case
of $\Hred_{\bullet}$, the formulae are considerably simpler. The
exactness in the case of $\Hfrom_{\bullet}$ could also be deduced from
the other two cases, by chasing the square diagram in which the
columns are the $i,j,p$ exact sequences and the rows are the surgery
cobordism sequences.

\subsection{Local coefficients}
\label{sec:Twisted}

We describe here a refinement of the long exact sequence, with local
coefficients. As in Section~\ref{subsec:Surgery-A}, 
we consider a $3$-manifold $M$  with torus boundary, and let
$\gamma_{0}$, $\gamma_{1}$, $\gamma_{2}$ be three oriented simple
closed curves on $\partial M$ with algebraic  intersection
\[
(\gamma_{0}\cdot\gamma_{1}) = (\gamma_{1}\cdot \gamma_{2})  =
(\gamma_{2}\cdot \gamma_{0}) = -1.
\]
We again write
$W_n\colon Y_n\longrightarrow
Y_{n+1}$ for the $2$-handle cobordisms.

The interior $M^o$ can be viewed as an open subset of $Y_i$
for all $i$. Fix local coefficient systems $\Gamma_i$ over $Y_i$
which are supported in $M^o\subset Y_i$, in the sense of
Definition~\ref{defn:Support}. 
Moreover, the set $[0,1]\times M^o$ can be viewed as an open subset of
$W_n$ (where here $\{0\}\times M^o\subset Y_n$ and $\{1\}\times
M^o\subset Y_{n+1}$).  Let $$\Gamma_{W_n}\colon
\Gamma_{n}\longrightarrow \Gamma_{n+1}$$ be a $W_n$-morphism of the
local system which is supported in $[0,1]\times M^o$.

\begin{theorem}
\label{thm:TwoHandleSurgeries-Local}
	Let $\Gamma_n$ be local systems on the $Y_n$, and let
        $\Gamma_{W_n} :\Gamma_{n}\to\Gamma_{n+1}$ be morphisms of
	local systems. Suppose these satisfy the support condition
        just described.
	Then, the induced maps with local coefficients 
	$\Fto_n={\Hto}(W_n;\Gamma_{W_n})$ fit into a long 
	exact sequence of the form
	\[ \cdots \longrightarrow
		\tHto_{\bullet}(Y_{n-1};\Gamma_{n-1})
		\stackrel{\Fto_{n-1}}
		{\longrightarrow} 
		\tHto_{\bullet}(Y_{n};\Gamma_{n})
	\stackrel{\Fto_n}{\longrightarrow}
	\Hto_{\bullet}(Y_{n+1};\Gamma_{n+1})\longrightarrow \cdots.  \] 
	There are also corresponding long exact sequences for the other
	two variants of Floer homology.
\end{theorem}

\begin{proof}
To prove exactness, we once again appeal to
Lemma~\ref{lemma:MCone}. Indeed, the homotopies ${\check H}_n'\colon
\Cto(Y_n;\Gamma_n) \longrightarrow
\Cto(Y_{n+2};\Gamma_{n+2})$ are constructed as before, only now 
the entries contain also the $X_n$-morphisms gotten by composing
the morphisms $\Gamma_{W_n}$ and $\Gamma_{W_{n+1}}$
(we add the primes to distinguish the homotopies here from 
the ones appearing in the discussion in Subsection~\ref{subsec:FirstHomotopy}).
Observe that
this composite morphism of local systems
is supported in the complement of $Z_n$
(so as in Lemma~\ref{lem:R1-lemma1}, the contributions from
$Z_i$ still drop out in pairs). 
In fact, the proof of Proposition~\ref{prop:Homotopy2}
carries over, as well, since the triple-composite morphism
of local systems is supported in a complement of
the $N_i$. 

In the same way, the map $\Lto'$ is seen to be multiplication
gotten by multiplying the chain map induced by the triple-composite 
$\Gamma_{W_{n+2}}\circ \Gamma_{W_{n+1}}\circ \Gamma_{W_n}$
(which is supported in the complement of $N_n$) 
by the power series $\Lto$.
In particular, the map
$\Lto'$ is a quasi-isomorphism, too.
\end{proof}

Specializing to the local system determined by cycles
(Example~\ref{ex:LocalSystem}), we get the following:

\begin{cor}
	\label{cor:TwoHandleTwist} 
	Let $Y_0$, $Y_1$, $Y_2$ be as above, with
	the additional property that $H_1(Y_0;\Z)\cong \Z$ and
	$H_1(Y_1;\Z)=H_1(Y_2;\Z)=0$. Fix a cycle $\eta$ in $M^o$ which
	generates the image of $H_1(M;\Z)$ in $H_1(M;\R)$.  In this
	case, we have a long exact sequence \[ \cdots \longrightarrow
	\tHto_{\bullet}(Y_{-1})\otimes \Kfield
	\stackrel{{\Fto}_{-1}}{\longrightarrow}
	\tHto_{\bullet}(Y_{0};\Gamma_{\cycle_0})
	\stackrel{\Fto_{0}}{\longrightarrow}
	\Hto_{\bullet}(Y_{1})\otimes \Kfield
	\stackrel{{\Fto}_{1}}{\longrightarrow} \cdots.  \] in which
	the map ${\Fto_1}$ can be expressed in terms of the usual maps
	induced by cobordisms, by the following formula:
	\[{\Fto}_1=\sum_{\fs\in\SpinC(W_1)} \mu({\langle
	c_1(\spinc),[h]\rangle})\cm {\Hto}(W_1,\spinc),\] where
	$[h]\in H_2(W_1;\Z) \cong \Z$ is a generator.
\end{cor}

\begin{proof}
We apply Theorem~\ref{thm:TwoHandleSurgeries-Local} in the following
setting.  We let $\Gamma_n$ be the local system on $Y_n$ induced by
the chain $\eta\subset M^o\subset Y_n$. Indeed, in the cobordisms
$W_n\colon Y_n\to Y_{n+1}$, we choose two-chains $\nu_n$ which are
products $\nu_n=[0,1]\times \eta\subset [0,1]\times M^o \subset W_n$.
The chain $\nu_n$ induces a $W_n$-morphism of the local system
$\Gamma_{W_n,\nu_n}$ which is supported in $[0,1]\times M^o$.

Recall that the isomorphism class of $\Hto(Y;\Gamma_{\cycle})$ depends
only on the homology class of $\cycle$. In fact, since both $Y_{1}$
and $Y_2$ are homology three-spheres, the cycle is
null-homologous, so it follows at once that for $i=1,2$,
\begin{equation}
\label{eq:CanonIdent}
\Hto(Y_i;\Gamma_{\cycle_i})\cong \Hto(Y_i)\otimes \Kfield.
\end{equation}

We argue that the chain $\nu_1\subset W_1$ represents a generator of
$H_2(W_1,\partial W_1;\Z)\cong \Z$. To see why, recall that inside $Y_1$,
$\gamma_0$ can be viewed as a knot, with a Seifert surface
$\Sigma$. Pushing the interior of the Seifert surface into
$[0,1]\times Y_1$, which we then cap off inside the added two-handle,
we obtain a generator $[{\widehat \Sigma}]$ for $H_2(W_1;\Z)$.
On the other hand, our chain $\nu_1$ is gotten by $[0,1]\times
\cycle$, and $\cycle$ links the knot $\gamma_0$ once. Thus, it follows
that the oriented intersection number of ${\widehat \Sigma}$ with
$\nu_1$ is $\pm 1$, and hence $[\nu_1]$ is also a generator of
$H_2(W_1,\partial W_1)$.  Now, for each $\SpinC$ structure on $W_1$,
we see that composing with the identification
from Equation~\eqref{eq:CanonIdent}, we see that
the map 
$\HFto(W_1;\Gamma_{W_1,\nu},\spinc)$ is
identified with
$\mu({\langle c_1(\spinc),[\nu]]\rangle}) \cm \Hto(W_1,\spinc)$,
where here $[\nu]\in H_2(W_1;\Z)$ is the unique homology class
corresponding to $\nu_1$.
The result now follows.
\end{proof}

%
%
%
%

\section{Proof of the non-vanishing theorem}
\label{sec:NonVanishingProof}

\subsection{Statement of the sharper result}
\label{subsec:SharpNonVanishStatement}

We now turn to the proof of Theorem~\ref{thm:NonVanishing}. There is a
sharper version of this non-vanishing theorem, which involves the
Floer groups with local coefficients. 

\begin{theorem}
\label{thm:NonVanishing-twist}
    Suppose $Y$ admits a taut foliation $\Fol$ and is not
    $S^{1}\times S^{2}$. Let $\eta$ be a $C^{\infty}$ singular
    $1$-cycle in $Y$ whose homology class $[\eta]$ satisfies
    \[
               \mathrm{P.D.}[\eta] = [\omega] + t \,e(\Fol) \in
               H^{2}(Y;\R)
    \]
    where $\omega$ is closed $2$-form which is positive on the leaves
    of $\Fol$ and $t\in \R$. Then the image of the map
    \[
                 j_{*} : \Hto_{k}(Y;\Gamma_{\eta}) \to
                 \Hfrom_{k}(Y;\Gamma_{\eta})
    \]
    is non-zero, where $k\in J(Y)$ is the homotopy class of the
    $2$-plane given  by the tangents planes to $\Fol$
\end{theorem}

As an application, we consider the case of the manifold $Y =
S^{3}_{0}(K)$, where $K\ne U$. By the results of \cite{Gabai} again, this manifold has
a taut foliation $\Fol$; and if $\omega$ is closed and positive on
the leaves, then the cohomology class $[\omega]$ will be non-zero,
because Gabai's foliation has a compact leaf. We therefore have the
following corollary. Unlike Corollary~\ref{cor:NonVanishing}, this
result applies also to genus one knots:

\begin{cor}
     \label{cor:NonVanishingTwist} Suppose $K \ne U$, let $g$ be the
     Seifert genus of $K$, and let $\spinc$ be a $\SpinC$ structure on
     $S^{3}_{0}(K)$ for which $c_{1}(\spinc)$ is $2g-2$ times a
     generator for $H^{2}(S^{3}_{0}(K);\Z)$. Then there is a $k\in
     J(Y, \spinc)$ such that the image of 
     \[ j_{*} :
	\Hto_{k}(S^{3}_{0}(K);\Gamma_{\eta}) \to
	\Hfrom_{k}(S^{3}_{0}(K);\Gamma_{\eta}) 
     \] 
     is non-zero whenever
     the homology class $[\eta]$ is non-zero. By contrast, $j_{*}$ is
     zero for $S^{3}_{0}(U)$, for all $\eta$. \qed
\end{cor}
                                                                                
The proof of the theorem is based on the results of
\cite{KMcontact}. Let $X$ be a compact oriented $4$-manifold with
oriented boundary $Y$. We assume $X$ is connected, but may allow $Y$
to be disconnected. Let $\xi$ be an oriented contact structure on $Y$,
compatible with the orientation of $Y$. If $\alpha$ is a $1$-form on
$Y$ whose kernel is the field of $2$-planes $\xi$, then the
orientation condition can be expressed as the condition that $\alpha\wedge
d\alpha>0$.  Let $\spinc$ be the $\SpinC$ structure on $Y$ determined
by $\xi$, and let $\spinc_{X}$ be any extension of $\spinc$ to
$X$. Note that the space of isomorphism classes of such extensions,
denoted by $\SpinC(X,\xi)$, is
an affine space for the group $H^2(X,Y;\Z)$. Generalizing the monopole
invariants of closed $4$-manifolds, the paper
\cite{KMcontact} defines an invariant $\km(X, \xi, \spinc_{X})$ associated to
such data. Neglecting orientations, we can take it to be an
element of $\Field = \Z/2$.  We review some of the properties of this
invariant.

The $2$-plane field $\xi$ picks out not just a $\SpinC$ structure
$\spinc$ on $Y$, but also a preferred nowhere-vanishing section
$\Phi_{0}$ of the spin bundle $S\to Y$. When $\spinc_{X}$ is given, we can
interpret $\Phi_{0}$ as a section of $S^{+}_{X}|_{Y}$, and there is a
relative second Chern class (or Euler class) which we use as the
definition of 
$\gr(X, \xi, \spinc_{X})$:
\[
          \gr(X,\xi,\spinc_{X}) = \langle c_{2}(S^{+}_{X}, \Phi_{0}) , [X,\partial
          X]\rangle \in \Z.
\]
The condition $\gr(X,\xi,\spinc_{X})=0$ is equivalent to the existence of
an almost complex structure on $X$ for which the $2$-planes $\xi$ are
complex and such that the associated $\SpinC$ structure is $\spinc_{X}$.

\begin{theorem}[\cite{KMcontact}]\label{thm:KMCont}
    \ \par
    \begin{enumerate}
        \item The invariant $\km(X,\xi,\spinc_{X})$ is non-zero only if
        $\gr(X,\xi,\spinc_{X}) = 0$, and vanishes for all but finitely
        many $\spinc_{X}\in\SpinC(X,\xi)$.

        \item Suppose $X$ carries a symplectic form $\omega$ that is
        positive on the oriented $2$-plane field $\xi$. Let
        $\spinc_{\omega}\in\SpinC(X,\xi)$ 
	be the canonical $\SpinC$ structure which
        $\omega$ determines on $X$. Then
        \[
                  \km(X,\xi,\spinc_{\omega}) = 1.
        \]

        \item Let $\omega$ and $\spinc_{\omega}$ be as in the previous
        item, let $e \in H^{2}(X,\partial X;\Z)$, and let $\spinc_{\omega}
        + e\in\SpinC(X,\omega)$ 
	denote the $\SpinC$ structure with spin bundle $S =
        S_{\omega}\otimes L$, where $L$ is the line bundle,
        trivialized on $\partial X$, with relative first Chern class
        $c_{1}(L)=e$. Then, if $\km(X,\xi,\spinc_{\omega}+e) \ne 0$, it
        follows that
        \[
              \langle e \smile \omega, [X,\partial X] \rangle \ge 0,
        \]
        with equality only if $e=0$. \qed
    \end{enumerate}
\end{theorem}

We combine the individual invariants $\km(X,\xi,\spinc_{X})$ into a
generating function. Recall that $\Kfield$ is the field of
fractions of $\Field[\R]$, and that $\mu: \R \to \Kfield^{\times}$ is
the canonical homomorphism. Given a reference $\SpinC$ structure
$\spinc_{0}$ on $X$ extending $\spinc$, we define a function
\[
   \km^{*}(X,\xi,\spinc_{0}) : H_{2}(X,Y;\R) \to \Kfield
\]
by the formula
\[
    \km^{*}(X,\xi,\spinc_{0})(h) = \sum_{e} \km(X,\xi,\spinc_{0}+e) \, \mu \langle
    2e, h \rangle.
\]
As a corollary of Theorem~\ref{thm:KMCont} we have:

\begin{cor}
\label{cor:KMCont}
If $X$ carries
a symplectic form $\omega$ positive on $\xi$ then
\[
\km^{*}(X,\xi,\spinc_{\omega})(\mathrm{P.D.}[\omega])\ne 0.
\]
Further, if the intersection form on $H^{2}(X,\partial X)$ is trivial then
\[
\km^{*}(X,\xi,\spinc_{\omega})(\mathrm{P.D.}([\omega]+t\,c_{1}(\spinc_{\omega})))\ne 0
\]
for all $t\in \R$.  
\end{cor}

\begin{proof}
  The first statement is an immediate consequence of Theorem~\ref{thm:KMCont}.
  For the second statement, we note that the condition
  $\gr(X,\xi,\spinc_{\omega}+e) = 0$ is equivalent to
  \[
             e \smile (e  + c_{1}(\spinc_{\omega})) = 0,
  \]
  or simply to $e\smile c_{1}(\spinc_{\omega})=0$ when the cup product
  on the relative cohomology is zero; so
  \[
 \km^{*}(X,\xi,\spinc_{\omega})(\mathrm{P.D.}([\omega]+t\,c_{1}(\spinc_{\omega})))
    = \km^{*}(X,\xi,\spinc_{\omega})(\mathrm{P.D.}[\omega]).
  \]
\end{proof}

\subsection{Construction of the invariants of $(X,\xi)$}

We review the construction of the invariant $\km(X,\xi,\spinc_{X})$
from \cite{KMcontact}.  Let
$\alpha$ again be the $1$-form defining the $2$-plane field $\xi$ on
$Y$, 
let $\omega_{0}$ be the symplectic form $d(t^{2}\alpha/2)$ on
$[1,\infty)\times Y$, and let $g_{0}$ be a compatible metric.
Attaching $[1,\infty)\times Y$ to $X$, we obtain a complete Riemannian manifold
$X^{+}$ with an expanding conical end. On the end, there is a
canonical $\SpinC$ connection $A_{0}$ and spinor $\Phi_{0}$ of unit
length. We choose a reference $\SpinC$ structure $\spinc_{0}$ on $X$,
extending the $\SpinC$ structure on the end, and extend $A_{0}$ and
$\Phi_{0}$ arbitrarily.  For a pair $(A,\Phi)$ consisting of a
$\SpinC$ connection in $\spinc_{X}$ and a section of $S^{+}_{X}$, we
consider the equations on $X^+$:
\begin{equation}\label{eq:XplusSW}
	\begin{aligned}
		\frac{1}{2}\rho(F^{+}_{ A^{\ttr}}) -
                (\Phi\Phi^{*})_{0} &=
                    \frac{1}{2}\rho(F^{+}_{ A_{0}^{\ttr}}) -
                (\Phi_{0}\Phi_{0}^{*})_{0} + \epsilon,
                \\
		D^{+}_{A} \Phi &= 0,
	\end{aligned}
\end{equation}
where $\epsilon$ is a perturbation term: an exponentially decaying
section of $i\mathfrak{su}(S^{+}_{X})$.  There is a moduli space
$M(X^{+}, \spinc_{X})$ consisting of all gauge-equivalence classes of
solutions $(A,\Phi)$ on $X^{+}$ which are asymptotically equal to
$(A_{0},\Phi_{0})$, in that
\begin{equation}\label{eq:AFAK-asymptotics}
\begin{aligned}
    A - A_{0} &\in L^{2}_{k} \\
    \Phi - \Phi_{0} &\in L^{2}_{k,A_{0}}.
\end{aligned}
\end{equation}
For generic $\epsilon$, this moduli space is a smooth manifold, and
\[
           \mathrm{dim} M(X^{+},\spinc_{X}) = \gr(X,\xi,\spinc_{X})
\]
if the moduli space is non-empty. Note that the asymptotic conditions
mean that $\Phi$ is non-zero, so there are no reducible solutions
in the moduli space. The moduli space is compact, and
$\km(X,\xi,\spinc_{X})$ is defined as the number of points in the moduli
space, mod 2, or as zero if the dimension is positive.

\subsection{Floer chains from contact structures}

Let $Z$ now denote the half-infinite cylinder $\R^{+} \times (-Y)$,
which has oriented boundary $\{0\} \times Y$. Given a contact
structure $\xi$ defined by a $1$-form $\alpha$ as above, we again form
the symplectic cone $[1,\infty)\times Y$, with oriented boundary $\{1\}\times
(-Y)$, and attach this to $Z$. The result is a complete Riemannian
manifold $Z^{+}$ with one cylindrical end and one expanding conical end; the
latter carries the symplectic form $\omega_{0}$.

On $Z^{+}$, we write down monopole equations which resemble the
equations \eqref{eq:XplusSW} on the conical end and resemble the
perturbed equations $\SW_{\pertX}=0$ on the cylindrical end. A
convenient way to make this construction is  as a fiber product. To do
this, we
choose a regular perturbation $\pertY$ for the equations on $Y$, and
let $\pertX$ be a $t$-dependent perturbation on the cylinder $Z$ that
is equal to $\pertY$ on the end and is zero near the boundary. For
each critical point $\fa$ in $\bonf^{\sigma}(-Y)$, there
is a moduli space
\[
    M(\R^{+}\times (-Y), \fa) \subset \bonf^{\sigma}(\R^{+}\times
    (-Y)).
\]
This is a Banach manifold with a restriction map
\[
    r_{0} :  M(\R^{+}\times (-Y), \fa) \to
    \bonf^{\sigma}(\{0\}\times(-Y)) =
    \bonf^{\sigma}(Y).
\]
There is also the moduli space $M_{1}$ of solutions  to the
equations \eqref{eq:XplusSW} on the conical manifold $[1,\infty)\times
Y$, with asymptotic conditions \eqref{eq:AFAK-asymptotics}. We choose
the term $\epsilon$ so that the right-hand side of the first equation
vanishes near the boundary $\{1\}\times Y$. This is
another Banach manifold; and because the solutions are irreducible,
there is a restriction map to blown-up configuration space of the
boundary, because of a unique continuation argument:
\[
    r_{1} :  M_{1} \to \bonf^{\sigma}(\{1\}\times Y) =
    \bonf^{\sigma}(Y).
\]

\begin{defn}
    We define the moduli space $M(Z^{+}, \fa)$ as the fiber product of
    the maps $r_{0}$ and $r_{1}$.
\end{defn}

Although the fiber product makes a convenient definition, we can also
regard $M(Z^{+}, \fa)$ as a subspace of
$\bonf^{\sigma}_{\loc}(Z^{+})$. Given  $\gamma$ in $M(Z^{+}, \fa)$,
we can define a path $\check{\gamma}$ in $\bonf^{\sigma}(Y)$ by
restricting $\gamma$ to the slices $\{t\}\times Y$, first in the
cylindrical end then in the conical end. If $\gamma$, $\gamma'$ are
two solutions, then the corresponding paths $\check{\gamma}$ and
$\check{\gamma}'$ both have limit point $\fa$ on the cylindrical
end and have the same asymptotics on the conical end; so there is a
well-defined difference element in $\pi_{1}(\bonf^{\sigma}(Y), \fa)$.
In this way, we partition $M(Z^{+},\fa)$ into components of different
topological type:
\[
            M(Z^{+},\fa) = \bigcup_{z} M_{z}(Z^{+},\fa).
\]
Once again, we can count points in zero-dimensional moduli spaces, to
define:
\[
            m_{z}(Z^{+},\fa) = 
             \begin{cases}
                | M_{z}(Z^{+},\fa)| \bmod{2},
                &\text{if $\dim M_{z}(Z^{+},\fa)=0$,}\\
                0, &\text{otherwise.}
             \end{cases}
\]

The compactness results of \cite{KMcontact} tell us that
$m_{z}(Z^{+},\fa)$ is zero for all but finitely many $\fa$ and $z$.
Because $M_{z}(Z^{+},\fa)$ consists only of irreducibles, $\fa$ must
be either irreducible or boundary-stable on $-Y$ if the moduli space is to be
non-empty. (Note that the notion of ``boundary-stable'' for a critical
point $\fa$ depends on the orientation: $\Crit^{s}(-Y)$ is the same as
$\Crit^{u}(Y)$.) Taking the irreducible and boundary-stable elements
in turn, we define an element of the complex $\Cto_{*}(-Y) = C^{o}(-Y)
\oplus C^{s}(-Y)$, by
\begin{equation}\label{eq:Check-psi}
        \check{\psi} = (\psi^{o} , \psi^{s}),
\end{equation}
where
\[
       \psi^{o} = \sum_{\fa\in \Crit^{o}} \sum_{z} m_{z}(Z^{+},\fa)\, e_{\fa}
\]
and
\[
       \psi^{s} = \sum_{\fa\in \Crit^{s}} \sum_{z} m_{z}(Z^{+},\fa)\, e_{\fa}.
\]
We have:

\begin{lemma}
    The element $\check{\psi}$ in $\Cto_{*}(-Y)$ is closed: that is,
    $\dto\check{\psi} = 0$.
\end{lemma}

\begin{proof}
    As usual, this is proved by  counting the boundary points of
    $1$-dimensional moduli spaces $M_{z}(Z^{+},\fb)$, augmented by the
    observation that there are no reducible solutions in these moduli
    spaces. Specifically, counting boundary points in the case that
    $\fb$ is irreducible gives the identity
    \[
                     \doo  \psi^{o} + \duo \dsu \psi^{s} = 0,
    \]
     while the case that $\fb$ is boundary-stable provides the
     identity
     \[
            \dos \psi^{o} + \dss \psi^{s} + \dIus \dsu \psi^{s} = 0.
     \]
     Together these tell us that $\dto\check{\psi}=0$.
\end{proof}

Next we extend the construction of $\check{\psi}$ to the Floer complex
with local coefficients.
Let $\eta$ be a $C^{\infty}$ real $1$-cycle in $Y$, and let $\eta^{+}$
be the corresponding non-compact $2$-chain in $Z^{+}$. Like $Z^{+}$,
the $2$-chain $\eta^{+}$ has one cylindrical end and one expanding
conical end; it is oriented
so that the cylindrical end coincides with $-(\R^{+}\times\eta)$.
Extend the connection $A_{0}$ from the conical end to all
of $Z^{+}$ in such a way that it is translation-invariant and in
temporal gauge on the cylindrical end.
Let $\gamma$ belong to $M_{z}(Z^{+},\fa)$, and let $A$ be
the corresponding $\SpinC$ connection on $Z^{+}$. The integral
\[
            f(z) = (i/2\pi) \int_{\eta^{+}} (F_{A^{\ttr}} -
            F_{A_{0}^{\ttr}})
\]
is finite, and depends on $\gamma$ through only its homotopy class
$z$. Let $\Gamma_{-\eta}$ be the local system on $\bonf^{\sigma}(-Y)$
defined in Example~\ref{ex:LocalSystem}, and denote the generator $1$
in $\Gamma_{\fa} \cong \Kfield$ by $e_{\fa}$.
Then we can
define an element
\[
        \check{\psi}_{\eta} = (\psi^{o}_{\eta},
                \psi^{s}_{\eta}) 
                \in \Cto(-Y;\Gamma_{-\eta})
\]
by the formulae
\[
\begin{aligned}
\psi^{o}_{\eta} &= \sum_{\fa\in \Crit^{o}} \sum_{z}
       m_{z}(Z^{+},\fa) \mu(f(z)) \, e_{\fa} \\
       \psi^{s}_{\eta} &= \sum_{\fa\in \Crit^{s}} \sum_{z}
       m_{z}(Z^{+},\fa) \mu(f(z)) \, e_{\fa}.
       \end{aligned}
\]
As in the previous case, we have:
\begin{lemma}
    The element $\check{\psi}_{\eta}$ in $\Cto_{*}(-Y;\Gamma_{-\eta})$ is
    closed. \qed
\end{lemma}

\subsection{Proof of Theorem~\ref{thm:NonVanishing-twist}.}

We recall Eliashberg and Thurston's construction
\cite{EliashbergThurston}, whereby a taut
foliation on $Y$ leads to a symplectic form $\omega_{W}$ on the cylinder
\[
        W = [-1,1] \times Y
\]
and contact structures $\xi_{\pm}$ on the boundary components. Let
$\alpha$ be a $1$-form defining the tangents to the foliation $\Fol$,
and let $\omega$ be a closed $2$-form positive on the leaves. Set
\[
        \omega_{W} = d(t\alpha) + \omega
\]
on $W$. This form is symplectic. According to
\cite{EliashbergThurston}, there exist smooth contact structures
$\xi_{+}$ and $\xi_{-}$, compatible with the orientations of $Y$ and
$-Y$ respectively, which are $C^{0}$ close to the tangent plane field
of the foliation. We regard these as contact structures on the two
boundary components $\{1\}\times Y$ and $\{-1\}\times Y$ of $W$; the
$C^{0}$-close condition means that $\omega_{W}$ will be positive on
these $2$-plane fields.

If we regard $W$ as a manifold with a contact structure $\xi =
(\xi_{-}, \xi_{+})$ on its boundary $\{-1,1\}\times Y$, then we can
construct the invariants
\[
          \km(W,\xi, \spinc_{W})
\]
for $\SpinC$ structures $\spinc_{W}$ extending the standard $\SpinC$
structure determined by the $2$-plane field on the boundary, as
above. On the other hand, the contact structure $\xi_{+}$ provides a
cycle
\[
           \check{\psi}( \xi_{+}) = \bigl(\psi^{o}(\xi^{+}),
           \psi^{s}(\xi^{+})\bigr) \in \Cto_{*}(-Y),
\]
and from $\xi_{-}$ we similarly obtain a cycle
\[
           \check{\psi}( \xi_{-}) = \bigl(\psi^{o}(\xi^{-}),
           \psi^{s}(\xi^{-}) \bigr)\in \Cto_{*}(Y).
\]
Let
\[
\begin{aligned}
        \check{\Psi}(\xi_{+})& \in \Hto_{*}(-Y) &
        \check{\Psi}(\xi_{-})& \in \Hto_{*}(Y)
\end{aligned}
\]
be the homology classes of these.

\begin{prop}
        We have the pairing formula
        \[
                  \sum_{\spinc_{W}} \km( W, \xi, \spinc_{W}) =
                  \bigl\langle
                       j_{*}  \check{\Psi}(\xi_{+}),
                       \check{\Psi}(\xi_{-})
                  \bigr\rangle_{\dual},
        \]
        where $\langle\;\mathord{-}\;, \;\mathord{-}\;\rangle_{\dual}$
        is the duality pairing
        \[
                    \Hfrom_{*}(-Y) \otimes \Hto_{*}(Y) \to \Field,
        \]
        from Section~\ref{subsec:Duality}.
\end{prop}

\begin{proof}
        Let $Z^{+}(\xi_{+})$ be the manifold with one cylindrical end
        and one conical end, obtained by applying the
        construction of the previous subsection to $\xi_{+}$ on $Y$,
        and let $Z^{+}(\xi_{-})$ be constructed similarly using
        $\xi_{-}$. As above, we have counting-invariants
        \[
        \begin{aligned}
        m_{z}(Z^{+}(\xi_{+}), \fa) &\in \Field, & (\fa &\in
               \Crit^{o}(-Y) \cup \Crit^{s}(-Y)), \\
               m_{z}(Z^{+}(\xi_{-}), \fa) &\in \Field, & (\fa &\in
               \Crit^{o}(Y) \cup \Crit^{s}(Y)).
               \end{aligned}
        \]
        In the case of the first of these two, we can regard $\fa$ as
        a critical point in $\bonf^{\sigma}(Y)$ via the identification
        \[
               \Crit^{o}(-Y) \cup \Crit^{s}(-Y)) =
               \Crit^{o}(Y) \cup \Crit^{u}(Y)).
        \]
        If we unravel the pairing on the right-hand side of the formula
        in the proposition, we find it is equal to
        \begin{multline}
        \label{eq:ET-split-sum}
               \sum_{\fa\in \Crit^{o}(Y)}   m_{z_{1}}(Z^{+}(\xi_{-}), \fa)
               m_{z_{2}}(Z^{+}(\xi_{+}), \fa)
              \\ +   \sum_{\fa\in \Crit^{s}(Y)}
                \sum_{\fb\in \Crit^{u}(Y)}
                m_{z_{1}}(Z^{+}(\xi_{-}), \fa)
                n_{z_{2}}(\fa,\fb)
               m_{z_{3}}(Z^{+}(\xi_{+}), \fb)
        \end{multline}
        where $n_{z}(\fa,\fb)$ counts unparametrized
        boundary-obstructed trajectories as in
        Section~\ref{subsec:ThreeMorse}. On the other hand, we can
        consider the $1$-parameter family of metrics on $W$
        parametrized by $Q=[0,\infty)$ in which the length of the
        cylinder is increased. There is a corresponding parametrized moduli space
        \[
                           M(W^{+}, \spinc_{W})_{Q}
        \]
        associated to the manifold $W^{+}$ with two conical ends. The
        map to $Q$ is proper, and there is a compactification
        $M(W^{+},\spinc_{W})_{\bar{Q}}$ over $\bar{Q} = [0,\infty]$,
        where at $\infty$ the manifold $W^{+}$ becomes the disjoint
        union of $Z^{+}(\xi_{-})$ an $Z^{+}(\xi_{+})$. If we look at
        the union of all $1$-dimensional moduli spaces $M(W^{+},
        \spinc_{W})_{\bar{Q}}$ and count the endpoints of these, then the
        contributions from endpoints lying over $0\in \bar{Q}$ is
        equal to the left-hand side in the proposition, while the
        contribution from the endpoints lying over $\infty$ is the sum
        \eqref{eq:ET-split-sum}.
        \end{proof}
    
We now reformulate this proposition for local coefficients.  Because
$\xi_{-}$ and $\xi_{+}$ are both $C^{0}$-close to $\Fol$, there is a
canonical choice of reference $\SpinC$ structure $\spinc_{0}$ on $W$,
and we can write an arbitrary
$\spinc_{W}$ as
\[
          \spinc_{W} = \spinc_{0} + e,
\]
for some $e\in H^{2}(W,\partial W;\Z)$.  Let $\eta$ be a $1$-cycle in
$Y$, and set
\[
       h_{\eta} = [-1,1] \times [\eta]
\]
regarded as a class in $H_{2}(W,\partial W;\R)$. Let
\[
 \km^{*}(W,
\xi, \spinc_{0}) : H_{2}(W,\partial W;\R) \to \Kfield
\]
be as in Section~\ref{subsec:SharpNonVanishStatement}.  The
construction of $\check{\psi}_{\eta}$ leads to classes
\[
\begin{aligned}
        \check{\Psi}_{\eta}(\xi_{+})& \in \Hto_{*}(-Y;\Gamma_{-\eta}) &
        \check{\Psi}_{\eta}(\xi_{-})& \in \Hto_{*}(Y;\Gamma_{\eta}).
\end{aligned}
\]
Just as in the case of coefficients $\Field$, we have:

\begin{prop}
\label{prop:ET-pairing}
        We have the pairing formula
        \[
                 \km^{*}( W, \xi, \spinc_{0})(h_{\eta}) =
                  \bigl\langle
                       j_{*}  \check{\Psi}_{\eta}(\xi_{+}),
                       \check{\Psi}_{\eta}(\xi_{-})
                  \bigr\rangle_{\dual},
        \]
        where $\langle\;\mathord{-}\;, \;\mathord{-}\;\rangle_{\dual}$
        is the duality pairing
        \[
                    \Hfrom_{*}(-Y;\Gamma_{-\eta}) \otimes
                    \Hto_{*}(Y;\Gamma_{\eta}) \to \Field,
        \]
        from Section~\ref{subsec:Duality}, and $j_{*}$ is the map
         $\Hto_{*}(-Y;\Gamma_{-\eta}) \to
         \Hfrom_{*}(-Y;\Gamma_{-\eta})$. \qed
\end{prop}

Now we conclude the proof of Theorem~\ref{thm:NonVanishing-twist}.
Suppose the class $[\eta]$ in $H_{1}(Y;\R)$ satisfies
\[
        \mathrm{P.D.}[\eta] = [\omega] + t\, c_{1}(\Fol),
\]
so that the class $h_{\eta}$ in $H_{2}(W,\partial W;\R)$ satisfies
\[
        \mathrm{P.D.}[h_{\eta}] = [\omega_{W}] + t\,
        c_{1}(\spinc_{\omega_{W}}).
\]
The intersection form on $H^{2}(W,\partial W)$ is trivial, so
Corollary~\ref{cor:KMCont} tells us that
$\km^{*}(W,\xi,\spinc_{0})(h_{\eta})$ is non-zero. (Note that
$\spinc_{0}$ and $\spinc_{\omega_{W}}$ are the same.) From the
proposition above, it follows that $j_{*}\check{\Psi}_{\eta}(\xi^{+})$
is non-zero; and in particular the image of $j_{*}$ is non-trivial in
$\Hfrom_{*}(-Y;\Gamma_{-\eta})$. The hypotheses of the theorem are
symmetrical with respect to orientation, so $j_{*}$ has non-zero
image also in $\Hfrom_{*}(Y;\Gamma_{\eta})$. \qed

\subsection{Application to knots of genus one.}
\label{subsec:GenusOneProof}

Using the non-vanishing theorem with local coefficients, we can now
complete the proof of Theorem~\ref{thm:RP3} for knots of genus 1, in
the case that the surgery coefficient is an integer. The arguments
of Section~\ref{sec:ProofOutline} continue to show that if $K$ is
$p$-standard then $K$ is $(p-1)$-standard, for integers $p\ge 2$. We
need therefore only prove the following result, by a method applicable
to genus $1$.

\begin{prop}
If $K$ is $p$-standard for $p=1$, then $K$ is the unknot.
\end{prop}

\begin{proof}
Consider the long exact sequence sequence with local coefficients,
in the form given in Corollary~\ref{cor:TwoHandleTwist}, applied to
$3$-manifolds $S^{3}$, $S^{3}_{1}(K)$ and $S^{3}_{0}(K)$:
	\[ \cdots \longrightarrow
	\tHto_{\bullet}(S^{3})\otimes \Kfield
	\stackrel{{\Fto}_{-1}}{\longrightarrow}
	\tHto_{\bullet}(S^{3}_{0}(K);\Gamma_{\cycle_0})
	\stackrel{\Fto_{0}}{\longrightarrow}
	\Hto_{\bullet}(S^{3}_{1}(K))\otimes \Kfield
		\stackrel{{\Fto}_{1}}{\longrightarrow} \cdots.  \] 
If $K$ is $1$-standard, then
the map $\Fto_1$ is given by multiplication by
\begin{equation*}
\sum_{n\geq 0} U^{n(n+1)/2}\cm (\mu({2n+1})+\mu({-2n-1}))
 ,
\end{equation*}
thought of as a map from $\Kfield[U^{-1},U]/\Kfield[U]$ ($\cong
\Hto(S^{3}_{1}(K))\otimes \Kfield\cong
\Hto(S^3)\otimes\Kfield$) to itself. 
Since the coefficient of
$U^0$ is a non-zero element of $\Kfield$, this map is an
isomorphism.  By the long exact sequence, we can conclude that
$\Hto(S^{3}_{0}(K);\Gamma_{\cycle_0})=0$. It follows now from
Corollary~\ref{cor:NonVanishingTwist} that $K$ is the unknot.
\end{proof}

\section{The case of non-integral $r$}
\label{sec:Rational}

We now turn to the proof of Theorem~\ref{thm:RP3} in the case where
$r$ is non-integral. Some of the results proved along the way apply in
more general settings, and will be used later. 

Our strategy here is to show that if there is an
orientation-preserving diffeomorphism $S^3_r(K)\cong S^3_r(U)$ for
$r>0$, then $K$ is $p$-standard where $p$ is the smallest integer
greater than $r$. In this way, we reduce to the case of integral $r$,
which is proved earlier in the paper. Of course, the case, where $r<0$
once again follows from the case where $r>0$, by reflecting the knot.

\begin{lemma}
\label{lemma:SOneSTwo}
Let $W\colon Y_1\to Y_2$ be a cobordism which contains a sphere
with self-intersection number zero 
$S\subset W$ which represents a non-trivial homology class in
$H_2(W,\partial W;\Z)$.  Then, for each $\SpinC$ structure
$\spinc\in\SpinC(W)$, the induced map on Floer homologies are trivial.
\end{lemma}

\begin{proof}
Since $S^1\times S^2$ admits a metric of positive scalar curvature, it
follows that if $\Hto(W,\spinc)$ is non-trivial, then $\langle
c_1(\spinc), [S]\rangle =0$. 

We must now prove
that even if $\langle c_1(\spinc),[S]\rangle =0$, then the 
induced map $\Hto(W,\spinc)$ is trivial.
To see this, we pass to the blow-up ${\widehat W}=W\# \mCP$.
Fix a two-sphere $E$ in ${\widehat W}$ supported in $\mCP$ with
square $-1$.  Given $\spinc\in\SpinC(W)$, there is a $\SpinC$
structure ${\widehat \spinc}\in\SpinC({\widehat W})$ 
which extends $\spinc$
and with the additional property
that $\langle c_1({\widehat \spinc}),[E]\rangle = +1$.  By the blow-up
formula, $\Hto(W,\spinc)=\Hto({\widehat W},{\widehat \spinc})
= \Hto({\widehat W},{\widehat \spinc}+\PD(E))$. Since $[S]-[E]$
can also be represented by a sphere with self-intersection number $-1$,
$\Hto({\widehat W},{\widehat \spinc}+\PD(E))=
\Hto({\widehat W},{\widehat \spinc}+\PD[S]))$. By repeatedly
using the same argument , we have that 
$\Hto(W,\spinc)=\Hto({\widehat W},{\widehat
\spinc}+ k\PD[S])$ for all $k\in\Z$. However, our hypotheses on $S$
ensure that $\{{\widehat \spinc}+k\PD[S]\}_{k\in\Z}$ is an infinite
collection of $\SpinC$ structures; but for any fixed $\eta\in
\Hto(Y_1)$, there can be only finitely many $\SpinC$ structures for
which the map $\Hto({\widehat W},\spinct)(\eta)$ is non-trivial.  It
now follows that $\Hto(W,\spinc)=0$.
Since $b^+(W)>0$, it follows
from Proposition~\ref{prop:Reducibles-II} that
the map $\Hred(W,\spinc)=0$ as well, and an easy diagram chase now
also shows that $\Hfrom(W,\spinc)=0$.
\end{proof}

\begin{prop}
\label{prop:SOneSTwo}
Let $M$ be a three-manifold with torus boundary, and choose oriented
curves $\gamma_1$ and $\gamma_2$ with $\gamma_1\cdot\gamma_2 = -1$.
Fix also a cycle $\eta\subset M$.  Letting $\gamma_3$ be a curve
representing $[\gamma_1]+[\gamma_2]$ and $\gamma_4$ be a curve
representing $[\gamma_1]-[\gamma_2]$.  Consider the surgery long exact
sequences
\[
\begin{gathered}
\cdots
\longrightarrow \Hto_{\bullet}(Y_1;\Gamma_\eta)
\stackrel{\Fto_{1}}{\longrightarrow}
\Hto_{\bullet}(Y_{2};\Gamma_\eta) 
\stackrel{\Fto_2}{\longrightarrow} \Hto_{\bullet}(Y_{3};\Gamma_\eta)\longrightarrow \cdots \\
\cdots\longrightarrow \Hto_{\bullet}(Y_2;\Gamma_\eta)
\stackrel{\Gto_2}\longrightarrow \Hto_{\bullet}(Y_{1};\Gamma_\eta) \stackrel{\Gto_1}{\longrightarrow}
\Hto_{\bullet}(Y_{4};\Gamma_\eta)\stackrel{\Gto_4}{\longrightarrow}\cdots
\end{gathered}
\]
where here $Y_i$ is obtained
from $M$ by filling $\gamma_i$
and $\Fto_i$ and $\Gto_i$ are maps induced
by the cobordisms 	
equipped with the product cycles $[0,1]\times \eta$, thought
of as supported in the complement of the two-handle additions,
as in Theorem~\ref{thm:TwoHandleSurgeries-Local}.
Then we have that $\Fto_1\circ \Gto_2=0=\Gto_2\circ \Fto_1$.
For composites of the maps belonging to 
the long exact sequences for the other two Floer homologies
$\Hred$ and $\Hfrom$, we have an analogous vanishing results.
(e.g. $\Ffrom_1\circ \Gfrom_2=0$).
\end{prop}

\begin{proof}
Let $A_1\colon Y_1\longrightarrow Y_2$ denote the cobordism inducing
the map $\Fto_1$, and and $B_2\colon Y_2\longrightarrow Y_1$ denote the
cobordism inducing the map $\Gto_2$, so that the union $B_2\cup_{Y_2}
A_1\colon Y_1\longrightarrow Y_1$ is a cobordism
(i.e. $\Fto_1=\Hto(A_1;\Gamma_{A_1,\nu})$ and
$\Gto_2=\Hto(B_2;\Gamma_{B_2,\nu})$, where here the chains  $\nu$ are
induced by the product cycles $[0,1]\times \eta$, thought of as
supported  the complement of two-handle additions).  

Inside the composite cobordism $W=B_2\cup_{Y_2} A_{1}$, one can find a sphere
with self-intersection number $S$ equal to zero, which represents a
non-trivial homology class in $H_2(W,\partial W;\Z)$. Specifically,
suppose that $A_{1}$ is built from $Y_1$ by attaching a two-handle along
$K_1$ (with some framing) and $B_2$ is obtained by then attaching a
two-handle along $K_2$ (with some other framing), then the two-sphere
$S$ corresponds to $K_2$, and it is homologically non-trivial since
the homology class corresponding to $K_1$ intersects  it once
(c.f. Figure~\ref{fig:TwoHandleSwap}). 

It follows from the composition law for cobordisms, together with
Lemma~\ref{lemma:SOneSTwo} that $$ 0 =\Hto(W;\Gamma_{W,\nu})=
\Hto(B_2;\Gamma_{B_2,\nu})\circ \Hto(A_1;\Gamma_{A_1,\nu})=\Gto_2\circ \Fto_1.
$$
The composite $\Fto_1\circ \Gto_2$ vanishes  in the same way.

\begin{figure}
\mbox{\vbox{\includegraphics{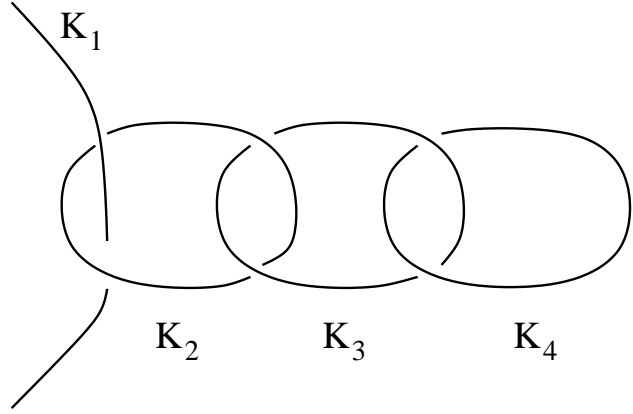}}}
\caption{\label{fig:TwoHandleSwap}
{\bf{Handle decomposition for Proposition~\ref{prop:SOneSTwo}.}}  Consider the link pictured above, where
here $K_1$ is thought of as any initial framed knot in the 
three-manifold $Y_1$, and $K_2$, $K_3$, and $K_4$ are unknots
with the property that $K_i$ links $K_{i-1}$ and $K_{i+1}$
geometrically once (for $i=2,3$). The three cobordisms 
$A_1\colon Y_1\to Y_2$, 
$A_2\colon Y_2\to Y_3$,
and $A_3\colon Y_3\to Y_1$
which induce maps fitting into the long exact sequence are given as
follows. $A_1$ is specified by the framed knot $K_1$, $A_2$ is specified
by the framed knot $K_2$ (thought of as a knot in $Y_2$) with
framing $-1$, while $A_3$ is specified by the framed knot $K_3$ with
framing $-1$.  Three cobordisms $B_2\colon Y_2\to Y_1$,
$B_1\colon Y_1\to Y_4$, $B_4\colon Y_4\to Y_2$ are specified as follows.
$B_2$ is specified by $K_2$
(thought of as a knot inside $Y_2$) with framing $0$,
while $B_1$ is specified by $K_3$ with framing $-1$, 
and $B_4$ is specified by $K_4$ with  framing $-1$. 
In particular, the cobordism
$B_2\cup_{Y_2} A_1\colon Y_1\to Y_1$ 
is specified by the link $K_1\cup K_2$, where here 
$K_2$ is given framing
$0$. In this composite cobordism, $K_2$
corresponds to a homologically non-trivial sphere with self-intersection
number zero.}
\end{figure}

\end{proof}

\begin{prop}
\label{prop:Rationals}
Let $K$ be a knot in $S^3$, and fix a cycle $\cycle\in S^3-K$ whose
homology class generates $H_1(S^3-K;\R)$. Let $r_0,r_1\in \Q\cup
\{\infty\}$ with $r_0, r_1$ non-negative. 
Suppose moreover that if we write $r_0$ and $r_1$ as fractions in
their lowest terms $r_i = p_i/q_i$ (where here all $p_i$, $q_i$ are
non-negative integers), then $p_0 q_1-p_1 q_0 = 1$.  Then,
we have a short exact sequence of the form: 
\[
0\longrightarrow \Hto_{\bullet}(S^3_{r_1}(K);\Gamma_\cycle)
\longrightarrow
\Hto_{\bullet}(S^3_{r_2}(K);\Gamma_\cycle) 
\longrightarrow \Hto_{\bullet}(S^3_{r_0}(K);\Gamma_\cycle)\longrightarrow 0,
\]
where here $r_2=\pfrac{p_0+p_1}{q_0+q_1}$.
\end{prop}

\begin{proof}
We prove the result by induction on $q_2=q_0+q_1$.

In the case where $q_0+q_1=1$, it follows that $q_0=0$ and $q_1=p_0=1$.
Now, Theorem~\ref{thm:TwoHandleSurgeries-Local}
gives us an exact sequence
\[
\longrightarrow  \Hto_{\bullet}(S^3_{p_1+1}(K);\Gamma_\eta) \longrightarrow 
\Hto_{\bullet}(S^3;\Gamma_\eta)\stackrel{\Hto(W;\Gamma_{W,\nu})}{\longrightarrow}\Hto_{\bullet}(S^3_{p_1}(K);\Gamma_\eta)\longrightarrow \cdots.
\]
We claim that the map $\Hto(W;\Gamma_\nu)\equiv 0$. This follows
from commutativity
of the diagram, 
$$
\begin{CD}
\Hred_\bullet(S^3;\Gamma_\eta)@>{\Hred_\bullet(W;\Gamma_{W,\eta})}>>
\Hred_\bullet(S^3_{p_1}(K);\Gamma_\eta) \\
@Vi_{*}VV @V{i'_{*}}VV \\
\Hto_{\bullet}(S^3;\Gamma_\eta))@>{\Hto(W)}>> \Hto_{\bullet}(S^3_{p_1}(K);\Gamma_\eta) 
\end{CD}
$$ bearing in mind that $i_*\colon \Hred_\bullet(S^3;\Gamma_\eta)\longrightarrow
\Hto_{\bullet}(S^3;\Gamma_\eta)$ is surjective, 
together with the fact that $\Hred_\bullet(W;\Gamma_{W,\eta})\equiv 0$, which
is analyzed in two cases. In the case where $p_1>0$, $\Hred(W)\equiv
0$ since $b_2^+(W)=1$, in view of
Proposition~\ref{prop:Reducibles-II}; while in the case where $p_1=0$,
it follows from the fact that $\Hred_\bullet(S^3_0(K);\Gamma_\eta)=0$, in view
of Lemma~\ref{lemma:VanishTwist}. (Note that the present case of
Lemma~\ref{lemma:VanishTwist} follows at once from
Proposition~\ref{prop:Reducibles}, together with the surgery long
exact sequence on the level of $\Hred$ with local coefficients.)

For the inductive step, Theorem~\ref{thm:TwoHandleSurgeries-Local}
gives a long exact sequence
$$
		\longrightarrow
                \Hto_{\bullet}(S^3_{r_2};\Gamma_\eta)\longrightarrow
                \Hto_{\bullet}(S^3_{r_0};\Gamma_\eta) 
		\stackrel{F}{\longrightarrow}\Hto_{\bullet}(S^3_{r_1};\Gamma_\eta)\longrightarrow \cdots
$$
Let $r_3 =\pfrac{p_0-p_1}{q_0-q_1}$. 
If $q_0>q_1$, then by induction on the denominator, we have a short
exact sequence
\[
0\longrightarrow \Hto_{\bullet}(S^3_{r_1};\Gamma_\eta)
\stackrel{G}{\longrightarrow}
\Hto_{\bullet}(S^3_{r_0};\Gamma_\eta)
\longrightarrow \Hto_{\bullet}(S^3_{r_3};\Gamma_\eta)\longrightarrow  0.
\]
By Proposition~\ref{prop:SOneSTwo}, it follows that $G\circ F=0$. Since
$G$ is injective, it follows that $F$ is the trivial map.

In the case where $q_0<q_1$, by induction on the denominator, we have
$$
\begin{CD}
0\longrightarrow \Hto_{\bullet}(S^3_{r_3}(K);\Gamma_\eta) \longrightarrow \Hto_{\bullet}(S^3_{r_1}(K);\Gamma_\eta) \stackrel{G}{\longrightarrow}\Hto_{\bullet}(S^3_{r_0}(K);\Gamma_\eta)\longrightarrow  0.
\end{CD}
$$
Now, since $F\circ G=0$ and $G$ is surjective, it follows that $F\equiv 0$.

In the final case, where $q_0=q_1$, it follows that $q_0=q_1=1$ and
that $p_0=p_1+1$. Exactness now follows from the short exact sequence
\[
0\longrightarrow
\Hto_{\bullet}(S^3_{p_1}(K);\Gamma_\nu)\longrightarrow
\Hto_{\bullet}(S^3_{p_1+1}(K);\Gamma_\nu)
\longrightarrow \Hto_{\bullet}(S^3;\Gamma_\nu)
\longrightarrow 0
\]
which was established earlier.
\end{proof}

\begin{prop}
\label{prop:LSpaceRat}
Let $K$ be a knot in $S^3$ and suppose that 
$$j\colon \Hto_{\bullet}(S^3_r(K))\longrightarrow \Hfrom_\bullet(S^3_r(K))$$
is trivial for some non-integral, rational $r>0$. Let $p$ be the smallest integer
greater than or equal to $r$, then
$$j\colon \Hto_\bullet(S^3_p(K))\longrightarrow \Hfrom_\bullet(S^3_p(K))$$
is trivial, as well.
\end{prop}

\begin{proof}
Note that for $r>0$, if $\eta$ is any real cycle
in the rational homology sphere $S^3_r(K)$, then 
the map 
$$j_\cycle\colon
\Hto_\bullet(S^3_r(K);\Gamma_\cycle)\longrightarrow \Hfrom_\bullet(S^3_r(K);\Gamma_\cycle)$$ is
non-trivial if and only if the corresponding map $$j\colon
\Hto_\bullet(S^3_r(K))\longrightarrow \Hfrom_\bullet(S^3_r(K))$$ is; 
in fact since $\eta$ is null-homologous, 
we have identifications 
\begin{eqnarray*}
\Hto_\bullet(S^3_r(K);\Gamma_\eta)\cong
\Hto_\bullet(S^3_r(K))\otimes\Kfield
&{\text{and}}&
\Hfrom_\bullet(S^3_r(K);\Gamma_\eta)\cong
\Hfrom_\bullet(S^3_r(K))\otimes\Kfield,
\end{eqnarray*}
under which the map $j\otimes \Id_{\Kfield}$ is identified with
$j_\eta$.

Write $r=p/q$ in its lowest terms.  Since $p$ and $q$ are relatively
prime, we can find a pair of integers $a$ and $b$ with the property
that $a q - b p = \pm 1$. Since $p>0$ and $q>1$, it follows that $a$
and $b$ must have the same sign, or $a=0$.  Without loss of
generality, we can  assume that $a$ and $b$ are both
non-negative. By simultaneously subtracting multiples of $p$ off from
$a$ and multiples of $q$ off from $b$, we can arrange for
$0\leq a<p$ and $0<b<q$. If $a q - b p = +1$, let $r_0 =
\nofrac{a}{b}$, $r_1=\pfrac{p-a}{q-b}$ and $r_2 =r= \nofrac{p}{q}$, while
if $a q - b p = -1$, we let $r_0 = \pfrac{p-a}{q-b}$ and
$r_1=\nofrac{a}{b}$. In both cases, the short exact sequence from
Proposition~\ref{prop:Rationals} holds, giving us the following
diagram, where the rows and columns are exact
\[
\begin{small}
\begin{CD}
\Hred_\bullet(S^3_{r})@>{\Fred}>> \Hred_\bullet(S^3_{r_0}) @>>> 0 \\
@V{i^r_*}VV @V{i^{r_0}_{*}}VV   \\
\Hto_\bullet(S^3_{r}) @>{\Fto}>>\Hto_\bullet(S^3_{r_0}) @>>> 0  \\
@VVV \\
0 \\
\end{CD}
\end{small}
\]
It follows at once that $i^{r_0}_{*}$ is surjective as well. Note
that the denominator of $r_0$ is smaller than that of $r$, 
and there are no integers between $r$ and $r_0$; 
hence by induction on this denominator, the result follows
from the long exact sequence which connects $p_{*}$, $i_{*}$, and
$j_{*}$.
\end{proof}

Let $K\subset S^3$ be a knot and $r>0$ be a rational number.  
We can construct a map
$${\pref_r}\colon \SpinC(S^3_r(K))\longrightarrow
\SpinC(S^3_r(U))$$
as follows.
Consider the Hirzebruch-Jung continued fractions expansion of $r$
\begin{equation}
\label{eq:HJContFrac}
r= [a_1,..,a_n]= a_1 - \bigfrac{1}{a_2 - \bigfrac{1}{\ddots -\bigfrac{1}{a_{n}}}},
\end{equation}
where $a_1\geq 1$ and $a_i\geq 2$ for $i>1$. 
Consider the four-manifold whose Kirby calculus
picture is given $K=K_1$ followed by a chain of unknots $K_2,\ldots,K_n$,
where $K_i$ links $K_{i-1}$ and $K_{i+1}$ once; and the framing of
$K_i$ is $-a_i$.  After deleting a ball, this gives a cobordism
$$W_r(K)\colon S^3_r(K) \longrightarrow S^3.$$ It is straightforward
to see that any $\SpinC$ structure $\spinct\in\SpinC(S^3_r(K))$ can be
extended to a $\SpinC$ structure over $W$. Let $\spinc$ denote such an
extension. Next, let $W_r(U)$ denote the corresponding cobordism for the
unknot $$W_r(U)\colon L(p,q)
\longrightarrow S^3.$$ 
By the construction of these cobordisms, there is a distinguished identification
$$\tau\colon H^2(W_r(U);\Z) \stackrel{\cong}{\longrightarrow} H^2(W_r(K);\Z).$$
Let, $\spinc'\in\SpinC(W_r(U))$ denote the $\SpinC$ structure
structure with $\tau(c_1(\spinc'))=c_1(\spinc)$.
It is straightforward to see that the correspondence which sends
$\spinct$ to the restriction of $\spinc'$ to $S^3_r(U)\subset \partial W_r(U)$
induces a well-defined map ${\pref_r}$ as stated in the beginning of the paragraph.

\begin{prop}
\label{prop:Inequality}
Let $K$ be a knot in $S^3$. Then, for all $\spinct\in\SpinC(S^3_r(K))$, we have
that 
$$\max_{\{\spinc\in\SpinC(W_r(U))\big| \spinc|_{S^3_r(U)}=\pref_r(\spinct)\}}
c_1(\spinc)^2+\rk H_2(W_r(U)) \leq -4 \Froy(\spinct),$$
with equality when $K=U$.
\end{prop}

\begin{proof}
Throughout this proof, we let $\alpha\in\Hred_{\xi_-}(S^3)$
be the non-zero class supported in the in the summand corresponding
to the two-plane field $\xi_-\in J(S^3)$ (recall that  $h(\xi_-)=0$).
Note that $i_*(\alpha)\in\Hto_{\bullet}(S^3)$ is non-trivial.
Note that $W_r(K)\colon S^3_r(K) \longrightarrow S^3$ is a cobordism
between rational homology three-spheres with $b_2^+(W_r(K))=0$.
Since for each $\spinc\in\SpinC(W_r(K))$,
the induced map $\Hred(W_r(K),\spinc)$ is an isomorphism
(c.f. Proposition~\ref{prop:Reducibles}), so we get a diagram
\[
\begin{CD}
\Hred_*(S^3_r(K),\spinct) @>{\Hred(W_r(K),\spinc)}>{\cong}>\Hred_*(S^3) \\
@V{i_*}VV	@VV{i_*}V \\
\Hto_*(S^3_r(K),\spinct) @>{\Hto(W_r(K),\spinc)}>>\Hto_*(S^3), \\
\end{CD}
\]
where here $\spinct=\spinc|_{S^3_r(K)}$.
It follows that the map $i_*\circ \Hto(W_r(K),\spinc)$ is surjective,
and hence there is a two-plane-field $j\in J(\spinct)$ and
an element $\beta\in \Hred_{j}(S^3_r(K))$ with the property that
$\Hred(W_r(K),\spinc)(\beta)=\alpha$.

Now, by the dimension formula $$-4h(j) =
{c_1(\spinc)^2-2\chi(W_r(K))-3\sigma(W_r(K))}={c_1(\spinc)^2+\rk
H_2(W_r(K))}.$$ But since $i_*(\beta)\in \Hto_j(S^3_r(K))$ is a
non-zero-homogeneous element in the image of $i_*$, it follows that
$h(j)\geq \Froy(\spinct)$. Putting these together, we have  shown that
\begin{equation}
	\label{eq:FroyshovInequality}
	c_1(\spinc)^2+\rk H_2(W_r(K)) \leq - 4 \Froy(\spinct),
\end{equation}
for any $\spinc$ which extends $\spinct$ over $W_r(K)$. Note that
the left-hand-side of this equation depends only on the homological
properties of $W_r(K)$, and hence can be replaced by $W_r(U)$ as
in the statement of the proposition.

It remains to show that for given $\spinct\in\SpinC(S^3_r(U))$, there
is a $\SpinC$ structure for which equality holds in
Equation~\eqref{eq:FroyshovInequality}.

To this end, we claim that there is a cobordism $$V_r(U)\colon
S^3\longrightarrow S^3_r(U),$$ with the property that
$X=W_r(U)\cup_{S^3_r(U)} V_r(U)\colon S^3 \longrightarrow S^3$ is a
negative-definite four-manifold (indeed, it is obtained from the
cylinder $[0,1]\times S^3$ by a sequence of blow-ups). Moreover, each
$\SpinC$ structure $\spinct\in\SpinC(S^3_r(U))$ can be extended to a
$\SpinC$ structure $\fu\in\SpinC(X)$ with the property that
$c_1(\fu)^2+\rk H_2(X)=0$ (i.e. so that its square is maximal).
Concretely, $V_r(U)$ is constructed from a plumbing of spheres
with multiplicities $\{b_1,\ldots,b_m\}$ chosen so that 
$$ 1 = [a_1,\ldots,a_n,1,b_1,\ldots,b_m].$$
The property of $\SpinC$ structures with
minimal square can be proved by induction on the size of the expansion.

Given this fact, note that the composite $$i_*'\circ
\Hred(V_r(U),\fu|_{V_r(U)})\colon \Hred_{\bullet}(S^3) \longrightarrow
\Hto_{\bullet}(S^3_r(U))$$ is
once again surjective for all $\fu\in\SpinC(X)$
(as it is the composite of an isomorphism with a surjection).
It follows that if there were no 
no extension of $\spinct\in \SpinC(S^3_r(U))$ to $W_r(U)$ for
which equality holds in Equation~\eqref{eq:FroyshovInequality},
then the composite map 
$\Hto(X,\fu)=\Hto(W_r(U),\fu|_{W_r(U)})\circ \Hto(V_r(U),\fu_{V_r(U)})$
would have kernel for any $\fu\in\SpinC(U)$ with $\fu|_{S^3_r(U)}=\spinct$.
But for any choice of $\fu\in\SpinC(X)$ with
$c_1(\fu|_{W_r(U)})^2+\rk H_2(W_r(U);\Z)=0$, this map is an isomorphism.
\end{proof}

The following result reduces Theorem~\ref{thm:RP3} in the case where
$r$ is non-integral to the integral case:

\begin{theorem}
Suppose that there is an orientation-preserving diffeomorphism
$S^3_r(K)\cong S^3_r(U)$ for some non-integral $r>0$. Then $K$ is
$p$-standard, where $p$ is the smallest integer greater than $r$.
\end{theorem}

\begin{proof}
Since 
$S^3_r(K)$ is orientation-preserving diffeomorphic to
$S^3_r(U)$, it follows that 
$$\sum_{\spinct\in \SpinC(S^3_{r}(K))} \Froy(S^3_r(K),\spinct) = 
\sum_{\spinct\in \SpinC(S^3_{r}(U))} \Froy(S^3_r(U),\spinct).$$
Since $\pref_r$ is a bijection, it follows from this equation
together with 
Proposition~\ref{prop:Inequality} that in fact
$\Froy(S^3_r(K),\spinct)=\Froy(S^3_r(U),\pref_r(\spinct))$ for all
$\spinct\in\SpinC(S^3_r(K))$.  Since $j$ is trivial on
$\Hto_\bullet(S^3_r(U))=\Hto_\bullet(S^3_r(K))$, it follows from
Proposition~\ref{prop:LSpaceRat} that $j$ is trivial on
$\Hto_\bullet(S^3_p(K))$ as well. 

We claim that
$\Froy(S^3_p(K),\spinct)=\Froy(S^3_p(U),\pref_p(\spinct))$ for all
$\spinct\in\SpinC(S^3_p(U))$. To see this, note that by 
construction, we can decompose 
$$W_r(K) = V\cup_{S^3_p(K)} W_p(K),$$
where $V\colon S^3_r(K) \longrightarrow S^3_p(K)$ is
obtained from the $n-1$ two-handle additions (specified by $K_2,\ldots,K_n$).
We claim that
any ${\mathfrak u}\in\SpinC(S^3_p(K))$
admits an extension $\spinc$ over all of $W_r(K)$, so that the induced map
$$\Hto(W_r(K),\spinc)\colon \Hto_\bullet(S^3_r(K),\spinc|W_r(K)) 
\longrightarrow \Hto_\bullet(S^3)$$ is
an isomorphism. Now, corresponding to the decomposition, we can write
$$\Hto(W_r(K),\spinc)=\Hto(W_p(K),\spinc|_{W_p(K)})\circ
\Hto(V,\spinc|_V).$$
Since $b_2^+(V)=0$,
$\Hred(V,\spinc|_V)$ induces an isomorphism (c.f. Proposition~\ref{prop:Reducibles}), 
and $j_*$ 
is trivial for both $S^3_r(K)$ and $S^3_p(K)$, it follows easily that
$\Hto(V,\spinc|_V)$ is surjective. It follows  that
$$\Hto(W_p,\spinc|_{W_p(K)})\colon \Hto_\bullet(S^3_p(K),\spinct)\longrightarrow \Hto_\bullet(S^3)$$
is an isomorphism, and hence that 
$\Froy(S^3_p(K),\spinct)=\Froy(S^3_p(U),\pref_p(\spinct))$.
From this, it follows readily that $K$ is $p$-standard.
\end{proof}

\section{Further applications for lens space surgeries}
\label{sec:FurtherApps}

In this section, we use the surgery long exact sequence, along with
some earlier results, to study the more general problem of lens space
surgeries.

Recall the following result of Meng and Taubes~\cite{MengTaubes}
(reformulated in the context of monopole Floer homology):

\begin{theorem}
\label{thm:MengTaubes}
Let $K$ be a knot in $S^3$, and write its symmetrized
Alexander polynomial as
$$\Delta_K(T)=a_0 + \sum_{i} a_i(T^i+T^{-i}).$$
We fix a  generator $h\in H_2(S^3_0(K);\Z)$, and let 
$\Hto_\bullet(S^3_0(K),i)$ denote the Floer homology of $S^3_0(K)$
with local coefficients determined  by any 
cycle $\nu$ which generates $H_1(S^3_0(K);\R)$),
evaluated in the $\SpinC$ stucture $\spinc$ with
$\langle c_1(\spinc),[h]\rangle = 2i$. 
Then
$$\chi_{\Kfield} (\Hto(S^3_0(K),i))=-\sum_{j=1}^\infty j a_{|i|+j}.$$
(Where here the left-hand-side is the Euler characteristic over
the field  $\Kfield$ of fractions of the group-ring $\Field[\R]$, and
the sign is determined by the canonical mod two grading on Floer homology
described in Subsection~\ref{subsec:CanonModTwo}.)
\end{theorem}

Combining this with the results of this paper, we obtain the following
necessary criterion for a $S^{3}_{p}(K)$ to be a lens space:

\begin{theorem}
\label{thm:PreciseGenusBounds}
Let $K$ be a knot in $S^3$ with the property that some integer surgery
on $K$ gives a lens space,
then the Seifert genus of $K$ coincides with the degree of the symmetrized
Alexander polynomial of $K$.
\end{theorem}

\begin{proof}
If the genus of $K$ is bigger than the degree of the Alexander
polynomial, then the according to the Meng-Taubes theorem, we know
that $\chi(\HFto(S^3_0(K),g-1))=0$, while by the non-vanishing result,
Corollary~\ref{cor:NonVanishingTwist}, $\HFto(S^3_0(K),g-1)\neq 0$,
and hence $\HFto_{\odd}(S^3_0(K))\neq 0$. Thus, by
Theorem~\ref{thm:TwoHandleSurgeries-Local}, we have the exact sequence
$$\longrightarrow \HFto_{\odd}(S^3;\Gamma_\nu) \longrightarrow 
\HFto_{\odd}(S^3_0(K);\Gamma_\nu) \longrightarrow 
\HFto_{\odd}(S^3_1(K);\Gamma_\nu)\longrightarrow 
$$
Since $\HFto_{\odd}(S^3)=0$, it follows that $\HFto_{\odd}(S^3_1(K))\neq 0$
(whether or not we use the local coefficient system $\Gamma_\nu$).

Indeed, if $\HFto_{\odd}(S^3_p(K))\neq 0$, it follows easily that
$\HFto_{\odd}(S^3_{p+1}(K))\neq 0$, since in this case, 
Theorem~\ref{thm:TwoHandleSurgeries-A}
sequence takes the form
$$
\longrightarrow \Hto_\bullet(S^3) \longrightarrow  \Hto_\bullet(S^3_p(K)) \longrightarrow  
\Hto_\bullet(S^3_{p+1}(K))\longrightarrow ,
$$
where here $F_{p}$ and $F_{p+1}$ preserve the absolute $\Zmod{2}$ grading
(c.f. Section~\ref{subsec:CanonModTwo}).
In particular, this makes it impossible for $S^3_n(K)$ to be a lens space
(for integral $n>0$).
\end{proof}

In~\cite{AbsGraded}, it is shown that if $K$ is any knot with the
property that $S^3_p(K)=L(p,q)$, then the Alexander polynomial of $K$
is uniquely determined up to a finite indeterminacy. Indeed, this algorithm
is concrete: its input is the Milnor torsion for $L(p,q)$ and $L(p,1)$, 
and the finite indeterminacy depends on the homology class of the induced
knot in $L(p,q)$. 

Consider the rational numbers $d(-L(p,q),i)$ associated to a lens space
$L(p,q)$
(recall
that we have fixed here the orientation convention that
$L(p,q)=S^3_{p/q}(U)$), and an element $i\in\Zmod{p}$, determined by the recursive
formula
\begin{eqnarray*}
d(-L(1,1),0)&=&0\\
d(-L(p,q),i)&=&\left(\frac{pq-(2i+1-p-q)^2}{4pq}\right)-d(-L(q,r),j), 
\end{eqnarray*}
where $r$ and $j$ are the reductions modulo $q$ of $p$ and $i$
respectively. 
(Note that these numbers turn out to agree with the Fr{\o}yshov invariants
of the lens space $-L(p,q)$, under a particular identification
$\SpinC(L(p,q))\cong \Zmod{p}$.)

The following can be found in Corollary~\ref{AbsGraded:cor:AlexLens} 
of~\cite{AbsGraded}:

\begin{theorem}
\label{thm:AlexLens}
The lens space $L(p,q)$ is obtained as surgery on a knot $K\subset S^3$
only if there is a one-to-one correspondence 
$$\sigma\colon \Zmod{p}\longrightarrow \SpinC(L(p,q))$$
with the following symmetries:
\begin{itemize}
\item $\sigma(-[i])={\overline{\sigma([i])}}$
\item there is an isomorphism $\phi\colon \Zmod{p}
\longrightarrow \Zmod{p}$
with the property that 
$$\sigma([i])-\sigma([j])=\phi([i-j]),$$
\end{itemize}
with the following properties.  For $i\in\Z$, let $[i]$ denote its reduction 
modulo $p$,
and define
$$t_i=
\left\{\begin{array}{ll}
-d(L(p,q),\sigma[i]) + d(L(p,1),[i]) &{\text{if $2|i|\leq p$}} \\ \\
0 & {\text{otherwise,}}
\end{array}
\right. 
$$ then the Laurent polynomial
$$L_\sigma(T)=1 + \sum_{i} \left(\frac{t_{i-1}}{2}-t_i+ \frac{t_{i+1}}{2}\right)
T^i =\sum_i a_i \cm T^i$$ has integral coefficients,
and all the $t_i\leq 0$. Indeed, if $S^3_p(K)\cong L(p,q)$,
then its Alexander polynomial has the form
$L_\sigma(T)$ for some choice of $\sigma$ as above.
\end{theorem}
 
By combining Theorem~\ref{thm:PreciseGenusBounds}, results
of~\cite{AbsGraded}, and work of Goda and
Teragaito~\cite{GodaTeragaito}, we obtain the following:

\begin{cor}
\label{cor:Trefoil}
If $K$ is a knot with the property that for some $p\in\Z$,
$S^3_{p}(K)$ is a lens space and $|p|<9$, then $K$ is either the
unknot or the trefoil.
\end{cor}

\begin{proof}
In view of Theorem~\ref{thm:AlexLens}, it is now an experiment in
numerology to see that if $S^3_p(K)$ is a lens space with $|p|<9$,
then Alexander polynomial of $K$ is either trivial or $T-1+T^{-1}$
(see the list at the end of Section~\ref{AbsGraded:sec:Examples}
of~\cite{AbsGraded}).  In view of
Theorem~\ref{thm:PreciseGenusBounds}, it follows that the genus of $K$
is zero or one. Combining this with a theorem of Goda and
Teragaito~\cite{GodaTeragaito}, according to which the only genus one
knot which admits lens space surgeries is the trefoil, the corollary
is complete.
\end{proof}

As another application, we obtain the following bound on the Seifert
genus $g$ of $K$ in terms of the order of the lens space.

\begin{cor}
\label{cor:GenusBound}
Let $K$ be a knot in $S^3$ with the property that 
for some integer $p$, $S^3_p(K)$ is a lens space, 
then $2g-1\leq p$. 
\end{cor}

\begin{proof}
The bound $2d-1\leq p$ where $d$ is the degree of the Alexander
polynomial of $K$ follows immediately from
the algorithm described in Theorem~\ref{thm:AlexLens}; the 
rest follows from Theorem~\ref{thm:PreciseGenusBounds}.
\end{proof}

The bound on the Seifert genus stated above is still fairly coarse,
and can usually be improved for fixed $p$ and $q$ using
Theorem~\ref{thm:PreciseGenusBounds}, combined with the algorithm
for determining the Alexander polynomial of $K$
given in Theorem~\ref{thm:AlexLens}.

It is interesting to compare Corollary~\ref{cor:GenusBound} with a
conjecture of Goda and Teragaito for hyperbolic knots which admit lens
space surgeries. They conjecture that for such a knot, the order $p$
of the fundamental group is related with the Seifert genus $g$ by the
inequalities $2g+8\leq p\leq 4g-1$.  Indeed they prove (Theorem~1.1
of~\cite{GodaTeragaito}) that if $K$ is a hyperbolic knot in $S^3$,
and if $S^3_p(K)$ is a lens space, then $|p|\leq 12g-7$.  They
restrict to the hyperbolic case, since the case of non-hyperbolic
knots yielding lens space surgeries is completely understood,
c.f.~\cite{BleilerLitherland},
\cite{Wang}, \cite{Wu}. The only such knots with lens space
surgeries are torus knots, and the $(2,2pq\pm 1)$-cable of a $(p,q)$
torus knot (in which case the resulting lens space
is $L(4pq\pm 1, 4q^2)$). 

The condition that $t_i\leq 0$ from Theorem~\ref{thm:AlexLens} has an
improvement, using the Floer homology for knots
(see~\cite{HolDiskKnots} and~\cite{RasmussenThesis}). Specifically, in
Corollary~\ref{NoteLens:cor:StructAlex} of~\cite{NoteLens}, it is
shown that if $K$ is a knot on which some integral surgery is a lens
space, then all the non-zero coefficients of its Alexander polynomial
are $\pm 1$, and they alternate in sign.

Combining all this information, we can give stronger constraints on
the lens spaces which can be obtained by surgeries on knots with a
fixed Seifert genus. As an illustration, we have the following:

\begin{cor}
\label{cor:GenusTwo}
The only lens spaces which can be obtained by positive integer surgery on a
knot in $S^3$ with Seifert genus $2$ are orientation-preserving
diffeomorphic to $L(9,7)$ and $L(11,4)$. 
\end{cor}

\begin{proof}
According to Theorem~\ref{thm:PreciseGenusBounds}, combined with
Corollary~\ref{NoteLens:cor:StructAlex} of~\cite{NoteLens}
(which states that the non-zero coefficients of the Alexander polynomial
all have absolute value one and alternate in sign), we see
that if $K$ has genus two and some integral surgery on it gives a lens
space, then $\Delta_K$ is either $T^{-2}-1+T^2$ or
$T^{-2}-T^{-1}+1-T+T^2$. In the first case, $K$ is neither a torus
knot nor the $(2,2pq\pm 1)$-cable of a $(p,q)$ torus knot, and hence
it must be hyperbolic (c.f.~\cite{BleilerLitherland},
\cite{Wang}, \cite{Wu}). According to Goda and Teragaito's bound,
$|p|\leq 17$. But this is now ruled out by
Corollary~\ref{AbsGraded:cor:AlexLens} of~\cite{AbsGraded}. In the
second case, we can rule out the possibility that $K$ is hyperbolic
and $p\neq 9, 11$ in the same manner. In these remaining cases, the
algorithm of Theorem~\ref{thm:AlexLens} forces
$S^3_p(K)$ to be one of the two listed possibilities.
\end{proof}

The above procedure is purely algorithmic, and can be repeated for
higher genera. For instance, if positive integral surgery on a genus
three knot gives a lens space, then that lens space is contained in
the list $$L(11,9),\hskip .5in L(13,10), \hskip .5in L(13,9),\hskip .5in
L(15,4)$$ (all of which are
realized by torus knots); in the genus four case, the list is
$$L(14,11),\hskip .5in
L(16,9),\hskip .5in L(17,13),\hskip .5in L(19,5)$$ (again, all of these are
realized by torus knots). In the genus five case, the list reads
$$L(18,13),\hskip .5in L(19,11),\hskip .5in L(21,16),\hskip .5in L(23,6),$$ where now the last two
examples are realized by a torus knot, and the first two are realized
by the $(-2,3, 7)$ pretzel knot (c.f.~\cite{FSLensSurgeries}).

In a different direction, we can combine
Theorem~\ref{thm:PreciseGenusBounds} with properties of the Heegaard
Floer homology for knots (see~\cite{HolDiskKnots}
or~\cite{RasmussenThesis}) to obtain the following result on the
four-ball genera of knots admitting lens space surgeries, compare
also~\cite{KMMilnor}:

\begin{cor}
Let $K$ be a knot in $S^3$ with the property that $S^3_p(K)=L(p,q)$.
Then, the Seifert genus, the four-ball genus, and the
degree of the Alexander polynomial all coincide.
\end{cor}

\begin{proof}
According to~\cite{NoteLens}
(c.f. Corollary~\ref{NoteLens:cor:CalcTau} in that reference), the
four-ball genus of $K$ is bounded below by the degree of the Alexander
polynomial of $K$. The equality of the three quantities follows from
the fact that the four-ball genus is less than or equal to the Seifert
genus of $K$, together with Theorem~\ref{thm:PreciseGenusBounds}.
\end{proof}

\subsection{Seifert fibered surgeries}
\label{subsec:SurgSeif}

We give an application of the long exact sequence to the question of 
when a knot in $S^3$ admits a Seifert fibered surgery.
To state the strongest form, it is useful to pin down orientations.

Let $Y$ be a Seifert fibered space with $b_1(Y)=0$ or $1$. Such a
manifold can be realized as the boundary of a four-manifold
$W(\Gamma)$ obtained by plumbing two-spheres according to a weighted
tree $\Gamma$. Here, the weights are thought of as a map $m$ from the
set of vertices of $\Gamma$ to $\Z$.

\begin{defn}
\label{def:PositiveSeifert}
Let $Y$ be an oriented Seifert fibered three-manifold with $b_1(Y)=0$
or $1$. We say that $Y$ has a positive Seifert fibered orientation if
it can be presented as the oriented boundary of a plumbing of spheres
$W(\Gamma)$ along a weighted tree $\Gamma$ so that with
$b^-(W(\Gamma))=0$.  If $Y$ does not have a positive Seifert fibered
orientation, then $-Y$ does, and we say that $Y$ is {\em negatively
oriented}.
\end{defn}

Moreover, either orientation on any lens space is a positive Seifert
orientation; similarly, either orientation on a Seifert fibered space
with $b_1(Y)=1$ is a positive Seifert orientation.  Finally, if $Y$ is
the quotient of a circle bundle $\pi\colon N\longrightarrow \Sigma$
over a Riemann surface by a finite group of orientation-preserving
automorphisms $G$, and if $N$ is oriented as a circle bundle with
positive degree, then the induced orientation on $Y$ is a positive
Seifert orientation.

The basic property of the monopole Floer homology of Seifert fibered
spaces we will use is the following result, which follows quickly
from~\cite{MOY}. Or, alternatively, using
Theorem~\ref{thm:TwoHandleSurgeries-A}, one can adapt the proof of the
corresponding result for Heegaard Floer homology
(c.f. Corollary~\ref{Seifert:cor:EvenDegrees} of~\cite{SomePlumb}),
to the context of Seiberg-Witten monopoles.

\begin{theorem}
\label{thm:SeifertFibered}
If $Y$ is a positively oriented Seifert fibered
rational homology three-sphere, then 
$\Hto_{\bullet}(Y)$ is supported entirely in even degrees.
\end{theorem}

Sometimes, we will consider Seifert fibered spaces with first Betti
number equal to one. One could adapt techniques of~\cite{MOY} to this
situation, as well, but it is quicker now to appeal to the surgery
long exact sequence. The relevant fact in this case is the following:

\begin{cor}
\label{cor:BOneSeifert}
If $Y_0$ is a Seifert fibered space with $b_1(Y_0)=1$, and $\cycle$ is a
generator of $H_1(Y_0;\R)$, then
$\Hto_{\bullet}(Y_0;\Gamma_\cycle)$ is supported entirely in odd degrees.
\end{cor}

\begin{proof}
Express $Y_0$ as the boundary of a plumbing $W$ of spheres with
$b^-(W)=0$, and let $Y_1$ denote the new Seifert fibered space
obtained by increasing the multiplicity of the central node by one,
and let $Y_2$ denote the plumbing of spheres obtained by deleting the
central node. (The latter space, $Y_2$ is a connected sum of lens
spaces.)  We have that the Floer homology groups $Y_0$, $Y_1$, and
$Y_2$ fit into a long exact sequence as in
Theorem~\ref{thm:TwoHandleSurgeries-Local}. In fact, since $Y_2$ can
be given a positive scalar curvature metric, $i_*\colon
\Hred(Y_2;\Gamma_\eta)\longrightarrow \Hto_{\bullet}(Y_2;\Gamma_\eta)$ is surjective, and hence, since
$\Hred(Y_0;\Gamma_\cycle)=0$ (c.f. Lemma~\ref{lemma:VanishTwist}),
it follows that the map from
$\Hto_{\bullet}(Y_2;\Gamma_\eta)\longrightarrow
\Hto_{\bullet}(Y_0;\Gamma_\eta)$ is
trivial. Thus, we get the short exact sequence
$$
0\longrightarrow \Hto_{\bullet}(Y_0;\Gamma_\eta)
\stackrel{\Hto(W_0)}{\longrightarrow}\Hto_{\bullet}(Y_1;\Gamma_\eta)
\stackrel{\Hto(W_1)}{\longrightarrow}\Hto_{\bullet}(Y_2;\Gamma_\eta)\longrightarrow 0.
$$ Since $\Hto_{\bullet}(Y_1;\Gamma_\eta)$ is supported in even degrees and the map
$\Hto(W_0)$ reverses the canonical mod two grading
(c.f. Proposition~\ref{prop:CanonicalMod2}), the result follows.
\end{proof}

Another application of the surgery long exact sequence gives the following:

\begin{prop}
\label{prop:OddDegrees}
Let $K$ be a knot in $S^3$. Then, for all $r=p/q>0$, we have that
$$\rk \HFto_{\odd}(S^3_{p/q};\Gamma_\nu) = q\cm \rk \HFto_{\ev}(S^3_0(K);\Gamma_\nu).$$
\end{prop}

\begin{proof}
This follows immediately from Proposition~\ref{prop:Rationals}.
\end{proof}

We obtain the following direct generalization of
Theorem~\ref{thm:PreciseGenusBounds}:

\begin{theorem}
\label{thm:SeifSurgeryObstructions}
Let $K$ be a knot whose Alexander polynomial $\Delta_K(T)$ has degree
strictly less than its Seifert genus. Then, there is no rational
number $r\geq 0$ with the property that $S^3_r(K)$ is a 
positively oriented Seifert fibered space.
\end{theorem}

\begin{proof}
In view of Theorem~\ref{thm:MengTaubes}, the condition on $K$ ensures
that $\Hto_{\ev}(S^3_0(K);\Gamma_\eta)\neq 0$
(and also that $\Hto_{\odd}(S^3_0(K);\Gamma_\eta)\neq 0$, but we do not use this here).  
The case where $r=0$ now is ruled
out by the latter fact, together with
Corollary~\ref{cor:BOneSeifert}. For the case where $r>0$,
$\Hto_{\odd}(S^3_r(K);\Gamma_\eta)\neq 0$ in view of Proposition~\ref{prop:OddDegrees},
and hence Theorem~\ref{thm:SeifertFibered} shows that it is never a
positively oriented Seifert fibered space.
\end{proof}

It is a more subtle problem to detect whether $S^3_r(K)$ is
a negative Seifert fibered space for $r\geq 0$. We include the following:

\begin{theorem}
If $K$ is a knot whose Seifert genus $g$ is strictly greater than
the degree of its Alexander polynomial, and also $g>1$, then 
$S^3_{1/n}(K)$ is not a Seifert fibered space for any integer $n$.
\end{theorem}

\begin{proof}
As usual, by reflecting $K$ if necessary, it suffices to consider
the case where $n>0$. $S^3_{1/n}(K)$ is not a positively oriented
Seifert fibered space, according to Theorem~\ref{thm:SeifSurgeryObstructions}.
Thus, we are left with the case where $n>0$ and $Y$ is a negatively oriented
Seifert fibered space.

We claim that if $Y$ is a negatively oriented Seifert fibered space,
then the cokernel of $i_\bullet\colon \Hred_\bullet(Y)\longrightarrow 
\Hto_\bullet(Y)$
is supported entirely in odd degrees. This follows easily from
duality.

Now, in view of Corollary~\ref{cor:NonVanishing},
our hypothesis on the knot $K$ ensures that 
$\Hto_{\odd}(S^3_0(K))\neq 0$.
Indeed, letting $\xi$ be any non-trivial
element of $\Hto_{\odd}(S^3_0(K),g-1)\subset \Hto_{\odd}(S^3_0(K))$,
its conjugate $\overline \xi$ lies in the summand 
$\Hto_{\odd}(S^3_0(K),-g+1)$, and hence it is linearly independent
of $\xi$.
Consider the surgery
exact sequence, following Proposition~\ref{prop:Rationals}
$$
0\longrightarrow \Hto_{\bullet}(Y_0(K))
\stackrel{F}{\longrightarrow}\Hto_{\bullet}(Y_{1/n}(K))\longrightarrow
\Hto_{\bullet}(Y_{\frac{1}{n-1}})\longrightarrow 0
,$$
where here $F$ reverses the mod $2$ degree.
Since $F$ is injective, if $\eta=F(\xi)$, 
then $\overline \eta$ is linearly independent from $\eta$.
It follows that $\eta$ is not in the image of 
$i_*\colon \Hred(Y_{1/n}(K))\longrightarrow \Hto_{\bullet}(Y_{1/n}(K))$,
since conjugation acts trivially on $\Hred(Y_{1/n}(K))$. In view of the previous paragraph, $S^3_{1/n}(K)$
is not negatively Seifert fibered for $n>0$. 
\end{proof}

\section{Foliations}
\label{sec:Foliations}

The exact sequence, together with Theorem~\ref{thm:NonVanishing} can
be used to exhibit large classes of three-manifolds admitting no
(coorientable) taut foliations. (Recall that all foliations we
consider in this paper are coorientable, so we drop this modifier from
the statements of our results.)

\begin{defn}
A {\em monopole $L$-space}  is a 
rational homology three-sphere $Y$ for which
$$j_*\colon \HFto_\bullet(Y) \longrightarrow \HFfrom_\bullet(Y)$$ is
trivial. 
\end{defn}

\noindent
Examples include all lens spaces, and indeed all three-manifolds with
positive scalar curvature.  By Theorem~\ref{thm:NonVanishing}, a
monopole $L$-space admits no taut foliations.

\begin{prop}
\label{prop:LSpaces}
Let $M$ be a connected, oriented three-manifold with torus boundary,
equipped with three oriented, simple closed curves $\gamma_0$,
$\gamma_1$, and $\gamma_2$ as in
Theorem~\ref{thm:TwoHandleSurgeries-A}.  Suppose moreover that $Y_0$,
$Y_1$, and $Y_2$ are rational homology three-spheres, with the
property that $$|H_1(Y_2;\Z)|=|H_1(Y_0;\Z)|+|H_1(Y_1;\Z)|.$$ Suppose
also that $Y_0$ and $Y_1$ are monopole $L$-spaces. Then it follows
that $Y_2$ is a monopole $L$-space, too.
\end{prop}

\begin{proof}
It follows from the hypotheses that the cobordism 
$W_0\colon Y_0 \longrightarrow Y_1$
has $b_2^+(W_0)=1$, and hence
according to Section~\ref{subsec:Reducibles}
the map 
 $\Hred(W_0)\colon \Hred(Y_0)\longrightarrow \Hred(Y_1)$ is trivial.
Now, Theorem~\ref{thm:TwoHandleSurgeries-A}, and our hypothesis, gives the following
diagram (where all rows and columns are exact):
\[\small
   \begin{CD}
        @. @VV{0}V @VV{j_{*}}V @VV{0}V @. \\
    @>>> \Hfrom_{\bullet}(Y_1) @>\Hfrom(W_{1})>>
      \Hfrom_{\bullet}(Y_{2}) @>\Hfrom(W_{})>> \Hfrom_{\bullet}(Y_0) @ >>>
     \\
      @. @VV{p_{*}}V @VV{p_{*}}V @VV{p_{*}}V @.\\
      @>0>> \Hred_{\bullet}(Y_1) @>\Hred(W_{1})>>
      \Hred_{\bullet}(Y_{2}) @>\Hred(W_{})>> \Hred_{\bullet}(Y_0) @ >0>>
       \\
              @. @VV{i_{*}}V @VV{i_{*}}V @VV{i_{*}}V @. 
      \\
      @>>> \Hto_{\bullet}(Y_1) @>\Hto(W_{1})>>
      \Hto_{\bullet}(Y_{2}) @>\Hto(W_{})>> \Hto_{\bullet}(Y_0) @ >>>
       \\
        @. @VV{0}V @VV{j_{*}}V @VV{0}V @. 
   \end{CD}
\]
Now, a diagram-chase shows that $j_{*}=0$ as well.
\end{proof}

Using the above proposition, together with Theorem~\ref{thm:NonVanishing},
we can find large classes of three-manifolds which admit no taut foliations.
We list several here:

\begin{defn}
A {\em weighted graph} is a graph $G$ equipped with an integer-valued
function $m$ on its vertices. The {\em degree} of a vertex $v$,
written $d(v)$, is the number of edges which contain it.
\end{defn}

\begin{cor}
Let $G$ be a connected, weighted tree which satisfies the inequality
$m(v)\geq d(v)$ at each vertex $v$, and for which the inequality is strict
at at least one vertex. Let $Y(G)$ denote the three-manifold
which is the boundary of the sphere-plumbing associated to the graph $G$.
Then, $Y(G,m)$ admits no taut foliations.
\end{cor}

\begin{proof}
We prove that under the hypotheses on $(G,m)$, $Y(G,m)$ is a monopole
$L$-space by an induction on the number of vertices of $G$. If the
number of vertices is one, then $Y(G,m)$ is a lens space. Let $G$ be a
general tree satisfying the hypotheses, and fix a leaf $v$. We prove
the corollary by sub-induction on the weight $m(v)$. If $m(v)=1$, we can
form the graph $(G',m')$ where here $G'$ is obtained from $G$ by
deleting the vertex $G$, and $m'$ agrees with $m|G'$ except at the
neighbor $w$ of $v$, where $m'(w)=m(w)-1$.  It is easy to see that
$Y(G,m)\cong Y(G',m')$ and that $(G',m')$ also satisfies the
hypotheses, and hence it is an $L$-space by the inductive hypothesis
(on the number of leaves). Let $(G_0,m_0)$ be the weighted graph
obtained by deleting the vertex $v$ from $G$, and $(G_1,m_1)$ be the
weighted graph obtained from $G$ by decreasing the weight at $v$ by one.
It is easy to see that
$$|H_1(Y(G,m);\Z)|=|H_1(Y(G_0;\Z)|+|H_1(G_1;\Z)|.$$
The inductive step  now follows from
Proposition~\ref{prop:LSpaces}, together with
Theorem~\ref{thm:NonVanishing}. (For more details, see
the corresponding result in~\cite{SomePlumb}.)
\end{proof}

\begin{cor}
Let $L$ be a non-split, alternating link in $S^3$, and let $\Sigma(L)$ denote
the branched double-cover of $S^3$ along $L$. 
Then, $\Sigma(L)$ does not admit a taut foliation.
\end{cor}

\begin{figure}
\mbox{\vbox{\includegraphics{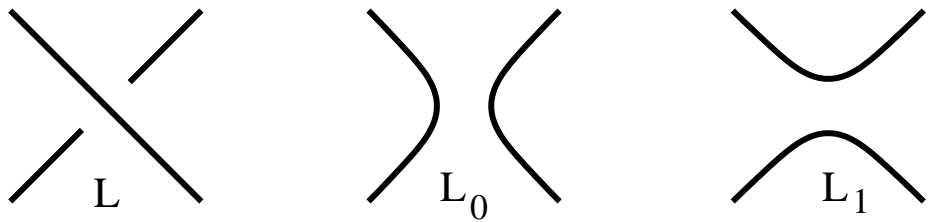}}}
\caption{\label{fig:LinkCrossings}
{\bf{Resolving link crossings.}}
Given a link with a crossing as labelled in $L$ above, we have two
resolutions $L_0$ and $L_1$, obtained by replacing the crossing by
the two simplifications pictured above.}
\end{figure}

\begin{proof}
For any link
diagram for $L$, choose a crossing, and let $L_0$ and $L_1$ denote the
two resolutions of $L$ at the crossing, as pictured in
Figure~\ref{fig:LinkCrossings}. It is easy to see that $\Sigma(L)$,
$\Sigma(L_0)$, and $\Sigma(L_1)$ fit into an exact triangle as in
Theorem~\ref{thm:TwoHandleSurgeries-A}. Furthermore, if we start with a connected,
reduced alternating projection for $L$, then both $L_0$ and $L_1$ are
connected, alternating projections and hence $\Sigma(L)$,
$\Sigma(L_0)$ and $\Sigma(L_1)$ are rational homology spheres. Indeed,
$$|H_1(\Sigma(L))|=|H_1(\Sigma(L_0))|+|H_1(\Sigma(L_1)|$$ by classical
knot theory (c.f.~\cite{Lickorish}). The result now follows from
Proposition~\ref{prop:LSpaces}, together with
Theorem~\ref{thm:NonVanishing}.
(For more details, see
the corresponding result in~\cite{BranchedDoubleCovers}.)
\end{proof}

\begin{cor}
\label{cor:LensSpaceSurgeries}
Let $K\subset S^3$ be a knot for which there is a positive rational
number $r$ with the property that $S^3_r(K)$ is a lens space or, more
generally, a monopole $L$-space. Then, for all rational numbers $s\geq
r$, $S^3_s(K)$ admits no taut foliation.
\end{cor}

\begin{proof}
By Proposition~\ref{prop:LSpaceRat}, if $p$ denotes the smallest
integer greater than $r$, then $S^3_p(K)$ is also a monopole
$L$-space. The result now follows from repeated applications of
Proposition~\ref{prop:LSpaces}.
\end{proof}

By a theorem of Thurston (c.f.~\cite{Thurston}, \cite{Thurston2}), if
$K$ is not a torus knot or a satellite knot, then $S^3_r(K)$ is
hyperbolic for all but finitely many $r$. 
One can now use the above corollary to construct infinitely many
hyperbolic three-manifolds with no taut foliation by, for example,
starting with a hyperbolic knot which admits some lens space surgery.
(A complete list of the known knots which admit lens space surgeries
can be found in~\cite{BergeUnpublished}, see also~\cite{Berge},
\cite{GabaiLens},
\cite{FSLensSurgeries}.) 
The first examples of hyperbolic three-manifolds which admit no taut foliation
were constructed in~\cite{RobertsShareshianStein}, see also~\cite{CalegariDunfield}

\begin{cor}
Fix an odd integer $n\geq 7$, and let $K$ be the $(-2,3,n)$ pretzel
knot.  For all $r\geq 2n+4$, then $S^3_r(K)$ admits no taut
foliations.
\end{cor}

\begin{proof}
When $n=7$, then $S^3_{18}(K)$ is a lens space. The result now follows
from Corollary~\ref{cor:LensSpaceSurgeries} (see also~\cite{Jun}). Indeed, for $n\geq 7$, it
is well-known (c.f.~\cite{BleilerHodgson}), $S^3_{2n+4}(K)$ is a
Seifert fibered space with Seifert invariants
$(-2,1/2,1/4,(n-8)/(n-6))$. Repeated applications of
Proposition~\ref{prop:LSpaces} can be used to show that this Seifert
fibered space is a monopole $L$-space, and hence again we can apply
Corollary~\ref{cor:LensSpaceSurgeries}.
\end{proof}

These bounds can be sharpened. For example, when $n=7$, $S^3_{17}(K)$
is a quotient of $S^3$ by a finite isometry group, and in particular,
it has positive scalar curvature. Thus, for all $r\geq 17$, $S^3_r(K)$
admits no taut foliation. In fact, with some extra work, one can
improve the bound in general to $r\geq n+2$. It is interesting to
compare this with results of Roberts~\cite{RobertsI},
\cite{RobertsII}, which constructs taut foliations on certain
surgeries on fibered knots. For example, when $K$ is the $(-2,3,7)$
pretzel knot, she shows that $S^3_r(K)$ admits a taut foliation for
all $r<1$.



In a similar vein, we have

\begin{prop}
For any three rational numbers $a, b, c \geq 1$, the three-manifold
$M(a,b,c)$
obtained by performing $a$, $b$, and $c$ surgery on the Borromean rings
carries no taut foliations.
\end{prop}

\begin{proof}
First, we prove the case where $a$, $b$, and $c$ are integers.  In the
basic case where $a=b=c$, the manifold $M(1,1,1)$ is the Poincar\'e
homology sphere, which admits a metric of positive scalar curvature,
and hence it is a monopole $L$-space. For the inductive step, we apply
Proposition~\ref{prop:LSpaces}, to see that the fact that $(p,1)\#
L(c,1)$ and $M(a,b,c)$ are monopole $L$-spaces implies that so is
$M(a+1,b,c)$. Repeated applications of the proposition also gives now
the result for all rational numbers $a$, $b$, and $c$ in the range.
\end{proof}

Note that this family includes the ``Weeks manifold'' $M(1,5/2,5)$,
which is known by other methods not to admit any taut foliations, 
see~\cite{CalegariDunfield}.

\bibliographystyle{plain}
\bibliography{biblio}

\end{document}